\documentclass{amsart}
\usepackage{amsmath,amssymb}
\usepackage{bbold}
\usepackage{leftidx}
\usepackage[all]{xy}

\newtheorem{thm}{Theorem}[section]
\newtheorem{cor}[thm]{Corollary}
\newtheorem{lem}[thm]{Lemma}
\newtheorem{prop}[thm]{Proposition}
\theoremstyle{definition}
\newtheorem{rem}[thm]{Remark}
\newtheorem{exa}[thm]{Example}

\newcommand{\co}{\colon}
\newcommand{\id}{\mathrm{id}}
\newcommand{\opp}{\mathrm{op}}
\newcommand{\cc}{\mathcal{C}}

\newcommand{\dd}{\mathcal{D}}
\newcommand{\ee}{\mathcal{E}}
\newcommand{\uu}{\mathcal{U}}
\newcommand{\lhs}{\mathcal{L}}
\newcommand{\rhs}{\mathcal{R}}
\newcommand{\biMod}[2]{\leftidx{_{#1}}{\mathrm{Mod}}{_{#2}}}
\newcommand{\bimod}[2]{\leftidx{_{#1}}{\mathrm{mod}}{_{#2}}}
\newcommand{\lMod}[1]{\leftidx{_{#1}}{\mathrm{Mod}}{}}
\newcommand{\lmod}[1]{\leftidx{_{#1}}{\mathrm{mod}}{}}
\newcommand{\biproj}[2]{\leftidx{_{#1}}{\mathrm{proj}}{_{#2}}}
\newcommand{\iso}{\stackrel{\sim}{\longrightarrow}}
\newcommand{\kk}{\Bbbk}
\newcommand{\kt}{$\Bbbk$\nobreakdash-\hspace{0pt}}
\newcommand{\trait}{\nobreakdash-\hspace{0pt}}
\newcommand{\Rt}{$\mathrm{R}$\nobreakdash-\hspace{0pt}}
\newcommand{\ti}{\mbox{-}\,}
\newcommand{\un}{\mathbb{1}}
\newcommand{\ad}{\mathrm{ad}}
\newcommand{\Ob}{\mathrm{Ob}}

\newcommand{\End}{\mathrm{End}}
\newcommand{\Hom}{\mathrm{Hom}}
\newcommand{\Nat}{{\textsc{Hom}}}
\newcommand{\Autnat}{{\textsc{Aut}}}
\newcommand{\Fun}{{\textsc{Fun}}}
\newcommand{\Vect}{{\mathrm{Vect}}}
\newcommand{\vect}{\mathrm{vect}}
\newcommand{\lev}{\mathrm{ev}}
\newcommand{\rev}{\widetilde{\mathrm{ev}}}
\newcommand{\lcoev}{\mathrm{coev}}
\newcommand{\rcoev}{\widetilde{\mathrm{coev}}}
\newcommand{\ldual}[1]{\leftidx{^\vee}{\!#1}{}}
\newcommand{\rdual}[1]{{#1}^\vee}
\newcommand{\lldual}[1]{\leftidx{^{\vee\vee}}{\!#1}{}}
\newcommand{\rrdual}[1]{{#1}^{\vee\vee}}
\newcommand{\lexcla}[1]{\leftidx{^!}{\!#1}{}}
\newcommand{\rexcla}[1]{\leftidx{}{#1}{^!}}
\providecommand{\bysame}{\leavevmode\hbox to3em{\hrulefill}\thinspace}

\begin{document}
\title{Hopf monads}
\author[A. Brugui\`eres]{Alain Brugui\`eres}
\thanks{The first author was partially supported by the CNRS-NSF
project n$^\circ$17149 (Algebraic and homologic methods in low dimensional topology).}
\author[A. Virelizier]{Alexis Virelizier}
\thanks{The second author was partially supported by the European grant
MOIF-CT-2005-008682. He also thanks the Research Institute for Mathematical Sciences in Kyoto, where part of this work was
carried out, for its invitation and hospitality.} \email{bruguier@math.univ-montp2.fr \and virelizi@math.berkeley.edu}
\date{\today}
\begin{abstract}
In this article, we introduce and study Hopf monads on auto\-no\-mous categories (i.e., monoidal categories with duals).
Hopf monads generalize Hopf algebras to a non-braided (and non-linear) setting. Indeed, any monoidal adjunction between
autonomous categories gives rise to a Hopf monad. We extend many fundamental results of the theory of Hopf algebras (such as
the decomposition of Hopf modules, the existence of integrals, Maschke's criterium of semisimplicity, etc...)\@ to Hopf
monads. We also introduce and study quasitriangular and ribbon Hopf monads (again defined in a non-braided setting).
\end{abstract} \maketitle

\setcounter{tocdepth}{1} \tableofcontents

\section*{Introduction}

In 1991, Reshetikhin and Turaev \cite{RT2} introduced  a new $3$-manifold invariant. Its construction consists in
representing the $3$-manifold by surgery along a framed link and  then  assigning a scalar to the link by applying a
suitable algorithm involving simple representations of a quantum group at a root of unity.  Since then, this construction
has been re-visited many times. In particular, Turaev \cite{Tur2} introduced the notion of a modular category, which is a
semisimple ribbon category satisfying conditions of finiteness and non-degeneracy, and showed that such a category defines a
$3$-manifold invariant and indeed a TQFT.

A more general setting for constructing quantum invariants of $3$\trait manifolds has been subsequently developed in
\cite{Lyu2}, and more recently in \cite{LyuKer,Vir}, where the input data is a  (non-necessarily semisimple)  ribbon
category $\cc$ which admits a coend $A=\int^{X \in \cc} \ldual{X} \otimes X$. This coend $A$ is endowed with a very rich
algebraic structure. In particular, it is a Hopf algebra in $\cc$. In fact, in this setting, the quantum invariants depend
only on certain structural morphisms of the coend $A$ (see \cite{BV1} for details).

Recall that a Hopf algebra in a braided category $\cc$ is an object $A$ of $\cc$ which is both an algebra and a coalgebra in
$\cc$, and has an antipode.  The  structural morphisms satisfy the traditional axioms of a Hopf algebra,  except that one
has to replace the usual flip of vector spaces with the braiding $\tau$ of $\cc$. More precisely, the axiom expressing the
compatibility between the product $m \co A \otimes A \to A$ and the coproduct $\Delta\co A \to A \otimes A$ of $A$ becomes:
\begin{equation*}
\Delta m= (m \otimes m)(\id_A \otimes \tau_{A,A} \otimes \id_A)(\Delta \otimes \Delta).
\end{equation*}
Hopf algebras in braided categories have been extensively studied (see \cite{BKLT} and references therein). Many results
about Hopf algebras have been extended to this setting.

However, general (non necessarily braided) monoidal categories also play an important role in quantum topology. Firstly,
they are the input data for another class of  $3$\trait manifold invariants,  the Turaev-Viro invariants (see
\cite{TV,BW}). Also, via the center construction due to Drinfeld, they lead to a large class of braided categories: if $\cc$
is a monoidal category, then its center $Z(\cc)$ is a braided category. Under mild hypotheses, the category $Z(\cc)$ admit a
coend $A$, which is a Hopf algebra in $Z(\cc)$. How can one describe this Hopf algebra $A$ in terms of the category $\cc$?
Note that, if the coend of $\cc$ exists, then it is a coalgebra but not a Hopf algebra (and in general not even an algebra)
and therefore is not sufficient to describe $A$. What we need is a suitable generalization of the notion of Hopf algebra to
a non-braided setting.

The aim of this paper is to provide such a generalization by introducing the notion of Hopf monad.  What is a Hopf monad?
Fix a category $\cc$. Recall that the category $\End(\cc)$ of endofunctors of $\cc$ is a monoidal category, with composition
of endofunctors for monoidal product and trivial endofunctor $1_\cc$  for unit object. A \emph{monad} on $\cc$ is an algebra
in the monoidal category $\End(\cc)$. In other words, it is an endofunctor $T$ of $\cc$ endowed with functorial morphisms
$\mu \co T^2 \to T$ (the product) and $\eta\co 1_\cc \to \cc$ (the unit) such that, for any  object $X$ of $\cc$,
\begin{equation*}
\mu_X \mu_{T(X)}=\mu_X T(\mu_X) \quad \text{and} \quad \mu_X \eta_{T(X)}=\id_{T(X)}=\mu_X T(\eta_X).
\end{equation*}
Note that monads are a very general notion: every monad comes from an adjunction and every adjunction defines a monad. Let
$T$ be a monad on $\cc$. A $T$-module (also called $T$-algebra) is a pair $(M,r)$, where $M$ is an object of $\cc$ and $r\co
T(M) \to M$ is a morphism in $\cc$, satisfying:
\begin{equation*}
r T(r)=r \mu_M \quad \text{and} \quad r \eta_M=\id_M.
\end{equation*}
Denote $T\ti\cc$ the category of $T$-modules and $U_T \co T\ti\cc \to \cc$ the forgetful functor defined by $U_T(M,r)=M$.
Suppose that $\cc$ is monoidal, and denote $\otimes$ its monoidal product and $\un$ its unit object. When does the monoidal
structure of $\cc$ lift to~$T\ti\cc$? The answer to this question lies in the notion of bimonad introduced by
Moerdijk~\cite{Moer}. A \emph{bimonad} is a monad $T$ which is comonoidal, that is, endowed with a functorial morphism
\begin{equation*}
T_2(X,Y)\co T(X\otimes Y) \to T(X) \otimes T(Y)
\end{equation*}
(which plays the role of a coproduct) and a morphism $T_0\co T(\un) \to \un$ (which plays the role of a counit) satisfying
certain compatibility axioms which generalize those of a bialgebra. For example, the axiom which corresponds with the
compatibility of the product and the coproduct is:
\begin{equation*}
T_2(X,Y) \,\mu_{X \otimes Y}= (\mu_X \otimes \mu_Y) \,T_2\bigl(T(X),T(Y)\bigr) \,T\bigl(T_2(X,Y)\bigr).
\end{equation*}
Note that  no braiding is needed  to write this down.

Now the next step is to define the notion of antipode for a bimonad. However, the usual axioms for an antipode cannot be
generalized to bimonads in a straightforward way. In order to bypass this difficulty, we use the categorical interpretation
of an antipode in terms of duality. Let~$\cc$ be a monoidal category which is autonomous, that is, such that each object $X$
of $\cc$ has a left dual $\ldual{X}$ and a right dual $\rdual{X}$. Given a bimonad $T$ on $\cc$, when is the monoidal
category $T\ti\cc$ autonomous? The answer to this question resides in the notion of Hopf monad. A \emph{Hopf monad} on $\cc$
is a bimonad $T$ endowed with a \emph{left antipode} and a \emph{right antipode}, that is, functorial morphisms
\begin{equation*}
s^l_X\co T(\ldual{T(X)}) \to \ldual{X} \quad \text{and} \quad s^r_X\co T(\rdual{T(X)}) \to \rdual{X}
\end{equation*}
which encode the autonomy of $T\ti\cc$. In particular, the left and right duals of a $T$-module $(M,r)$ are given by
\begin{equation*}
\ldual{(M,r)}=(\ldual{M},s^l_M T(\ldual{r})) \quad \text{and} \quad \rdual{(M,r)}=(\rdual{M},s^r_M T(\rdual{r})).
\end{equation*}

The notion of Hopf monad is very general. Firstly Hopf monads generalize Hopf algebras in braided autonomous categories.
Indeed, every Hopf algebra $A$ in a braided autonomous category $\cc$ yields a Hopf monad $T$ on $\cc$ which is defined by
$T(X)=A \otimes X$. In particular a finite-dimensional Hopf algebra $H$ over a field $\kk$ defines a Hopf monad $T(V)= H
\otimes V$ on $\vect(\kk)$. Secondly, Hopf monads do exist in a non-braided setting. In fact, any monoidal adjunction
between autonomous categories gives rise to a Hopf monad. In particular, we give a Tannaka reconstruction theorem in terms
of Hopf monads for fiber functors with values in a category of bimodules.

Quite surprisingly, many fundamental results of the theory of Hopf algebras (such as the decomposition of Hopf modules, the
existence of integrals, Maschke's criterium of  semisimplicity, etc...) extend to the very general setting of Hopf monads,
sometimes in a straightforward way, sometimes at the price of some technical trick. Also, many effective tools for the study
of Hopf algebras turn out to have a monadic counterpart. In particular, we introduce  the notions of sovereign grouplike
element,  \Rt matrix (and Drinfeld element), and twist for a Hopf monad $T$.  These  express the fact  that the
category $T\ti\cc$ of $T$-modules is sovereign, or braided, or ribbon (recall that  $\cc$ itself need not be braided).

In a subsequent paper \cite{BV2}, we  construct the \emph{double} of a Hopf monad on an autonomous category $\cc$, which is
a new Hopf monad on $\cc$.   In particular, the double of the trivial Hopf monad $1_\cc$ is a quasitriangular Hopf monad $D$
on $\cc$, whose category  of modules $D\ti\cc$ is canonically isomorphic (as a braided category) to the center $Z(\cc)$ of
$\cc$. Moreover, we use this quasitriangular Hopf monad to describe the coend of $Z(\cc)$.

The present paper is organized as follows. In Section~\ref{sect-Monads}, we review a few general facts about monads, which
we use intensively throughout the text. In Section~\ref{sect-Bimonads}, we recall the definition of bimonads. In
Section~\ref{sect-HMonads}, we define left and right antipodes and Hopf monads, and we establish their first properties. In
Section~\ref{sect-Hmod}, we introduce Hopf modules and prove the fundamental theorem for Hopf modules over a Hopf monad. We
apply this in Section~\ref{sect-int} to prove a theorem on the existence of universal left and right integrals for a Hopf
monad. In Section~\ref{sect-semsim}, we define semisimple and separable monads, and give a characterization of semisimple
Hopf monads (which generalizes Maschke's theorem). In Section~\ref{sect-sovinv}, we define sovereign grouplike elements and
involutory Hopf monads. In Section~\ref{sect-quasirib}, we define \Rt matrices and twists for a Hopf monad. Finally, in
Section~\ref{sect-perspec}, we give other examples which illustrate the generality of the notion of Hopf monad, including a
Tannaka reconstruction theorem.

\subsection*{Conventions and notations}
Unless otherwise specified, categories are assumed to be small, and monoidal categories are assumed to be strict.

If $\cc$ is a category, we denote $\Ob(\cc)$ the set of objects of $\cc$ and $\Hom_\cc(X,Y)$ the set of morphisms in $\cc$
from an object $X$ to an object $Y$.  The identity functor of $\cc$ will be denoted by $1_\cc$.  We denote $\cc^\opp$ the
\emph{opposite category}  (where arrows are reversed).

Let $\cc$, $\dd$ be two categories. Functors from $\cc$ to $\dd$ are the objects of a category $\Fun(\cc,\dd)$. Given two
functors $F, G\co \cc \to \dd$, a morphism  $\alpha \co F \to G$ is a family $\{\alpha_X\co F(X) \to G(X)\}_{X \in
\Ob(\cc)}$ of morphisms in $\dd$ satisfying the following functoriality condition: $\alpha_Y F(f)=G(f)\alpha_X$ for every
morphism $f\co X \to Y$ in $\cc$. Such a morphism is called a \emph{morphism of functors} or a \emph{functorial morphism}
(the term \emph{natural transformation} is also used in the literature). We denote $\Nat(F,G)$ the set $\Hom_{\Fun
(\cc,\dd)}(F,G)$  of functorial morphisms from $F$ to $G$.

If $\cc$, $\cc'$ are two categories, we denote $\sigma_{\cc,\cc'}$ the flip functor $\cc \times \cc' \to \cc' \times \cc$
defined by $(X, Y) \mapsto (Y, X)$  and $(f,g) \mapsto (g,f)$.

\section{Monads} \label{sect-Monads}

In this section, we review a few general facts about monads, which we use intensively throughout the text.

\subsection{Monads}\label{monad}
Let $\cc$ be a category. Recall that the category $\End(\cc)$ of endofunctors of $\cc$ is strict monoidal with composition
for monoidal product and identity functor $1_\cc$ for unit object. A \emph{monad} on $\cc$  (also called a \emph{triple}) is
an algebra in $\End(\cc)$, that is, a triple $(T,\mu,\eta)$ where $T\co \cc \to \cc$ is a functor, $\mu\co T^2 \to T$ and
$\eta\co 1_\cc \to T$ are functorial morphisms such that:
\begin{align}
& \mu_X T(\mu_X)=\mu_X\mu_{T(X)}; \label{mon-ass} \\
& \mu_X\eta_{T(X)}=\id_{T(X)}=\mu_X T(\eta_X); \label{mon-u}
\end{align}
for all object $X$ of $\cc$.

Let $(T,\mu,\eta)$ be a monad on $\cc$. A \emph{$T$-module} is a pair $(M,r)$ where $M$ is an object of $\cc$ and $r\co T(M)
\to M$ is a morphism in $\cc$ such that:
\begin{equation}\label{Tmoddef}
r T(r)= r \mu_M \quad \text{and} \quad r \eta_M= \id_M.
\end{equation}
Note that $T$-modules are called \emph{$T$-algebras} in \cite{ML1}.

Given two $T$-modules $(M,r)$ and $(N,s)$, a morphism $f\in \Hom_\cc(M,N)$ is said to be \emph{$T$-linear} if $f r=s T(f)$.
Such an $f$ is called a \emph{morphism of $T$-modules} from $(M,r)$ to $(N,s)$. This gives rise to the category $T\ti\cc$ of
$T$-modules (with composition inherited from $\cc$).

We will denote $U_T\co T\ti\cc \to \cc$ the forgetful functor of $T$ defined by $U_T(M,r)=M$  for any $T$-module $(M,r)$ and
$U_T(f)=f$ for any $T$-linear morphism~$f$. Then $U_T$ admits a left adjoint $F_T\co \cc \to T\ti\cc$, which is given by
$F_T(X)=(T(X),\mu_X)$ for any object $X$ of $\cc$ and $F_T(f)=T(f)$ for any morphism $f$ in $\cc$. Note that $T= U_T F_T$ is
the monad of this adjunction.

\begin{exa}\label{Ex1}
Let $\cc$ be a monoidal category and $A$ be an object of $\cc$. Let $A\otimes ?$ be the endofunctor of $\cc$ defined by
$(A\otimes ?)(X)=A \otimes X$ and $(A\otimes ?)(f)=\id_A \otimes f$. Let $m\co A \otimes A \to A$ and $u\co \un \to A$ be
morphisms in $\cc$. Define $\mu_X= m \otimes \id_X$ and $\eta_X=u \otimes \id_X$. Then $(A\otimes ?,\mu,\eta)$ is a monad on
$\cc$ if and only if $(A,m,u)$ is an algebra in $\cc$. If such is the case, then $(A\otimes ?)\ti\cc$ is nothing but the
category of left $A$-modules in $\cc$. Similarly, the endofunctor $? \otimes A$ is a monad on $\cc$ if and only if $A$ is an
algebra in $\cc$ and, if such is the case, $(? \otimes A)\ti \cc$ is the category of right $A$-modules in $\cc$.
\end{exa}
The following classical lemma will be useful later on.
\begin{lem} \label{key-lemma}
Let $T$ be a monad on a category $\cc$ and $U_T\co T\ti\cc \to \cc$ be the forgetful functor. Let $\dd$ be a second category
and $F,G\co \cc \to \dd$ be two functors. Then we have a canonical bijection
\begin{equation*}
?^\sharp \co \Nat(F,GT) \to \Nat(FU_T,G U_T), \; f \mapsto f^\sharp
\end{equation*}
defined by $f^\sharp_{(M,r)}=G(r)f_M$ for any $T$-module $(M,r)$. Its inverse
\begin{equation*}
?^\flat \co \Nat(FU_T,G U_T) \to \Nat(F,GT), \; g \mapsto g^\flat
\end{equation*}
is given by $g^\flat_X=g_{(T(X),\mu_X)}F(\eta_X)$ for any object $X$ of $\cc$. If $F,G$ are contravariant functors (that is,
functors from $\cc^\opp$ to $\dd$), then the bijection becomes:
\begin{equation*}
?^\sharp \co\Nat(GT^\opp,F) \to \Nat(GU_T^\opp,F U_T^\opp)
\end{equation*}
with $f^\sharp_{(M,r)}=f_M G(r)$ and $g^\flat_X=F(\eta_X) g_{(T(X),\mu_X)}$.
\end{lem}

\begin{proof}
Let us verify the covariant case. We first remark that $?^\sharp$ and $?^\flat$ are well-defined. Indeed $g^\flat$ is
clearly functorial, and $f^\sharp$ is functorial since, for any $T$-linear morphism $h\co (M,r) \to (N,s)$, we have
$f^\sharp_{(N,s)} F(h)= G(s) f_N F(h) = G(s T(h)) f_M= G(h r)f_M= G(h) f^\sharp_{(M,r)}$. Now
$f^{\sharp\flat}_{X}=G(\mu_X)f_{T(X)}F(\eta_X)=G(\mu_X)G(T\eta_X)f_X=f_X$ for any object $X$ of $\cc$, and
$g^{\flat\sharp}_{(M,r)}=G(r)g_{(T(M),\mu_M)}F(\eta_M)=g_{(M,r)} F(r)F(\eta_M)=g_{(M,r)}$ for any $T$-module $(M,r)$. Hence
$?^\sharp$ and $?^\flat$ are inverse each other. The contravariant case is a mere reformulation.
\end{proof}

Given a functor $F\co \cc \to \dd$ and a positive integer $n$, we denote $\cc^n=\cc \times \cdots \times \cc$ and $F^{\times
n}\co \cc^n \to \dd^n$ the $n$-uple cartesian product of $\cc$ and $F$. Note that if $T$ is a monad on a category $\cc$,
then $T^{\times n}$ is a monad on $\cc^n$, and we have  $T^{\times n} \ti \cc^n=(T\ti \cc)^n$ and $U_{T^{\times
n}}=(U_T)^{\times n}$. Re-writing Lemma~\ref{key-lemma} for this monad, we get:

\begin{lem} \label{key-lemmaN}
Let $T$ be a monad on a category $\cc$ and $U_T\co T\ti\cc \to \cc$ be the forgetful functor. Fix an positive integer $n$.
Let $\dd$ be a second category and $F,G\co \cc^n \to \dd$ be two functors. Then we have a canonical bijection
\begin{equation*}
?^\sharp \co \Nat(F,GT^{\times n}) \to \Nat(FU_T^{\times n},G U_T^{\times n}), \; f \mapsto f^\sharp
\end{equation*}
defined by $f^\sharp_{(M_1,r_1), \dots, (M_n,r_n)}=G(r_1, \dots, r_n)f_{M_1, \dots,M_n}$. Its inverse
\begin{equation*}
?^\flat \co \Nat(FU_T^{\times n},G U_T^{\times n}) \to \Nat(F,GT^{\times n}), \; g \mapsto g^\flat
\end{equation*}
is given by $g^\flat_{X_1, \dots, X_n}=g_{(T(X_1),\mu_{X_1}), \dots, (T(X_n),\mu_{X_n})}F(\eta_{X_1}, \dots, \eta_{X_n})$.
The contravariant case can be stated similarly (see Lemma~\ref{key-lemma}).
\end{lem}

\subsection{Convolution product} \label{convol}
Let $\cc$, $\dd$ be two categories and $(T,\mu,\eta)$ be a monad on $\cc$. Let $F, G, H$ be three functors $\cc^n \to\dd$.
Let $f \in \Nat(F, G T^{\times n})$ and $g \in \Nat(G, H T^{\times n})$. Define their \emph{convolution product} $g*f \in
\Nat(F, H T^{\times n})$ by setting, for any objects $X_1, \dots, X_n$ of $\cc$,
\begin{equation}\label{convoprodini}
(g*f)_{X_1, \dots, X_n}= H(\mu_{X_1}, \dots, \mu_{X_n}) \,g_{T(X_1), \dots, T(X_n)} \,f_{X_1, \dots, X_n}.
\end{equation}
This convolution product reflects the composition of morphisms in the category $\Fun(\cc^n, \dd)$ via the canonical
bijection $\Nat(F, G T^{\times n})\simeq \Nat(F U_T^{\times n}, G U_T^{\times n})$ given by Lemma \ref{key-lemmaN}.

We say that  $f \in \Nat(F, G T^{\times n})$  is \emph{$*$-invertible} if there exists  $g \in \Nat(G, F T^{\times n})$ such
that $g*f=F(\eta^{\times n})$ and $f*g= G(\eta^{\times n})$.  This means that $f^\sharp$ is an isomorphism of functors, with
inverse $g^\sharp$. If such a $g$ exists, then it is unique and we denote it~$f^{*-1}$.

\subsection{Central elements} \label{sect-centralelts}
Let $T$ be a monad on a category $\cc$. By Section~\ref{convol}, the set $\Nat(1_\cc,T)$ is a monoid, with  unit~$\eta$, for
the convolution product $*$ defined, for any $\phi,\psi\in\Nat(1_\cc,T)$, by
\begin{equation}\label{prodconvo1T}
(\phi * \psi)_X=\mu_X \phi_{T(X)} \psi_X =\mu_X T(\psi_X) \phi_X \co X \to T(X).
\end{equation}
Recall that, via the canonical bijection $?^\sharp\co\Nat(1_\cc,T) \to\Nat(U_T,U_T)$ of Lem\-ma~\ref{key-lemma}, this
convolution product corresponds with composition of endomorphisms of the forgetful functor $U_T\co T\ti\cc \to T$.

Given $a \in \Nat(1_\cc,T)$, let $L_a, R_a \in \Nat(T,T)$ be the functorial morphisms defined, for any object $X$ of $\cc$,
by:
\begin{equation}\label{defLR}
(L_a)_X=\mu_X a_{T(X)} \quad \text{and} \quad (R_a)_X=\mu_X T(a_X).
\end{equation}
Remark that $L_a b=a*b$ and $R_a b=b*a$ for all $a,b \in \Nat(1_\cc,T)$.

A \emph{central element} of $T$ is a functorial morphism $a\in\Nat(1_\cc,T)$ such that $L_a=R_a$. For example, by
\eqref{mon-u}, the unit $\eta$ of $T$ is a central element. Notice that any central element of $T$ is in particular central
in the monoid $(\Nat(1_\cc,T),*,\eta)$.

\begin{lem} \label{lemdefcentral}
Let $T$ be a monad on a category $\cc$ and $a\in\Nat(1_\cc,T)$. The following conditions are equivalent:
\begin{enumerate}
\renewcommand{\labelenumi}{{\rm (\roman{enumi})}}
\item  The morphism $a$ is a central element of $T$;
\item For any $T$-module $(M,r)$, the morphism $a^\sharp_{(M,r)} \co M \to M$ is $T$-linear;
\item There exists a (necessarily unique) functorial morphism $\tilde{a}\co 1_{T\ti \cc} \to 1_{T\ti\cc}$ such that $U_T (\tilde{a})=
a^\sharp$.
\end{enumerate}
\end{lem}
\begin{proof}
Clearly, (ii) is equivalent to (iii). Let $(M,r)$ be a $T$-module. Then, using~\eqref{Tmoddef},
\begin{equation*}
a^\sharp_{(M,r)} r =r a_M r = r T(r) a_{T(M)} = r \mu_M a_{T(M)} =  r (L_a)_M
\end{equation*}
and
\begin{equation*}
r T(a^\sharp_{(M,r)}) = r T(r)T(a_M)= r \mu_M T(a_M) = r (R_a)_M.
\end{equation*}
Therefore (i) is equivalent to (ii) by Lemma~\ref{key-lemma}.
\end{proof}

\subsection{The adjoint action} \label{sect-adj}
Let $T$ be a monad on a category $\cc$ and consider the maps $L,R\co \Nat(1_\cc,T) \to \Nat(T,T)$ defined as in
\eqref{defLR}.

\begin{lem}\label{lemcentr}
The maps $L$ and $R$ are respectively a homomorphism and an anti-ho\-mo\-mor\-phism of monoids from $(\Nat(1_\cc,T),*,\eta)$
to $(\Nat(T,T),\circ,\id_T)$. Moreover $L_aR_b=R_bL_a$ for all $a,b \in \Nat(1_\cc,T)$.
\end{lem}
\begin{proof}
For any object $X$ of $\cc$, we have $(L_\eta)_X=\mu_X\eta_{T(X)}=\id_{T(X)}$ by \eqref{mon-u} and, given
$a,b\in\Nat(1_\cc,T)$,
\begin{align*}
(L_{a *b})_X &=\mu_X\mu_{T(X)}a_{T^2(X)}b_{T(X)} \\
& =\mu_X T(\mu_X) a_{T^2(X)}b_{T(X)} \quad \text{by \eqref{mon-ass}}\\
&= \mu_X a_{T(X)}\mu_X b_{T(X)}=(L_a)_X(L_b)_X.
\end{align*}
Therefore $L$ is a homomorphism of monoids. Likewise one shows that $R$ anti-ho\-mo\-mor\-phism of monoids. Finally, given
$a,b\in\Nat(1_\cc,T)$ and an object $X$ of $\cc$, we have:
\begin{align*}
(L_a R_b)_X &=\mu_X a_{T(X)} \mu_X T(b_X) \\
& = \mu_X T(\mu_X) a_{T^2(X)} T(b_X)\\
& = \mu_X \mu_{T(X)} T^2(b_X) a_{T(X)} \quad \text{by \eqref{mon-ass}}\\
&= \mu_X T(b_X)\mu_X a_{T(X)}=(R_b)_X(L_a)_X,
\end{align*}
and so $L_aR_b=R_bL_a$.
\end{proof}
Given $a\in \Nat(1_\cc,T)$, we get from Lemma~\ref{lemcentr} that $L_a$ (resp.\@ $R_a$) is invertible if and only if $a$ is
$*$-invertible and, if such the case, $L^{-1}_a=L_{a^{*-1}}$ (resp.\@ $R^{-1}_a=R_{a^{*-1}}$). Denote $\Autnat(T)$ the group
of functorial automorphisms of $T$ and $\Nat(1_\cc,T)^\times$ the group of $*$-invertible elements of the monoid
$(\Nat(1_\cc,T),*,\eta)$. Define:
\begin{equation}\label{defad}
\ad\co \left \{
\begin{array}{ccl}
  \Nat(1_\cc,T)^\times & \to & \Autnat(T) \\
  a & \mapsto & \ad_a=L_aR_{a^{*-1}}=R_{a^{*-1}}L_a
\end{array}\right .
\end{equation}
The map $\ad$ is a group morphism (by Lemma~\ref{lemcentr}) and is called the \emph{adjoint action} of $T$. Its kernel is
made of the $*$-invertible central elements of $T$.

Notice that $\ad_a b=a*b*a^{*-1}$ for all $b \in \Nat(1_\cc,T)$.

\subsection{Morphisms of monads}\label{sect-morphmon}
A \emph{morphism of monads} between two monads $(T,\mu,\eta)$ and $(T',\mu',\eta')$ on a category $\cc$ is a functorial
morphism $f \co T \to T'$ such that:
\begin{equation}\label{morphmon}
f_X \mu_X=\mu'_X  f_{T'(X)}T(f_X) \quad \text{and} \quad f_X \eta_X=\eta'_X
\end{equation}
for any object $X$ of $\cc$.

\begin{lem}\label{lemmorphmon}
Let $T$ and $T'$ be two monads on a category $\cc$. Let $f\co T \to T'$ be a functorial morphism. The following conditions
are equivalent:
\begin{enumerate}
\renewcommand{\labelenumi}{{\rm (\roman{enumi})}}
\item $f\co T \to T'$ is a morphism of monads;
\item For all $T'$-module $(M,r)$, the pair $(M, r f_M)$ is a $T$-module.
\end{enumerate}
Moreover, if $f$ is a morphism of monads, then the assignment $(M,r) \mapsto (M, r f_M)$ defines a functor $f^* \co T'\ti\cc
\to T\ti \cc$ which satisfies $U_T f^* = U_{T'}$. Lastly, for any functor $F \co T'\ti\cc \to T\ti \cc$ such that $U_T F =
U_{T'}$, there exists a unique morphism of monads $f \co T \to T'$ such that $F=f^*$.
\end{lem}

\begin{proof}
Let $(M,r)$ be a $T'$-module. The pair $(M, r f_M)$ is a $T$-module if and only if $rf_M\mu_M= r f_M T(rf_M) = r T'(r)
f_{T'(M)} T(f_M)$ and $r f_M \eta_M=\id_M$. Using \eqref{Tmoddef}, this is equivalent to $rf_M\mu_M= r \mu'_M f_{T'(M)}
T(f_M)$ and $r f_M \eta_M=r \eta_M$. By Lemma~\ref{key-lemma}, this holds for all $T'$-module $(M,r)$ if and only if $f$ is
a morphism of monads. Clearly, if $f$ is a morphism of monads, then $f^*$ is a well-defined functor. Let $F \co T'\ti \cc
\to T\ti \cc$ be a functor such that $U_T F = U_{T'}$. For any object $X$ of $\cc$, let $\nu_X\co T(T'(X)) \to T'(X)$ be
such that $F(T'(X), \mu'_X)= (T'(X), \nu_X)$. Since $F$ is a functor, $\nu\co TT' \to T'$ is a functorial morphism. Firstly,
if $F=f^*$ for some morphism of monads $f\co T \to T'$, then $f$ is unique since $\nu_X= \mu'_X f_{T'(X)}$ and so
$f_X=\mu'_X T'(\eta'_X) f_X= \mu'_X f_{T'(X)}T(\eta'_X)=\nu_X T(\eta'_X)$ for all object $X$ of $\cc$. Conversely define
$f=\nu T(\eta')\co T \to T'$. Let $(M,r)$ be a $T'$-module and $s\co T(M) \to M$ such that $F(M,r)=(M,s)$. Let us show $s=r
f_M$. This will prove simultaneously that $f$ is a morphism of monads and $F=f^*$. Since $r\co (T'(M),\mu'_M)\to (M,r)$ is a
morphism of $T'$-modules, we have  $s T(r)= r \nu_M$ (by functoriality of $F$). Hence $r f_M= r \nu_M T(\eta'_M)=s
T(r\eta'_M)=s$.
\end{proof}

\section{Bimonads}  \label{sect-Bimonads}

In this section, we review the definition and properties of a bimonad. This notion was first introduced by
Moerdijk~\cite{Moer}.

\subsection{(Co-)monoidal functors}

Let $(\cc,\otimes,\un)$ and $(\dd, \otimes, \un)$ be two monoidal categories. A  \emph{monoidal functor} from $\cc$ to $\dd$
is a triple $(F,F_2,F_0)$ where $F\co \cc \to \dd$ is a functor, $F_2\co F\otimes F \to F \otimes$ is a morphism of
functors, and $F_0\co\un \to F(\un)$ is a morphism in $\dd$ such that:
\begin{align}
& F_2(X,Y \otimes Z) (\id_{F(X)} \otimes F_2(Y,Z))= F_2(X \otimes Y, Z)(F_2(X,Y) \otimes \id_{F(Z)}) ;\\
& F_2(X,\un)(\id_{F(X)} \otimes F_0)=\id_{F(X)}=F_2(\un,X)(F_0 \otimes \id_{F(X)}) ;
\end{align}
for all objects $X,Y,Z$ of $\cc$.

A \emph{comonoidal functor} from $\cc$ to $\dd$ is a triple $(F,F_2,F_0)$ where $F\co \cc \to \dd$ is a functor, $F_2\co F
\otimes \to F\otimes F$ is a morphism of functors, and $F_0\co F(\un) \to \un$ is a morphism in $\dd$ such that:
\begin{align}
& (\id_{F(X)} \otimes F_2(Y,Z))  F_2(X,Y \otimes Z)= (F_2(X,Y) \otimes \id_{F(Z)}) F_2(X \otimes Y, Z) ;\label{comoneq1}\\
& (\id_{F(X)} \otimes F_0) F_2(X,\un)=\id_{F(X)}=(F_0 \otimes \id_{F(X)}) F_2(\un,X) ;\label{comoneq2}
\end{align}
for all objects $X,Y,Z$ of $\cc$.  Comonoidal functors are sometimes called \emph{opmonoidal} in the literature.

A (co-)monoidal functor $(F,F_2,F_0)$ is said to be \emph{strong} (resp.\@ \emph{strict}) if $F_2$ and $F_0$ are
isomorphisms (resp.\@ identities). For example, the identity functor $1_\cc$ of a monoidal category $\cc$ is a strict
(co-)monoidal functor.

Given a functor $F:\cc \to \dd$, a functorial isomorphism $F_2\co F \otimes \to F\otimes F$, and an isomorphism $F_0\co\un
\to F(\un)$, the triple $(F,F_2,F_0)$ is a monoidal functor if and only if $(F,F_2^{-1},F_0^{-1})$ is a comonoidal functor.

\begin{lem}\label{compo}
Let $F$ and $G$ be two composable functors between monoidal categories.
\begin{enumerate}
  \renewcommand{\labelenumi}{{\rm (\alph{enumi})}}
  \item If $F$ and $G$ are monoidal functors, then $GF$
is a monoidal functor with $(GF)_2=G(F_2)G_2$ and $(GF)_0=G(F_0)G_0$.
  \item If $F$ and $G$ are comonoidal functors, then $GF$
is a comonoidal functor with $(GF)_2=G_2G(F_2)$ and $(GF)_0=G_0G(F_0)$.
\end{enumerate}
\end{lem}

\subsection{(Co-)monoidal morphisms of functors}
A functorial morphism $\varphi\co F \to G$ between two monoidal functors $F\co\cc \to \dd$ and $G\co\cc \to \dd$ is
\emph{monoidal} if it satisfies:
\begin{equation}
\varphi_{X \otimes Y}  F_2(X,Y)= G_2(X,Y)  (\varphi_X \otimes \varphi_Y) \quad \text{and} \quad G_0=\varphi_\un  F_0
\end{equation}
for all objects $X,Y$ of $\cc$.

Likewise, a functorial morphism $\varphi\co F \to G$ between two comonoidal functors $F\co\cc \to \dd$ and $G\co\cc \to \dd$
is \emph{comonoidal} if it satisfies:
\begin{equation}
G_2(X,Y)  \varphi_{X \otimes Y}= (\varphi_X \otimes \varphi_Y)  F_2(X,Y)\quad \text{and} \quad G_0  \varphi_\un= F_0
\end{equation}
for all objects $X,Y$ of $\cc$.
\subsection{Bimonads}
A \emph{bimonad} on a monoidal category $\cc$ is a monad $(T,\mu,\eta)$ on $\cc$ such that the functor $T\co \cc \to \cc$ is
comonoidal and the functorial  morphisms  $\mu\co T^2 \to T$ and $\eta\co 1_\cc \to T$ are comonoidal. Here $1_\cc$
is the (strict) comonoidal identity functor of $\cc$ and $T^2$ is the comonoidal functor obtained by composition of the
comonoidal functor $T$ with itself  as in Lemma~\ref{compo}(b).  Explicitly, $\mu$ and $\eta$ are comonoidal if they
satisfy, for any objects $X,Y$ of $\cc$,
\begin{align}
& T_2(X,Y) \mu_{X \otimes Y}= (\mu_X \otimes \mu_Y) T_2(T(X),T(Y)) T(T_2(X,Y)) ;\label{bimonad1}\\
& T_0 \mu_\un= T_0 T(T_0) ; \label{bimonad2}\\
& T_2(X,Y) \eta_{X \otimes Y}= (\eta_X \otimes \eta_Y) ; \label{bimonad3}\\
& T_0 \eta_\un= \id_\un. \label{bimonad4}
\end{align}

Our notion of bimonad coincides exactly with the notion of `Hopf monad' introduced in \cite{Moer}. However, by analogy with
the notions  of bialgebra and Hopf algebra, we prefer to reserve the term `Hopf monad' for bimonads with antipodes (see
Section~\ref{sect-antipodes}). This choice may be justified by the following example:

\begin{exa}\label{Ex2}
Let $\cc$ be a braided category, with braiding $\tau_{X,Y}\co X \otimes Y \to Y \otimes X$, and $A$ be an algebra in $\cc$.
Let $T=A\otimes ?$ be its associated monad on $\cc$, see Example~\ref{Ex1}.
 Let $\Delta\co A \to A \otimes A$ and
$\varepsilon \co A \to \un$ be morphisms in $\cc$. Set
\begin{align*}
& T_2(X,Y)=(\id_A \otimes \tau_{A,X} \otimes \id_Y)(\Delta \otimes \id_{X \otimes Y}) \co T(X \otimes Y) \to T(X) \otimes T(Y);\\
& T_0=\varepsilon \co T(\un) \to \un.
\end{align*}
Then $(T,T_2,T_0)$ is a bimonad on $\cc$ if and only if $(A,\Delta,\varepsilon)$ is a bialgebra in $\cc$. Similarly, given a
bialgebra $A$ in~$\cc$, the endofunctor $? \otimes A$ is a bimonad on $\cc$ with
\begin{equation*}
(? \otimes A)_2(X,Y)=(\id_X \otimes \tau_{Y,A} \otimes \id_A)(\id_{X \otimes Y} \otimes \Delta) \quad \text{and} \quad (?
\otimes A)_0=\varepsilon.
\end{equation*}
In particular, any bialgebra $H$ over a field $\kk$ defines bimonads $H \otimes_\kk ?$ and $? \otimes_\kk H$ on the category
$\Vect(\kk)$ of $\kk$-vector spaces.
\end{exa}

We can reformulate the main result of \cite{Moer} as follows (see also \cite{Crud1}):
\begin{thm}[Moerdijk, 2002]\label{moerthm}
Let $T$ be a monad on a monoidal category~$\cc$. If $T$ is a bimonad, then the category $T\ti\cc$ of $T$-modules is monoidal
by setting:
\begin{equation*}
(M,r) \otimes_{T\ti\cc} (N,s)=(M \otimes N, (r \otimes s) T_2(M,N)) \quad \text{and} \quad \un_{T\ti\cc}=(\un,T_0).
\end{equation*}
Moreover this gives a bijective correspondence between:
\begin{itemize}
 \item bimonad structures for the monad $T$;
 \item monoidal structures of $T\ti\cc$ such that the forgetful functor $U_T \co
T\ti\cc \to \cc$ is strict monoidal.
\end{itemize}
\end{thm}

\begin{exa}\label{exa-szl}
Szlach{\'a}nyi has shown that (left) bialgebroid, as defined in \cite{Szl}, may be interpreted in terms of bimonads. More
precisely, let $\kk$ is a commutative ring and $B$ a \kt algebra. Denote $\biMod{B}{B}$ the category of $B$-bimodules, which
is monoidal with tensor product $\otimes_B$ and unit object $\leftidx{_B}{B}{_B}$. Then following data are equivalent:
\begin{itemize}
\item bimonads on $\biMod{B}{B}$ which commute with inductive limits;
\item (left) bialgebroids with base $B$.
\end{itemize}
If $A$ is a (left) bialgebroid, then the corresponding bimonad is $T= A \otimes_B ?$ and the monoidal categories
$T\ti\biMod{B}{B}$ and $\lMod{A}$ are equivalent. Note that in general the monoidal category $\biMod{B}{B}$ is not braided.
\end{exa}

\begin{rem}\label{mon-opp}
Let $T$ be a bimonad on a monoidal category $\cc=(\cc,\otimes,\un)$. Then $T$ can be viewed as a bimonad $T^\opp$ on the
monoidal category $\cc^{\otimes \opp}=(\cc,\otimes^\opp,\un)$, with comonoidal structure $T_2^\opp=T_2\sigma_{\cc,\cc}$ and
$T_0^\opp=T_0$. The bimonad $T^\opp$  is called  the \emph{opposite} of the bimonad $T$. We have $T^\opp
\ti\cc^{\otimes\opp}=(T\ti \cc)^{\otimes \opp}$.
\end{rem}

\begin{rem}
Notice that the notion of bimonad is not `self-dual': one may define a \emph{bi-comonad} on a monoidal category $\cc$ to be
a bimonad of the opposite category~$\cc^\opp$.
\end{rem}

\subsection{Morphisms of bimonads}\label{sect-morphbimon}
A \emph{morphism of bimonads} between two bimonads $T$ and $T'$ on a monoidal category $\cc$ is a morphism of monads $f\co T
\to T'$ (see Section~\ref{sect-morphmon}) which is comonoidal.

\begin{lem}\label{lemmorphbimon}
Let $T$ and $T'$ be two bimonads on a monoidal category $\cc$. Let $f\co T \to T'$ be a morphism of monads. Then $f$ is a
morphism of bimonads if and only if the functor $f^* \co T'\ti\cc \to T\ti \cc$ induced by $f$ (see Lemma~\ref{lemmorphmon})
is monoidal strict. Moreover, for any strict monoidal functor $F \co T'\ti\cc \to T\ti \cc$ such that $U_T F = U_{T'}$,
there exists a unique morphism of bimonads $f \co T \to T'$ such that $F=f^*$.
\end{lem}

\begin{proof}
In view of Lemma~\ref{lemmorphmon}, we have to show that the functor $F=f^*$ is monoidal strict if and only if $f$ is
comonoidal. We have $F(\un, T'_0)= (\un, T'_0 f_\un)$ and
\begin{align*}
& F((M,r)\otimes(N,s))=(M \otimes N, (r \otimes s) T'_2(M,N) f_{M \otimes N}), \\
& F(M,r) \otimes F(N,s)=(M \otimes N, (r \otimes s)(f_M \otimes f_N)T_2(M,N)),
\end{align*}
for any $T'$-bimodules $(M,r)$ and $(N,s)$. We conclude by Lemma~\ref{key-lemmaN}.
\end{proof}

\section{Hopf monads}\label{sect-HMonads}
Let $T$ be a bimonad on a monoidal category $\cc$. By Theorem~\ref{moerthm}, the category $T\ti\cc$ of $T$-modules is
monoidal and the forgetful functor $U_T \co T\ti \cc \to \cc$ is strict monoidal. Assuming  $\cc$ is autonomous (i.e., has
duals), when is it true that $T\ti\cc$ is autonomous as well? The answer lies in the notions of antipode and Hopf monad,
which we introduce in this section. We first recall some properties of autonomous categories.

\subsection{Autonomous categories}\label{sect-autono}
Recall that a \emph{duality} in a monoidal category $\cc$ is a quadruple $(X,Y,e,d)$, where $X$, $Y$ are objects of $\cc$,
$e\co X \otimes Y \to \un$ (the \emph{evaluation}) and $d \co \un \to Y \otimes X$ (the \emph{coevaluation}) are morphisms
in $\cc$, such that:
\begin{equation}\label{axiomdual}
(e \otimes \id_X)(\id_X \otimes d)=\id_X \quad \text{and} \quad (\id_Y \otimes e)(d \otimes \id_Y)=\id_Y.
\end{equation}
Then $(X,e,d)$ is a \emph{left dual} of $Y$, and $(Y,e,d)$ is a \emph{right dual} of $X$.

If $D=(X,Y,e,d)$ and $D'=(X',Y',e',d')$ are two dualities, two morphisms $f\co X \to X'$ and $g\co Y' \to Y$ are \emph{in
duality with respect to $D$ and $D'$} if
\begin{equation*}
e'(f \otimes \id_{Y'})=e(\id_X \otimes g) \quad \bigl (\text{or, equivalently, $(\id_{Y'} \otimes f)d=(g \otimes \id_X)
d'$}\bigr).
\end{equation*}
In that case we write $f= \ldual{g}_{D,D'}$ and $g = \rdual{f}_{D,D'}$, or simply $f=\ldual{g}$ and $g=\rdual{f}$ if the
context justifies a more relaxed notation. Note that this defines a bijection between $\Hom_\cc(X,X')$ and $\Hom_\cc(Y',Y)$.

Left and right duals, if they exist, are essentially unique: if $(Y,e,d)$ and $(Y',e',d')$ are right duals of some object
$X$, then there exists a unique isomorphism $u\co Y \to Y'$ such that $e'=e(\id_X \otimes u^{-1})$ and $d'=(u \otimes
\id_X)d$.

A \emph{left autonomous} (resp.\@  \emph{right autonomous}, resp.\@ \emph{autonomous}) category is a monoidal category for
which every object admits a left dual (resp.\@  a right dual, resp.\@ both a left and a right dual).  Note that autonomous
categories are also called rigid categories in the literature.

Assume $\cc$ is a left autonomous category and, for each object $X$, pick a left dual $(\ldual{X},\lev_X,\lcoev_X)$. This
data defines a strong monoidal functor $\ldual{?}\co \cc^\opp \to \cc$, where $\cc^\opp$ is the opposite category to $\cc$
with opposite monoidal structure. This monoidal functor is called the \emph{left dual functor}. Notice that the actual
choice of left duals is innocuous in the sense that different choices of left duals define canonically isomorphic left dual
functors.

Likewise, if $\cc$ is a right autonomous category, picking a right dual $(\rdual{X},\rev_X,\rcoev_X)$ for each object $X$
defines a strong monoidal functor $\rdual{?}\co \cc^\opp\to \cc$, called the \emph{right dual functor}.

\begin{rem}
Subsequently, when dealing with left or right autonomous categories, we shall always assume tacitly that left duals or right
duals have been chosen. In formulae, we will often abstain (by abuse) from writing down the following canonical
isomorphisms:
\begin{align*}
& \ldual{?}_2(X,Y)\co \ldual{Y} \otimes \ldual{X} \to \ldual{(X \otimes Y)}, & \ldual{?}_0 \co \un \to \ldual{\un}, \\
& \rdual{?}_2(X,Y)\co \rdual{Y} \otimes \rdual{X} \to \rdual{(X \otimes Y)}, & \rdual{?}_0 \co \un \to \rdual{\un}.
\end{align*}
\end{rem}

\begin{rem}\label{qieq}
If $\cc$ is autonomous, then the functors ${\rdual{?}}^\opp$ and $\ldual{?}$ are canonically quasi-inverse. More precisely,
for any object $X$ of $\cc$, we have the following canonical functorial isomorphisms:
\begin{align*}
& (\rev_X \otimes \id_{\ldual{(\rdual{X})}})(\id_X \otimes \lcoev_{\rdual{X}}) \co X \to \ldual{(\rdual{X})},\\
&(\id_{\rdual{(\ldual{X})}} \otimes \lev_X)(\rcoev_{\ldual{X}} \otimes \id_X) \co X \to \rdual{(\ldual{X})}.
\end{align*}
Again, we will often abstain from writing down these isomorphisms.
\end{rem}

Let $\cc$, $\dd$ be autonomous categories (with chosen left and right duals). For any functor $F\co \cc \to\dd$, we define
two functors $\lexcla{F}\co \cc \to \dd$  and $\rexcla{F}\co \cc \to \dd$  by setting:
\begin{equation*}
\lexcla{F}(X)=\ldual{F(\rdual{X})}, \quad \lexcla{F}(f)=\ldual{F(\rdual{f})}, \quad \rexcla{F}(X)=\rdual{F(\ldual{X})},
\quad \rexcla{F}(f)=\rdual{F(\ldual{f})}
\end{equation*}
for all object $X$  and all morphism $f$ in $\cc$. For any functorial morphism $\alpha\co F \to G$ between functors from
$\cc$ to $\dd$, we define two functorial morphisms $\lexcla{\alpha}\co \lexcla{G} \to \lexcla{F} $ and $\rexcla{\alpha}\co
\rexcla{G} \to \rexcla{F} $ by setting:
\begin{equation*}
\lexcla{\alpha}_X=\ldual{(\alpha_{\rdual{X}})}\co \lexcla{G}(X) \to \lexcla{F}(X) \quad \text{and} \quad
\rexcla{\alpha}_X=\rdual{(\alpha_{\ldual{X}})}\co \rexcla{G}(X) \to \rexcla{F}(X)
\end{equation*}
for all object $X$ of $\cc$. These assignments lead to  functors
\begin{equation*}
\lexcla{?} \co \Fun(\cc,\dd)^\opp \to \Fun(\cc,\dd)\quad \text{and} \quad \rexcla{?} \co \Fun(\cc,\dd)^\opp \to
\Fun(\cc,\dd).
\end{equation*}
\begin{rem}\label{rem-excla}
The functors $\lexcla{?}$ and $\rexcla{?}$ enjoy the following properties: given two composable functors $F \co \cc \to \dd$
and $G \co \dd \to \ee$ between autonomous categories, we have $\lexcla{1_\cc}=1_\cc=\rexcla{1_\cc}$,
$\lexcla{(FG)}=\lexcla{F}\lexcla{G}$ and $\rexcla{(FG)}=\rexcla{F}\rexcla{G}$ (up to the canonical isomorphisms of
Remark~\ref{qieq}). Moreover, the functors $\lexcla{?}$ and $(\rexcla{?})^\opp$ are  quasi-inverse to each other.
\end{rem}

\subsection{Strong monoidal functors and duality}\label{strongfunctdual}
Let $F\co \cc \to \dd$ be a strong monoidal functor between monoidal categories. If $D=(X,Y,e,d)$ is a duality in $\cc$,
then
\begin{equation*}
 F(D)=(F(X),F(Y), F^{-1}_0 F(e)F_2(X,Y),F_2(Y,X)^{-1}F(h) F_0)
\end{equation*}
is a duality in $\dd$. In particular, if $\cc$ and $\dd$ are left (resp.\@  right) autonomous, we have a canonical
isomorphism:
\begin{equation*}
F^l_1(X)\co F(\ldual{X}) \iso \ldual{F(X)} \quad \bigl (\text{resp. }  F^r_1(X)\co F(\rdual{X}) \iso \rdual{F(X)} \bigr).
\end{equation*}

\begin{lem} \label{isomono}
Let $F, G  \co \cc \to \dd$ be strong monoidal functors and $\alpha\co F \to G$ be a monoidal morphism. If $\cc$ is left or
right autonomous, then $\alpha$ is an isomorphism. More precisely, if $\cc$ is left (resp.\@ right) autonomous then, for any
object $X$ of $\cc$,
\begin{equation*}
\alpha^{-1}_X=\rdual{(\alpha_{\ldual{X}})}_{G(D),F(D)} \quad \bigl(\text{resp.\@ }
\alpha^{-1}_X=\ldual{(\alpha_{\rdual{X}})}_{F(D),G(D)} \bigr )
\end{equation*}
where $D$ is the duality $(\ldual{X},X,\lev_X,\lcoev_X)$ (resp.\@ $(X,\rdual{X},\rev_X,\rcoev_X)$).  In particular, if $\cc$
is autonomous, then $\alpha^{-1}=\rexcla{\alpha}=\lexcla{\alpha}$ (up to the canonical isomorphisms of Remark~\ref{qieq}),
where $\lexcla{?}$ and $\rexcla{?}$ are the functors of Remark~\ref{rem-excla}.
\end{lem}

\begin{proof}
Let $D=(X,Y,e,d)$ be a duality in $\cc$. Denote $\ldual{\alpha_Y}\co G(X) \to F(X)$ the left dual of $\alpha_Y \co F(Y) \to
G(Y)$ with respect to the dualities $F(D)=(F(X),F(Y),e_F,d_F)$ and $G(D)=(G(X),G(Y),e_G,d_G)$. We have
\begin{align*}
e_G(\alpha_X \otimes \alpha_Y)& = G^{-1}_0 G(e) G_2(X,Y)(\alpha_X \otimes \alpha_Y)=G^{-1}_0 G(e) \alpha_{X \otimes
Y}F_2(X,Y)\\
& =G^{-1}_0 \alpha_\un F(e) F_2(X,Y)=F^{-1}_0 F(e) F_2(X,Y)=e_F
\end{align*}
and similarly $(\alpha_Y \otimes \alpha_X)d_F= d_G$. Now
\begin{align*}
\ldual{\alpha_Y} \alpha_X &= (e_G (\alpha_X \otimes \alpha_Y) \otimes \id_{F(X)})(\id_{F(X)} \otimes d_F)\\
&= (e_F \otimes \id_{F(X)})(\id_{F(X)} \otimes d_F)= \id_{F(X)}
\end{align*}
and
\begin{align*}
\alpha_X\ldual{\alpha_Y}& = (e_G \otimes \id_{G(X)})(\id_{G(X)} \otimes (\alpha_Y \otimes \alpha_X)d_F) \\
&= (e_G \otimes \id_{G(X)})(\id_{G(X)} \otimes d_G)= \id_{G(X)}.
\end{align*}
Hence $\ldual{\alpha_Y}$ is inverse to $\alpha_X$.
\end{proof}

\subsection{Antipodes}\label{sect-antipodes}
Let $(T,\mu,\eta)$ be a bimonad on a monoidal category $\cc$.

If $\cc$ is left autonomous, then a \emph{left antipode for $T$} is a functorial  morphism $s^l\co T \ldual{?} T \to
\ldual{?}$, that is, a functorial  family $s^l=\{s^l_X\co T(\ldual{T(X)}) \to \ldual{X}\}_{X \in \Ob(\cc)}$ of morphisms
in $\cc$, satisfying:
\begin{align}
& T_0 T(\lev_X)T(\ldual{\eta_X} \otimes \id_X)=\lev_{T(X)}(s^l_{T(X)}T(\ldual{\mu}_X) \otimes
\id_{T(X)})T_2(\ldual{T(X)},X); \label{lant1} \\
& (\eta_X \otimes \id_{\ldual{X}})\lcoev_X T_0=(\mu_X \otimes s^l_X) T_2(T(X),\ldual{T(X)})T(\lcoev_{T(X)}); \label{lant2}
\end{align}
for all object $X$ of $\cc$.

If $\cc$ is right autonomous, then a \emph{right antipode for $T$} is a functorial  morphism $s^r\co T \rdual{?} T \to
\rdual{?}$, that is, a functorial  family $s^r=\{s^r_X\co T(\rdual{T(X)}) \to \rdual X\}_{X \in \Ob(\cc)}$ of morphisms
in $\cc$ satisfying:
\begin{align}
& T_0 T(\rev_X)T(\id_X \otimes \eta_X^\vee)=\rev_{T(X)}(\id_{T(X)}
\otimes s^r_{T(X)}T(\mu_X^\vee))T_2(X,\rdual{T(X)}); \label{rant1}\\
& (\id_{X^\vee}\otimes \eta_X )\rcoev_X T_0=(s^r_X \otimes\mu_X) T_2(\rdual{T(X)},T(X))T(\rcoev_{T(X)});  \label{rant2}
\end{align}
for all object $X$ of $\cc$.
\begin{rem}
This apparently complicated definition is justified by Theorem \ref{dual-ant}. The notion of left and right antipodes
generalize the classical notion of an antipode and its inverse for a bialgebra. For details, see Example~\ref{ExHopf} below.
\end{rem}

\begin{rem}\label{remopp}
Let $T$ be a bimonad on a left autonomous category $\cc$,  endowed with a left antipode $s^l$. Then $s^l$ is a right
antipode for the bimonad $T^\opp$  on the right autonomous category $\cc^{\otimes \opp}$ (as defined in Remark
\ref{mon-opp}). Likewise, if $T$ is a bimonad on a right autonomous category $\cc$ endowed with a right antipode $s^r$, then
$s^r$ is a left antipode for $T^\opp$.
\end{rem}

The next theorem translates the fact that the antipode of a (classical) Hopf algebra is an anti-homomorphism of bialgebras.
\begin{thm}\label{deadsiths}
Let $T$ be a bimonad on a monoidal category $\cc$. If $s^l$ is a left antipode of $T$ (assuming $\cc$ is left autonomous),
then we have:
\begin{align}
 &s^l_X \mu_{\ldual{T(X)}}=s^l_X T(s^l_{T(X)}) T^2(\ldual{\mu_X}); \label{lsithm}\\
 &s^l_X\eta_{\ldual{T(X)}}=\ldual{\eta_X};\label{lsitha}\\
 &s^l_{X\otimes Y}T(\ldual{T_2(X,Y)})=(s^l_Y \otimes s^l_X)T_2(\ldual{T(Y)},\ldual{T(X)});\label{lantcomult}\\
 &s^l_\un T(\ldual{T_0})=T_0.\label{lantcounit}
\end{align}
Likewise, if  $s^r$ is a  right antipode of $T$ (assuming $\cc$ is right autonomous), then we have:
\begin{align}
 &s^r_X \mu_{\rdual{T(X)}}=s^r_X T(s^r_{T(X)}) T^2(\rdual{\mu_X}); \label{rsithm}\\
 &s^r_X \eta_{\rdual{T(X)}}=\rdual{\eta_X};\label{rsitha}\\
 &s^r_{X \otimes Y}T(\rdual{T_2(X,Y)})=(s^r_Y \otimes s^r_X)T_2(\rdual{T(Y)},\rdual{T(X)});\label{rantcomult}\\
 &s^r_\un T(\rdual{T_0})=T_0.\label{rantcounit}
\end{align}
\end{thm}
\begin{proof}
The `right part' can be deduced from the `left part' by  Remark~\ref{remopp}. We  prove here  the `multiplicative'
assertions \eqref{lsithm} and \eqref{lsitha}. The `comultiplicative' assertions \eqref{lantcomult}, \eqref{lantcounit} (and
\eqref{rantcomult}, \eqref{rantcounit}), which are stated here for convenience, will be proved in Section~\ref{at-last}.
Note that we will not use these assertions until then!

Assume $\cc$ is left autonomous and $T$ has a left antipode~$s^l$. Let us show \eqref{lsithm}. Fix an object $X$ of $\cc$.
Setting $\lhs_X=s^l_X \mu_{\ldual{T(X)}}$ and $\rhs_X=s^l_X T(s^l_{T(X)}) T^2(\ldual{\mu_X})$, we must prove
$\lhs_X=\rhs_X$. Recall (see Lemma~\ref{compo}(b)) that $T^2$ is a comonoidal functor. Define $\mu^{(2)}_X : T^3(X) \to
T(X)$ and $D_X\co T^2(\un) \to T^3(X) \otimes T^2(\ldual{T(X)})$ by
\begin{equation*}
\mu^{(2)}_X=\mu_X T(\mu_X) \quad\text{and}\quad D_X=T^2_2(T(X),\ldual{T(X)})T^2(\lcoev_{T(X)}).
\end{equation*}
Firstly, we have:
\begin{equation}\label{eqsith1}
(\mu^{(2)}_X \otimes \lhs_X)D_X=(\mu^{(2)} \otimes \rhs_X)D_X.
\end{equation}
Indeed
\begin{align*}
(\mu^{(2)}_X &\otimes \lhs_X)D_X \\
&=(\mu_X \mu_{T(X)}\otimes s^l_X\mu_{\ldual{T(X)}})\, T^2_2(T(X),\ldual{T(X)})\,
T^2(\lcoev_{T(X)})  \quad \text{by \eqref{mon-ass}}\\
&=(\mu_X \otimes s^l_X) \, T_2(T(X),\ldual{T(X)}) \,\mu_{T(X) \otimes \ldual{T(X)}}\,
T^2(\lcoev_{T(X)}) \quad \text{by \eqref{bimonad1}}\\
&=(\mu_X \otimes s^l_X) \, T_2(T(X),\ldual{T(X)}) \, T(\lcoev_{T(X)})\, \mu_\un\\
&= (\eta_X \otimes \id_{\ldual{X}})\lcoev_X T_0 \mu_\un \quad \text{by \eqref{lant2}}\\
&= (\eta_X \otimes \id_{\ldual{X}})\lcoev_X T_0 T(T_0) \quad \text{by \eqref{bimonad2}}\\
&=(\mu_X \otimes s^l_X) \, T_2(T(X),\ldual{T(X)}) \, T(\lcoev_{T(X)}T_0) \quad \text{by \eqref{lant2}}\\
& = (\mu_X \otimes s^l_X)\, T_2(T(X),\ldual{T(X)}) \, T\bigl((\mu_X\eta_{T(X)} \otimes
 \id_{\ldual{T(X)}})\lcoev_{T(X)}T_0 \bigr) \quad \text{by \eqref{mon-u}}\\
 &= (\mu_X \otimes s^l_X)\, T_2(T(X),\ldual{T(X)} ) \\
 & \phantom{XXXX} T\bigl((\mu_X\mu_{T(X)} \otimes s^l_{T(X)})\,  T_2(T^2(X),\ldual{T^2(X)}) \, T(\lcoev_{T^2(X)}) \bigr)
 \quad \text{by \eqref{lant2}}\\
 &= (\mu_X \otimes s^l_X)\, T_2(T(X),\ldual{T(X)})\\
 & \phantom{XXXX} T\bigl ((\mu_X T(\mu_{X}) \otimes s^l_{T(X)})\, T_2(T^2(X),\ldual{T^2(X)}) \, T(\lcoev_{T^2(X)}) \bigr) \quad \text{by \eqref{mon-ass}}\\
 &= \bigl(\mu_X T(\mu_X) \otimes s^l_X T(s^l_{T(X)})T^2(\ldual{\mu_X})\bigr)\, T_2(T^2(X),T(\ldual{T(X)}) )\\
 & \phantom{XXXX} T(T_2(T(X),\ldual{T(X)})) \, T^2(\lcoev_{T(X)})\\
 &= (\mu^{(2)}_X \otimes \rhs_X)D_X.
\end{align*}
Secondly, setting:
\begin{equation*}
\nu_X=s^l_{T(X)} T(\ldual{\mu_X}),\quad \nu^{(2)}_X=\nu_X T(\nu_X),\quad\text{and}\quad E_X=\lev_{T(X)}(\nu^{(2)}_X \otimes
\mu^{(2)}_X),
\end{equation*}
we have
\begin{equation}\label{eqsith2}
(E_X\otimes \id_{T^2(\ldual{T(X)})})(\id_{T^2(\ldual{T(X)})} \otimes D_X)T^2_2(\ldual{T(X)},\un)=\id_{T^2(\ldual{T(X)})}.
\end{equation}
Indeed, on one hand, we have:
\begin{align*}
(\id_{T^2(\ldual{T(X)})} & \otimes D_X)T^2_2(\ldual{T(X)},\un)\\
&=(\id_{T^2(\ldual{T(X)})} \otimes T^2_2(T(X),\ldual{T(X)}))\\
&\phantom{XXX} T^2_2(\ldual{T(X)},T(X)\otimes \ldual{T(X)})\,
T^2(\id_{\ldual{T(X)}} \otimes \lcoev_{T(X)})\\
&=(T^2_2(\ldual{T(X)},T(X)) \otimes \id_{T^3(X)})\\
& \phantom{XXX} T^2_2(\ldual{T(X)} \otimes T(X),\ldual{T(X)})\,T^2(\id_{\ldual{T(X)}}\otimes \lcoev_{T(X)}) \quad \text{by
\eqref{comoneq1}.}
\end{align*}
On the other hand, from \eqref{lant1} and  \eqref{mon-u}, we obtain:
\begin{equation*}
\lev_{T(X)}(\nu_X \otimes \mu_X)T_2(\ldual{T(X)} \otimes T(X))=T_0 T(\lev_{T(X)})
\end{equation*}
and so, using this twice,
\begin{align*}
E_X & T^2_2(\ldual{T(X)},T(X))\\
&= \lev_{T(X)}(\nu_X \otimes \mu_X )T_2(\ldual{T(X)},T(X))T(\nu_X\otimes \mu_X)\,T(T_2(\ldual{T(X)},T(X)))\\
&= T_0 T(\lev_{T_X})T(\nu_X\otimes \mu_X)\,T(T_2(\ldual{T(X)},T(X)))=T^2_0 T^2(\lev_{T(X)}).
\end{align*}
Hence
\begin{align*}
(E_X \otimes \id&_{T^2(\ldual{T(X)})}) (\id_{T^2(\ldual{T(X)})} \otimes D_X)T^2_2(\ldual{T(X)},\un)\\
&=(T^2_0 T^2(\lev_{T(X)}) \otimes \id_{T^2(\ldual{T(X)})}) \\
& \phantom{XXXXX} T^2_2(\ldual{T(X)} \otimes T(X),\ldual{T(X)})
T^2(\id_{\ldual{T(X)}}\otimes \lcoev_{T(X)})\\
&=(T^2_0 \otimes \id_{T^2(\ldual{T(X)})}) T^2_2(\un,\ldual{T(X)})\\
& \phantom{XXXXX}
T^2((\lev_{T(X)} \otimes \id_{\ldual{T(X)}}) (\id_{\ldual{T(X)}}\otimes \lcoev_{T(X)}))\\
&=\id_{T^2(\ldual{T(X)})} \quad \text{by  \eqref{comoneq2} and \eqref{axiomdual},}
\end{align*}
that is \eqref{eqsith2}. Finally, we conclude:
\begin{align*}
\lhs_X & = (E_X\otimes \lhs_X)(\id_{T^2(\ldual{T(X)})}\otimes
D_X)T^2_2(\ldual{T(X)},\un) \quad \text{by \eqref{eqsith2}}\\
& = (E_X\otimes \rhs_X)(\id_{T^2(\ldual{T(X)})}\otimes
D_X)T^2_2(\ldual{T(X)},\un) \quad \text{by \eqref{eqsith1}}\\
& = \rhs_X \quad \text{by \eqref{eqsith2}.}
\end{align*}

Let us prove \eqref{lsitha}. For any object $X$ of $\cc$, we have:
\begin{align*}
(\id_{T(X)} \otimes \ldual{\eta_X})\lcoev_{T(X)}
  &= (\eta_X \otimes \id_{\ldual{X}})\lcoev_X T_0 \eta_\un \quad \text{by \eqref{bimonad4}}\\
  &= (\mu_X \otimes s^l_X) T_2(T(X),\ldual{T(X)})T(\lcoev_{T(X)})\eta_\un \quad \text{by \eqref{lant2}}\\
  &= (\mu_X \otimes s^l_X) T_2(T(X),\ldual{T(X)})\eta_{T(X) \otimes \ldual{T(X)}}\lcoev_{T(X)} \\
  &= (\mu_X \eta_{T(X)}\otimes s^l_X \eta_{\ldual{T(X)}}) \lcoev_{T(X)} \quad \text{by \eqref{bimonad3}}\\
  &= (\id_{T(X)}\otimes s^l_X \eta_{\ldual{T(X)}}) \lcoev_{T(X)} \quad \text{by \eqref{mon-u}.}
\end{align*}
Hence $s^l_X\eta_{\ldual{T(X)}}=\ldual{\eta_X}$ by \eqref{axiomdual}.
\end{proof}

The following theorem relates the existence of a left (resp.\@ right) antipode for a bimonad to the existence of left
(resp.\@ right) duals for the category of modules over the bimonad.
\begin{thm}\label{dual-ant}
Let $T$ be a bimonad on a monoidal category $\cc$.
\begin{enumerate}
  \renewcommand{\labelenumi}{{\rm (\alph{enumi})}}
  \item Assume $\cc$ is left autonomous. Then $T$ has a left antipode $s^l$ if and only if the category $T\ti\cc$ of $T$-modules is
left autonomous. Moreover, if a left antipode exists, then it is unique. In terms of a left antipode $s^l$, left duals in
$T\ti\cc$ are given by:
\begin{equation*}
\ldual{(M,r)}=(\ldual{M}, s^l_M T(\ldual{r})),\quad \lev_{(M,r)}=\lev_M, \quad\lcoev_{(M,r)}=\lcoev_M.
\end{equation*}
  \item Assume $\cc$ is right autonomous. Then $T$ has a right antipode $s^r$ if and only if the category $T\ti\cc$ of $T$-modules is
right autonomous. Moreover, if a right antipode exists, then it is unique. In terms of a right antipode $s^r$, right duals
in $T\ti\cc$ are given by:
\begin{equation*}
\rdual{(M,r)}=(\rdual{M}, s^r_M T(\rdual{r})),\quad \rev_{(M,r)}=\rev_M, \quad\rcoev_{(M,r)}=\rcoev_M.
\end{equation*}
\end{enumerate}
\end{thm}
We prove  Theorem~\ref{dual-ant} in the next section.

In Section~\ref{sect-perspec} (see Theorem~\ref{thm-gen-ex}), we show that any strong monoidal functor $U\co \dd \to \cc$
between monoidal categories which admits a left adjoint $F$ defines a bimonad $T=UF$ on $\cc$. Furthermore, if $\cc$ and
$\dd$ are left (resp.\@ right) autonomous, then the bimonad $T$ admits a left (resp.\@ right) antipode. This generalizes the
`if' assertions of Theorem \ref{dual-ant}.

\subsection{Proof of Theorem~\ref{dual-ant}}
We first establish the following lemma:
\begin{lem}\label{prefer}
Let $T$ be a bimonad on a left autonomous category $\cc$ and $(M,r)$ be a $T$-module. If $(M,r)$ has a left dual, then it
has a unique left dual of the form $((\ldual{M},r'),\lev_M, \lcoev_M)$. We will call it the \emph{preferred} left dual of
$(M,r)$.
\end{lem}
\begin{proof}
Assume  $(M,r)$ has a left dual $((N,\rho),e,d)$. The forgetful functor $U_T\co T\ti\cc \to \cc$ being strict monoidal,
there is a unique isomorphism $u \co N \to \ldual{M}$ such that  $e=\lev_M(u \otimes \id_M)$ and $d=(\id_M \otimes
u^{-1})\lcoev_M$.  Define $r'\co T (\ldual{M}) \to \ldual{M}$ by $r'= u \rho T(u^{-1})$. Then clearly $((\ldual{M},r'),
\lev_M,\lcoev_M)$ is a left dual of $(M,r)$. Now if we have another left dual of this form, say $((\ldual{M},r''),
\lev_M,\lcoev_M)$, then we have an isomorphism $v \co (\ldual{M},r') \to (\ldual{M},r'')$ such that  $\lev_M=\lev_M(v
\otimes \id_M)$,  and hence $v=\id_{\ldual{M}}$  and $r''=r'$.
\end{proof}

Let us now prove Theorem~\ref{dual-ant}. Part (b) is just a re-writing of Part (a) applied to the bimonad $T^\opp$
 (see Remark~\ref{remopp}). Let us show Part (a). Let $T$ be a bimonad on a left autonomous category $\cc$.
Recall that each object $X$ of $\cc$ has a left dual $(\ldual{X},\lev_X, \lcoev_X)$. By Lemma \ref{key-lemma} (contravariant
case) we have a canonical bijection
\begin{equation*}
?^\sharp\co \Nat(T\ldual{?}T^\opp, \ldual{?}) \to \Nat(T\ldual{?} U_T^\opp, \ldual{?}U_T^\opp), \; f \mapsto f^\sharp.
\end{equation*}
Denote its inverse by $?^\flat$. Recall that they are given by $f^\sharp_{(M,r)}=f_M T(\ldual{r})$ for any $T$-module
$(M,r)$ and $g^\flat= \ldual{\eta_X} g_{(T(X),\mu_X)}$ for any object $X$ of $\cc$.

Assume $T\ti\cc$ is left autonomous. For any $T$-module $(M,r)$, let $\delta_{(M,r)} \co T(\ldual{M}) \to \ldual{M}$ be
defined by the condition that $((\ldual{M},\delta_{(M,r)}),\lev_M,\lcoev_M)$ is the preferred left dual of $(M,r)$. By
Lemma~\ref{prefer}, this determines $\delta$ uniquely. Note that if $f \co M \to N$ is a $T$-linear morphism between two
$T$-modules $(M,r)$ and $(N,s)$, then $\ldual{f}\co \ldual{N} \to \ldual{M}$ is $T$-linear too. Hence $\delta$ satisfies the
following functoriality property:  $\ldual{f} \delta_{(N,s)} = \delta_{(M,r)} T(\ldual{f})$. Synthetically,
$\delta\in\Nat(T\ldual{?} U_T^\opp, \ldual{?}U_T^\opp)$. In particular $\delta$ leads to a morphism of functors $s^l =
\delta^\flat\in\Nat(T\ldual{?}T^\opp, \ldual{?})$ defined, for any object $X$ of $\cc$, by $s^l_X=\ldual{\eta_X}
\delta_{(T(X),\mu_X)}$. Conversely, $\delta$ can be recovered from $s^l$ by $\delta=(s^l)^\sharp$, that is,
$\delta_{(M,r)}=s^l_M T(\ldual{r})$ for any $T$-module $(M,r)$.

Now let $\delta \in \Nat(T\ldual{?} U_T^\opp, \ldual{?}U_T^\opp)$ and set $s^l= \delta^\flat\in \Nat(T\ldual{?}T^\opp,
\ldual{?})$. We have the following equivalences:
\begin{enumerate}
\renewcommand{\labelenumi}{{\rm (\roman{enumi})}}
\item  The pair  $(\ldual{M},\delta_{(M,r)})$ is a $T$-module for all $T$-module $(M,r)$ if and only if $s^l$ satisfies
 \eqref{lsithm} and \eqref{lsitha};
\end{enumerate}
and, assuming the equivalent assertions of (i),
\begin{enumerate}
\renewcommand{\labelenumi}{{\rm (\roman{enumi})}}
\setcounter{enumi}{1}
\item  The evaluation  $\lev_M$ is $T$-linear for all $T$-module $(M,r)$ if and only if $s^l$ satisfies \eqref{lant1};
\item  The coevaluation  $\lcoev_M$ is $T$-linear  for all $T$-module $(M,r)$ if and only if $s^l$ satisfies \eqref{lant2}.
\end{enumerate}

Before verifying these equivalences, let us show that they suffice to prove the theorem. If $s^l$ is a left antipode then,
using only the `multiplicative' part of Theorem \ref{deadsiths}, which we have already proved, and setting
$\delta={s^l}^\sharp$, we see by (i)-(iii) that $((\ldual{M},\delta_{(M,r)}), \lev_M,\lcoev_M)$ is a left dual of $(M,r)$ in
$T\ti \cc$, and so $T\ti\cc$ is left autonomous. Conversely, if $T\ti\cc$ is left autonomous, then there exists $\delta$
such that $((\ldual{M},\delta_{(M,r)}), \lev_M,\lcoev_M)$ is a duality in $T\ti \cc$ and so, by (ii)-(iii),
$s^l=\delta^\flat$ satisfies the axioms of a left antipode. By Lemma~\ref{prefer}, such a $\delta$ is unique. Hence the
uniqueness of a left antipode, since the correspondence $\delta \leftrightarrow s^l$ is bijective.

Now let us show (i). Recall that $(\ldual{M},\delta_{(M,r)})$ is a $T$-module if and only if both identities $\delta_{(M,r)}
\mu_{\ldual{M}}= \delta_{(M,r)} T(\delta_{(M,r)})$ and $\delta_{(M,r)} \eta_{\ldual{M}}= \id_{\ldual{M}}$ hold. Replacing
$\delta_{(M,r)}$ by $s^l_M T(\ldual{r})$ in the first identity, we get:
\begin{equation*}
s^l_M T(\ldual{r})\mu_{\ldual{M}}=s^l_M T(\ldual{r})T(s^l_M)T^2(\ldual{r}).
\end{equation*}
The left-hand side may be rewritten as $s^l_M \mu_{\ldual{T(M)}} T^2(\ldual{r})$. The right-hand side may be rewritten as
 $s^l_M T(s^l_{T(M)})T^2(\ldual{T(r)} \ldual{r}) = s^l_M T(s^l_{T(M)})T^2(\ldual{\mu_M} \ldual{r})$.  Therefore we
finally get:
\begin{equation*}
s^l_M \mu_{\ldual{T(M)}} T^2(\ldual{r})=s_M T(s^l_{T(M)})T^2(\ldual{\mu_M}) T^2(\ldual{r}),
\end{equation*}
which is equivalent to \eqref{lsithm} by Lemma \ref{key-lemma}. Likewise, the second identity is equivalent to
\eqref{lsitha} by a straightforward application of Lemma \ref{key-lemma}.

Let us show (ii). Recall that $\lev_M$ is $T$-linear if an only if we have $T_0 T(\lev_M)= \lev_M(\delta_{(M,r)} \otimes
r)T_2(\ldual{M},M)$. Replacing $\delta_{(M,r)}$ by $s^l_M T(\ldual{r})$, we get:
\begin{equation}\label{eqdemthmantip}
T_0 T(\lev_M)= \lev_M(s^l_M T(\ldual{r})\otimes r)T_2(\ldual{M},M).
\end{equation}
Now, we have:
\begin{align*}
\lev_M(s^l_M T(&\ldual{r}) \otimes r)T_2(\ldual{M},M)\\
 & = \lev_{TM}(\ldual{r} s^l_M T(\ldual{r})\otimes \id_{T(M)})T_2(\ldual{M},M)\\
 & = \lev_{TM}( s^l_{T(M)} T(\ldual{T(r)})T(\ldual{r})\otimes \id_{T(M)})T_2(\ldual{M},M) \\
 & = \lev_{TM}( s^l_{T(M)} T(\ldual{\mu_M})T(\ldual{r})\otimes \id_{T(M)})T_2(\ldual{M},M)\quad \text{by \eqref{Tmoddef}} \\
 & = \lev_{TM}( s^l_{T(M)} T(\ldual{\mu_M}))\otimes \id_{T(M)})T_2(\ldual{T(M)},M)T(\ldual{r}\otimes \id_M).
\end{align*}
On the other hand, $T_0 T(\lev_M)= T_0 T(\lev_M)T(\ldual{\eta_M} \otimes \id_M)T(\ldual{r}\otimes \id_M)$  by
\eqref{Tmoddef}.  Therefore, by Lemma \ref{key-lemma} and duality, \eqref{eqdemthmantip} is equivalent to Axiom
\eqref{lant1}.

Finally, let us show (iii). Recall that $\lcoev_M$ is $T$-linear if and only if we have $\lcoev_M T_0= (r \otimes
\delta_{(M,r)}) T_2(M, \ldual{M}) T(\lcoev_M)$. By a computation similar to that of the proof of (ii), this is equivalent
to:
\begin{align*}
(r \otimes \id_{\ldual{M}})&(\eta_M \otimes \id_{\ldual{M}})\lcoev_M T_0 \\
& =(r \otimes \id_{\ldual{M}})(\mu_M \otimes s_M)T_2(T(M), \ldual{T(M)})T(\lcoev_{T(M)}),
\end{align*}
and so, by Lemma \ref{key-lemma} and duality, to Axiom \eqref{lant2}.  This completes the proof of Theorem~\ref{dual-ant}.

\subsection{End of the proof of Theorem \ref{deadsiths}}\label{at-last}
We still have to prove assertions \eqref{lantcomult} and \eqref{lantcounit} of Theorem \ref{deadsiths} (from which
\eqref{rantcomult} and \eqref{rantcounit} can be deduced via Remark~\ref{remopp}).

To show \eqref{lantcomult}, let $(M,r)$ and $(N,s)$ be two $T$-modules. Recall that
\begin{equation*}
\ldual{(N,s)} \otimes \ldual{(M,r)}=\bigl (\ldual{N} \otimes \ldual{M}, (s^l_N T(\ldual{s}) \otimes s^l_M
T(\ldual{r}))T_2(\ldual{N},\ldual{M}) \bigr )
\end{equation*}
and
\begin{equation*}
\ldual{\bigl ((M,r) \otimes (N,s) \bigr )}=\bigl (\ldual{(M \otimes N)}, s^l_{M \otimes N}T(\ldual{T_2(M,N)}\ldual{(r
\otimes s)}) \bigr )
\end{equation*}
are canonically isomorphic via the isomorphism  $\ldual{?}_2(M,N)\co\ldual{N} \otimes \ldual{M} \to \ldual{(M \otimes
N)}$.  By Lemma~\ref{prefer}, we get (up to suitable identifications):
\begin{equation*}
(s^l_N \otimes s^l_M )T_2(\ldual{T(N)},\ldual{T(M)})T(\ldual{(r \otimes s)}) =s^l_{M \otimes
N}T(\ldual{T_2(M,N)})T(\ldual{(r \otimes s)}).
\end{equation*}
Hence \eqref{lantcomult} by applying Lemma \ref{key-lemma}.

Finally, via the isomorphism $\ldual{?}_0\co \un \to \ldual{\un}$, the $T$-modules $(\un, T_0)$ and $\ldual(\un,
T_0)=(\ldual{\un}, s^l_\un T(\ldual{T_0}))$ are isomorphic. Hence $s^l_\un T(\ldual{T_0})=T_0$, that is, \eqref{lantcounit}.

\subsection{Hopf monads} A \emph{left (resp.\@ right) Hopf monad} is a bimonad on a left (resp.\@ right) autonomous
category which has a left (resp.\@ right) antipode.

A \emph{Hopf monad} is a bimonad on an autonomous category which has a left antipode and a right antipode. In particular, by
Theorem~\ref{dual-ant}, the category of modules  over  a Hopf monad is autonomous.

\begin{exa} \label{ExHopf}
Let $\cc$ be a braided autonomous category with braiding by $\tau$. Let $A$ be a bialgebra in $\cc$, with product $m$, unit
$u$, coproduct $\Delta$, and counit~$\varepsilon$. Consider the bimonad $A\otimes ?$ (see Example~\ref{Ex2}). Firstly, let
$S \co A \to A$ be a morphism in $\cc$ and define:
\begin{equation*}
s^l_X=( \lev_A \, \tau_{\ldual{A},A} \otimes \id_{\ldual{X}})(S \otimes \tau_{\ldual{A},\ldual{X}}^{-1})\co A \otimes
\ldual{X} \otimes \ldual{A} \to \ldual{X}.
\end{equation*}
Then $s^l$ is a left antipode for the bimonad $A \otimes ?$ if and only if $S$ is an antipode of the bialgebra $A$, that is,
if and only if $S$ satisfies:
\begin{equation*}
m(S \otimes \id_A)\Delta= u \varepsilon = m(\id_A \otimes S)\Delta.
\end{equation*}
Secondly, let $S' \co A \to A$ be another morphism in~$\cc$ and define:
\begin{equation*}
s^r_X=(\rev_A \otimes \id_{\rdual{X}})(S' \otimes \tau_{\rdual{A},\rdual{X}}^{-1})\co A \otimes \rdual{X} \otimes \rdual{A}
\to \rdual{X}.
\end{equation*}
Then $s^r$ is a right antipode for the bimonad $A \otimes ?$ if and only if $S'$ is an `inverse of the antipode', that is,
setting $m^\opp= m \tau^{-1}_{A,A}$, if and only if $S'$ satisfies:
\begin{equation*}
m^\opp(S' \otimes \id_A)\Delta= u \varepsilon = m^\opp(\id_A \otimes S')\Delta.
\end{equation*}
Thus $A \otimes ?$ is a Hopf monad  if and only if $A$ is a Hopf algebra in $\cc$ with invertible antipode. Similarly, a
right antipode for the bimonad $? \otimes A$ corresponds with an antipode for the bialgebra $A$, and a left antipode for $?
\otimes A$ corresponds with an `inverse of the antipode' for $A$. In particular, any \emph{finite-dimensional} Hopf algebra
$H$ over a field $\kk$ yields Hopf monads $H \otimes_\kk ?$ and $? \otimes_\kk H$ on the category $\vect(\kk)$ of
finite-dimensional $\kk$-vector spaces.
\end{exa}

\begin{prop} \label{anti-inv}
Let $T$ be a Hopf monad on an autonomous category $\cc$. Then its left antipode $s^l$ and its right antipode $s^r$ are
`inverse' to each other in the sense:
\begin{equation*}
\id_{T(X)}=s^r_{\ldual{T(X)}} T(\rdual{(s_X^l)})=s^l_{\rdual{T(X)}} T(\ldual{(s_X^r)})
\end{equation*}
for any object $X$ of $\cc$ (up to the canonical isomorphisms of Remark~\ref{qieq}).
\end{prop}
\begin{proof}
Let $(M,r)$ be a $T$-module. We have $\ldual{(M,r)}=(\ldual{M}, s^l_M T(\ldual{r}))$ and so $\rdual{(\ldual{(M,r)})} =
(\rdual{(\ldual{M})}, s^r_{\ldual{M}} T(\rdual{T(\ldual r)}) T(\rdual{(s^l_M)}))$. Via the canonical isomorphism of
Remark~\ref{qieq}, we have $r=s^r_{\ldual{M}} T(\rdual{T(\ldual r)}) T(\rdual{(s^l_M)})=r s^r_{\ldual{T(M)}}
T(\rdual{(s^l_M)})$. So, by Lemma \ref{key-lemma}, we have $\id_{T(X)}=s^r_{\ldual{T(X)}} T(\rdual{(s^l_X)})$. The second
identity is obtained by replacing $T$ with $T^\opp$.
\end{proof}

\begin{rem}\label{remTduals}
Let $T$ be a Hopf monad on an autonomous category $\cc$. Let $T^!$ the endofunctor of $\cc$ defined as in
Remark~\ref{rem-excla}. Then Proposition \ref{anti-inv} says that $T^!$ is right-adjoint to $T$, with adjunction morphisms
$e \co TT^! \to 1_\cc$ and $h \co 1_\cc \to T^!T$ given by $e_X = s^r_{\ldual{X}}\co T (T^! (X)) \to \rdual{(\ldual{X})}
\simeq X$ and $h_X =\rdual{(s^l_X)}\co X \simeq \rdual{(\ldual{X})} \to T^! (T(X))$. Likewise $\lexcla{T}$ is a
right-adjoint to $T$. An interesting consequence is that a Hopf monad always commutes with direct limits.
\end{rem}

\subsection{Morphisms of Hopf monads}\label{sect-morphhopfmon}
A \emph{morphism of Hopf monads} on an autonomous category is a morphism of their underlying bimonads (see
Section~\ref{sect-morphbimon}).

\begin{lem}\label{lemmorphhopfmon}
A morphism $f\co T \to T'$ of Hopf monads preserves the antipodes. More precisely, if $T$ has a left antipode $s^l$ and $T'$
has a left antipode $s'^{l}$, then $s^l_X T(\ldual{f_X})=s'^l_X f_{\ldual{T'(X)}}$ for any object $X$ of $\cc$, and
similarly for right antipodes.
\end{lem}
\begin{proof}
Let $f^*\co T'\ti \cc \to T\ti \cc$ be the strict monoidal functor induced by $f$. Recall it is given by $f^*(M,r) =(M,r
f_M)$.  Let $(M,r)$ be a $T'$-module. Since $f^*$ is monoidal strict,  $f^*(\ldual(M,r))=(\ldual{M},s'^l_M
T'(\ldual{r})f_{\ldual{M}})$ is a left dual of $(M,r)$ and so canonically isomorphic to $\ldual{f^*(M,r)}=(\ldual{M}, s^l_M
T(\ldual{f_M}\ldual{r}))$. Therefore $s^l_M T(\ldual{f_M})T(\ldual{r})=s'^l_M f_{\ldual{T'(M)}}T(\ldual{r})$. Hence, by
Lemma \ref{key-lemma}, we get $s^l_X T(\ldual{f_X})=s'^l_X f_{\ldual{T'(X)}}$ for any object $X$ of $\cc$.
\end{proof}

\subsection{Convolution product and antipodes} \label{sect-convoantip}
Let $T$ be a Hopf monad on an autonomous category~$\cc$. Let $?^\sharp \co \Nat(1_\cc,T) \to \Nat(U_T,U_T)$ be the
isomorphism of Lemma~\ref{key-lemma}, with inverse $?^\flat$, and let $\lexcla{?},\,\rexcla{?}\co \End(\cc)^\opp \to
\End(\cc)$ be the functors of Remark~\ref{rem-excla}. Define two maps:
\begin{equation}\label{defSL}
S\co \left \{
\begin{array}{ccl}
  \Nat(1_\cc,T)& \to & \Nat(1_\cc,T) \\
  f & \mapsto & S(f)=(\lexcla{(f^\sharp)})^\flat
\end{array}\right .
\end{equation}
and
\begin{equation}\label{defSR}
S^{-1}\co \left \{
\begin{array}{ccl}
  \Nat(1_\cc,T)& \to & \Nat(1_\cc,T) \\
  f & \mapsto & S(f)=(\rexcla{(f^\sharp)})^\flat
\end{array}\right .
\end{equation}
Explicitly, using Theorem~\ref{dual-ant} and Lemma~\ref{key-lemma}, we have:
\begin{equation*}
S(f)_X=\rdual{(s^l_X f_{\ldual{T(X)}})} \quad \text{and} \quad S^{-1}(f)_X=\ldual{(s^r_Xf_{\rdual{T(X)}})}
\end{equation*}
for all object $X$ of $\cc$ (up to the canonical isomorphisms of Remark~\ref{qieq}), where $s^l$ and $s^r$ are the left and
right antipodes of $T$ respectively.

\begin{lem}\label{lemsrsl}
Let $T$ be a Hopf monad on an autonomous category $\cc$. Then the map $S$ is an anti-automorphism of the monoid
$(\Nat(1_\cc,T),*,\eta)$, and $S^{-1}$ is its inverse.
\end{lem}
\begin{proof}
Since the convolution product $*$ corresponds to composition of endomorphisms of $U_T$, and since the functors $\lexcla{?},
\rexcla{?}\co \End(\cc) \to \End(\cc)^\opp$ are strong monoidal, the maps $S$ and $S^{-1}$ are anti-endomorphisms of
$\Nat(1_\cc,T)$. Since the functors $\lexcla{?}$ and $(\rexcla{?})^\opp$ are inverse to each other (up to the canonical
isomorphisms of Remark~\ref{qieq}), the maps $S$ and $S^{-1}$ are inverse to each other.
\end{proof}

\begin{exa}
For the Hopf monad $A\otimes ?$ (see Example~\ref{ExHopf}), where $A$ is a Hopf algebra in an autonomous braided category,
the maps $S$ and $S^{-1}$ are given by $S(f)=(S_A \otimes 1_\cc)f$ and $S^{-1}(f)=(S_A^{-1} \otimes 1_\cc)f$, where $S_A$ is
the antipode of $A$.
\end{exa}

\subsection{Grouplike elements}
A \emph{grouplike element} of a bimonad $T$ on a monoidal category $\cc$ is a functorial morphism $g \co 1_\cc \to T$
satisfying:
\begin{align}
&T_2(X,Y) g_{X\otimes Y}=g_X \otimes g_Y; \label{grl1} \\
&T_0 g_\un = \id_\un \label{grl2}.
\end{align}
We will denote by $G(T)$ the set of grouplike elements of $T$. Using \eqref{bimonad1}-\eqref{bimonad4}, we see that
$(G(T),*,\eta)$ is a monoid, where $*$ is the convolution product \eqref{prodconvo1T}.

\begin{lem}\label{gl-lem1}
Let $T$ be a bimonad on a monoidal category $\cc$. Via the canonical bijection $\Nat(1_\cc,T) \simeq \Nat(U_T,U_T)$ of
Lemma~\ref{key-lemma}, grouplike elements of~$T$ correspond exactly with monoidal endomorphisms of the strict monoidal
functor~$U_T$.
\end{lem}
\begin{proof}
Let $g \in \Nat(1_\cc,T)$. Then  $g^\sharp \in \Nat(U_T,U_T)$ is monoidal if and only if, for all $(M,r)$ and $(N,s)$ in
$T\ti\cc$, we have $(r \otimes s) T_2(X,Y)\,g_{X \otimes Y}=(r \otimes s) (g_X \otimes g_Y)$ and $T_0\,g_\un=\id_\un$, which
is equivalent to $g \in G(T)$ by Lemma \ref{key-lemma}.
\end{proof}

\begin{lem}\label{gl-lem2}
Let $T$ be a Hopf monad on an autonomous category $\cc$. Then  $(G(T),*,\eta)$ is a group. Moreover the inverse of $g \in
G(T)$ is  $g^{*-1}=S(g)=S^{-1}(g)$, with $S$ and $S^{-1}$ as in~\eqref{defSL} and \eqref{defSR}.
\end{lem}
\begin{proof}
Let $g \in G(T)$. By Lemma~\ref{isomono}, $g^\sharp \in \Nat(U_T,U_T)$ is a monoidal isomorphism with inverse
$\lexcla{(g^\sharp)}=\rexcla{(g^\sharp)}$. Hence, by Lemma~\ref{gl-lem1}, $g$ is invertible with inverse $S(g)=S^{-1}(g)$.
\end{proof}

\section{Hopf modules}\label{sect-Hmod}

In this section, we introduce Hopf modules and prove the fundamental theorem for Hopf modules over a Hopf monad.

\subsection{Comodules}
Let $C$ be a coalgebra in a monoidal category $\cc$, with coproduct $\Delta \co C \to C \otimes C$ and counit
$\varepsilon\co C \to \un$. Recall that a right $C$-comodule is a pair $(M,\rho)$, where $M$ is an object of $\cc$ and
$\rho\co M \to M \otimes C$ is a morphism in $\cc$, satisfying:
\begin{equation}\label{comddef}
(\rho \otimes \id_C)\rho=(\id_M \otimes \Delta)\rho \quad \text{and} \quad (\id_M \otimes \varepsilon)\rho=\id_M .
\end{equation}
A morphism $f\co (M,\rho) \to (N,\varrho)$ of right $C$-comodules is a morphism $f \co M \to N$ in~$\cc$ such that $\varrho
f=(f \otimes \id_C)\rho$. Thus the category of right $C$-comodules. Likewise one defines the category of left $C$-comodules.
\begin{lem}\label{lemcomod}
Let $C$ be a coalgebra in a monoidal category~$\cc$. If $(M,\rho)$ is a left $C$-comodule and $\cc$ is right autonomous,
then $\rdual{(M,\rho)}=(\rdual{M},\varrho^r)$ is a right $C$\ti comodule, where
\begin{equation*}
\varrho^r=(\id_{\rdual{M}\otimes C}\otimes \rev_{M})(\id_{\rdual{M}} \otimes \rho \otimes \id_{\rdual{M}}) (\rcoev_{M}
\otimes \id_{\rdual{M}}).
\end{equation*}
Moreover, this construction defines a contravariant functor form the category of left $C$-comodules to the category of right
$C$-comodules. Similarly, if $(M,\rho)$ is a right $C$-comodule and $\cc$ is left autonomous, then
$\ldual{(M,\rho)}=(\ldual{M},\varrho^l)$ is a left $C$-comodule, where
\begin{equation*}
\varrho^l= (\lev_M\otimes \id_{C\otimes \ldual{M}})(\id_{\ldual{M}} \otimes \rho \otimes \id_{\ldual{M}}) ( \id_{\ldual{M}}
\otimes \lcoev_{M})\,.
\end{equation*}
This construction is functorial too.
\end{lem}
\begin{proof}
Left to the reader.
\end{proof}

Let $T$ be a comonoidal endofunctor of a monoidal category $\cc$. By \eqref{comoneq1} and \eqref{comoneq2}, the object
$T(\un)$ is a coalgebra in $\cc$, with coproduct $T_2(\un,\un)$ and counit $T_0$. By a \emph{left} (resp.\@ \emph{right})
\emph{$T$-comodule}, we mean a left (resp.\@ right) $T(\un)$-comodule.

Note that if $T$ is a bimonad, then every object $X$ becomes a  left (resp.\@ right) $T$-comodule with \emph{trivial
coaction} given by $\eta_\un \otimes \id_X$ (resp.\@ $\id_X \otimes \eta_\un$).

\subsection{Hopf modules}
Let $T$ be a bimonad on a monoidal category $\cc$. The axioms of a bimonad ensure that $(T(\un),\mu_\un)$ is a coalgebra in
the category $T\ti \cc$ of $T$-modules, with coproduct $T_2(\un,\un)$ and counit $T_0$. A \emph{right Hopf $T$-module} is a
right $(T(\un),\mu_\un)$-comodule in $T\ti\cc$, that is, a triple $(M, r, \rho)$ such that $(M,r)$ is a $T$-module,
$(M,\rho)$ is a right $T$-comodule, and:
\begin{equation}\label{defHopfmod}
\rho r= (r \otimes \mu_\un) T_2(M,T(\un)) T(\rho).
\end{equation}
 A \emph{morphism of Hopf $T$-modules} between two right Hopf $T$-modules $(M, r, \rho)$ and $(N, s, \varrho)$ is a morphism of
$(T(\un),\mu_\un)$-comodules in $T\ti\cc$, that is, a morphism $f \co M \to N$ in $\cc$ such that
\begin{equation}
f r=s T(f) \text{\quad and \quad} (f \otimes \id_{T(\un)}) \rho= \varrho f.
\end{equation}

\begin{rem}\label{remisoHopf}
As is the classical case, any  morphism of Hopf $T$-modules  which is an isomorphism in $\cc$ is an isomorphism of Hopf
$T$-modules.
\end{rem}

Similarly, one can define the notion of \emph{left Hopf $T$-module}, which is a right Hopf $T$-module for the bimonad
$T^\opp$ (see Remark~\ref{mon-opp}).

\begin{lem}\label{lemexHopf}
Let $T$ be a bimonad on a monoidal category $\cc$. If $(M,\rho)$ is a right $T$-comodule, then the triple $(T(M), \mu_M,
\varrho)$ is a right Hopf $T$-module, where $\varrho=(\id_{T(M)} \otimes \mu_\un) T_2(M,T(\un)) T(\rho)$. In particular
$(T(X), \mu_X, T_2(X,\un))$  is a right Hopf $T$-module for any object $X$ of $\cc$.
\end{lem}
\begin{proof}
Let $(M,\rho)$ be a right $T$-comodule.  Firstly, we have:
\begin{align*}
(\varrho  \otimes \id&_{T(\un)})\varrho
  =\bigl ((\id_{T(M)} \otimes \mu_\un) T_2(M,T(\un)) T(\rho) \otimes \mu_\un \bigr) T_2(M,T(\un)) T(\rho)\\
 &=\bigl ((\id_{T(M)} \otimes \mu_\un) T_2(M,T(\un)) \otimes \mu_\un \bigr) T_2(M \otimes T(\un),T(\un))
    T((\rho\otimes \id_{T(\un)})\rho)\\
 &=\bigl (\id_{T(M)} \otimes (\mu_\un \otimes \mu_\un) T_2(T(\un),T(\un)) \bigr) T_2(M, T(\un)\otimes T(\un))\\
 & \phantom{XXXXXXXXXXXXXXX} T((\id_M \otimes T_2(\un,\un))\rho) \quad \text{by \eqref{comoneq1} and \eqref{comddef}}\\
 &=\bigl (\id_{T(M)} \otimes (\mu_\un \otimes \mu_\un) T^2_2(\un,\un) \bigr) T_2(M, T(\un))
    T(\rho)\\
 &=(\id_{T(M)} \otimes T_2(\un,\un)\mu_\un  ) T_2(M, T(\un))
    T(\rho) \quad \text{by \eqref{bimonad1}}\\
 &=(\id_{T(M)} \otimes T_2(\un,\un)  ) \varrho
\end{align*}
 and
\begin{align*}
(\id_{T(M)} \otimes T_0) \varrho &= (\id_{T(M)} \otimes T_0\mu_\un) T_2(M,T(\un)) T(\rho) \\
&= (\id_{T(M)} \otimes T_0T(T_0)) T_2(M,T(\un)) T(\rho) \quad \text{by \eqref{bimonad2}}\\
&= (\id_{T(M)} \otimes T_0) T_2(M,\un) T((\id_M \otimes T_0)\rho)\\
&= \id_{T(M)} \quad \text{by \eqref{comoneq2} and \eqref{comddef},}
\end{align*}
so that $(T(M),\varrho)$ is a right $T$-comodule. Secondly,
\begin{align*}
 \varrho  \mu_M &= (\id_{T(M)} \otimes \mu_\un) T_2(M,T(\un)) T(\rho) \mu_M\\
&= (\id_{T(M)} \otimes \mu_\un) T_2(M,T(\un)) \mu_{M \otimes T(\un)} T^2(\rho) \\
&= (\mu_M \otimes \mu_\un  \mu_{T(\un)}  ) T_2(T(M),T^2(\un)) T(T_2(M,\un))T^2(\rho) \quad \text{by \eqref{bimonad1}} \\
&= (\mu_M \otimes \mu_\un T(\mu_\un)) T_2(T(M),T^2(\un)) T(T_2(M,\un)T(\rho)) \quad \text{by \eqref{mon-ass}} \\
&= (\mu_M \otimes \mu_\un) T_2(T(M),T(\un)) T(\varrho).
\end{align*}
Hence $(T(M), \mu_M, \varrho)$ is a right Hopf $T$-module. Now, for any object $X$  of $\cc$, the pair  $(X,\id_X
\otimes \eta_\un)$ is a right $T$-comodule, so that $(T(X), \mu_X, \varrho)$ is a right Hopf $T$-module, with
\begin{align*}\varrho=
(\id_{T(X)} \otimes \mu_\un) T_2 & (X,T(\un)) T(\id_X \otimes \eta_\un)\\ &=(\id_{T(X)} \otimes \mu_\un T(\eta_\un))
T_2(X,\un) = T_2(X,\un)
\end{align*}
by  \eqref{mon-u}, which completes the proof of Lemma~\ref{lemexHopf}.
\end{proof}

\begin{lem}\label{leftrighthopfmod}
Let $T$ be a right Hopf monad on a right autonomous category $\cc$. If $M$ is a left Hopf module, then $\rdual{M}$ is a
right Hopf $T$-module, with the structure of $T$-module defined in Theorem~\ref{dual-ant}(b) and the structure of right
$T$-comodule defined in Lemma~\ref{lemcomod}. This defines a contravariant functor from the category of left Hopf
$T$-modules to the category of right Hopf $T$-modules.
\end{lem}
\begin{proof}
This results from Lemma~\ref{lemcomod}. Indeed, recall that  $(T(\un),\mu_\un)$ is a coalgebra in $T\ti\cc$, with coproduct
$T_2(\un,\un)$ and counit $T_0$.  Let $(M,r,\rho)$ be a left Hopf $T$\trait module. This means that $((M,r),\rho)$ is a left
$(T(\un),\mu_\un)$\trait comodule,  so $(\rdual{(M,r)}, \varrho^l)$ is a right $(T(\un),\mu_\un)$\trait co\-mo\-du\-le,  in
the notations of Lemma~\ref{lemcomod}. In other words, $(M, s^l_M T(\ldual{r}), \varrho^l)$ is a right Hopf $T$-module. This
construction is functorial since morphisms of Hopf $T$-modules are nothing but morphisms of $(T(\un),\mu_\un)$\trait
comodules.
\end{proof}

\subsection{Coinvariants}
Let $\dd$ be a category and $f,g\co X \to Y$ be parallel morphisms in $\dd$. A morphism $i\co E \to X$ in $\dd$
\emph{equalizes} the pair $(f,g)$ if $f\,i=g\,i$. An \emph{equalizer} (also called \emph{difference kernel}) of the pair
$(f,g)$ is a morphism $i\co E \to X$ which equalizes the pair $(f,g)$ and which is universal for this property in the
following sense:  for any morphism $j\co F \to X$ in $\cc$ equalizing the pair $(f,g)$, there exits a unique morphism $p\co
F \to E$ in $\cc$ such that  $j=pi$.  We say that \emph{equalizers exist in $\dd$} if each pair of parallel morphisms in
$\dd$ admits an equalizer.

We say that a functor $F\co \dd \to \dd'$ \emph{preserves equalizers} if, whenever $i$ is an equalizer of a pair
$(f,g)$ of parallel morphisms in $\dd$,  then $F(i)$ is an equalizer of the pair $(F(f),F(g))$. Notice that a left exact
functor preserves equalizers.\\

Let $T$ be a bimonad on a monoidal category $\cc$. We say that a right $T$-comodule $(M, \rho)$ \emph{admits coinvariants}
if the pair of parallel morphisms $(\rho,\id_M \otimes \eta_\un)$ admits an equalizer:
\begin{equation*}
\xymatrix@1{N \ar[r]^i & M \ar@<2pt>[r]^-{\rho}\ar@<-2pt>[r]_-{\id_M \otimes \eta_\un} &M \otimes T(\un).}
\end{equation*}
If such is the case, $N$ is called the \emph{coinvariant part of $M$,} and is denoted $M^{\mathrm{co} T}$. In fact
$M^{\mathrm{co} T}$ is a right $T$-comodule (with trivial coaction) and $i\co (N,\id_N \otimes \eta_\un) \to (M,\rho)$ is a
morphism of $T$-comodules.

Similarly, one defines the \emph{coinvariant part} of a left $T$-comodule $(M,\rho)$ which is, when it exists, an equalizer
of the pair $(\rho,\id_\un \otimes \id_M)$.

If a right or left $T$-comodule $(M,\rho)$ admits a coinvariant part $i\co M^{\mathrm{co} T} \to M$, we say that $T$
\emph{preserves the coinvariant part of $(M,\rho)$} if  $T(i)$ is an equalizer of the pair $(T(\rho),T(\id_M \otimes
\eta_\un))$ or $(T(\rho),T(\eta_\un \otimes \id_M))$ respectively. Note this is the case when $T$ preserves equalizers.

We say that a right (resp.\@ left) Hopf $T$-module $(M,r,\rho)$ \emph {admits coinvariants} if the underlying right (resp.\@
left) $T$-comodule $(M,\rho)$ admits coinvariants. If such the case, the \emph{coinvariant part of $(M,r,\rho)$} is the
coinvariant part of $(M,\rho)$.

\subsection{Decomposition of Hopf modules}
In this section we show  that, under certain assumptions on equalizers, Hopf modules can be decomposed as in the classical
case.

\begin{thm}\label{Hopfmodthm}
Let $T$ be a right Hopf monad on a right autonomous category. Let  $(M,r,\rho)$ be a right Hopf $T$-module such that
$(M,\rho)$ admits a coinvariant part $i\co M^{\mathrm{co} T} \to M$ which is preserved by $T$. Then $$rT(i)\co (M,r,\rho)
\to \bigl (T(M^{\mathrm{co} T}), \mu_{M^{\mathrm{co}T}}, T_2(M^{\mathrm{co} T},\un)\bigr)$$ is an isomorphism of right Hopf
$T$-modules.
\end{thm}
\begin{proof}
See Section~\ref{proofmodhopfthm}.
\end{proof}

Recall that a functor $F\co \cc \to \dd$ is said to be \emph{conservative} if any  morphism $f$ in~$\cc$ such that $F(f)$ is
an isomorphism in $\dd$, is an isomorphism in $\cc$.

\begin{thm}\label{Hopfmodcor}
Let $T$ be a right Hopf monad on a right autonomous category~$\cc$.  Suppose  that  right Hopf $T$-modules admit
coinvariants which are preserved by $T$. Then the assignments
\begin{equation*}
X \mapsto \bigl(T(X),\mu_X,T_2(X,\un) \bigr)\,, \quad f \mapsto T(f)
\end{equation*}
define a functor from $\cc$ to the category of right Hopf $T$-modules, which is an equivalence of categories if and only if
$T$ is conservative.
\end{thm}

\begin{proof}
See Section~\ref{proofmodhopfcor}.
\end{proof}

\begin{rem}
For a left Hopf monad $T$ over a left autonomous category,  one may formulate a similar decomposition theorem for left Hopf
$T$-modules, which may be deduced from  Theorem~\ref{Hopfmodcor}  applied to the Hopf monad $T^\opp$, in virtue of
Remark~\ref{remopp}.
\end{rem}

\begin{exa} \label{HM-exa}
Let $A$ be a Hopf algebra in a braided right autonomous category~$\cc$. Consider the right Hopf monad $T=?\otimes A$
on~$\cc$, see Example~\ref{Ex2}. A right Hopf $T$-module is a nothing but a right Hopf module over $A$ in the usual sense,
that is, a triple $(M,r\co M\otimes A \to M, \rho\co M \to M \otimes A)$ such that $(M,r)$ is a right $A$-module, $(M,\rho)$
is a right $A$-comodule, and $\rho r= (m \otimes r) (\id_A \otimes \tau_{A,A} \otimes \id_M) (\rho \otimes \Delta)$, where
$\tau$ is the braiding of $\cc$, $m$ is the product of $A$, and $\Delta$ is coproduct of $A$. Assume now that $\cc$ splits
idempotents (see \cite{BKLT}). Then $M$ admits a coinvariant part, which is the object splitting the idempotent $r(S_A
\otimes \id_M)\rho$, where $S_A$ denotes the antipode of $A$. Moreover, $T$ preserves coinvariants (because $\otimes$ is
exact) and $T$ is conservative (because $A$ is a bialgebra). Therefore Theorems~\ref{Hopfmodthm} and~\ref{Hopfmodcor} apply
in this setting: we obtain the fundamental theorem of Hopf modules for categorical Hopf algebras. In the case where $S_A$ is
invertible, it was first stated in \cite{BKLT}.
\end{exa}

\subsection{Proof of Theorem~\ref{Hopfmodthm}.}\label{proofmodhopfthm}
Let $T$ be a right Hopf monad  on  a right autonomous category $\cc$, with right antipode $s^r$. Our proof will rely very
strongly on the properties of the morphism of functors $\Gamma \co \,? \otimes T(\un) \to T^2$ defined by:
\begin{equation}\label{defgamma}
\Gamma_X =(\rev_{X}(\id_X \otimes s^r_X) \otimes \id_{T^2(X)})\bigl (\id_X \otimes T_2(\rdual{T(X)},T(X))T(\rcoev_{T(X)})
\bigr )
\end{equation}
for any object $X$ of $\cc$.

Notice that, if $T$ is of the form $? \otimes A$, where $A$ is a Hopf algebra in a braided autonomous category (see
Example~\ref{ExHopf}), then $\Gamma_X = \id_X \otimes (S \otimes \id_A)\Delta$.

\begin{lem}\label{propgamr}
For any object $X$ of $\cc$, we have
\begin{enumerate}
 \renewcommand{\labelenumi}{{\rm (\alph{enumi})}}
 \item $\mu_X \Gamma_X= \eta_X \otimes T_0$;
 \item $ T(\mu_X) \Gamma_{T(X)} T_2(X,\un)=T(\eta_X)$;
 \item $ T\bigl ( (\id_{T(X)} \otimes \mu_\un)T_2(X,T(\un))\bigr )\Gamma_{X \otimes T(\un)}
       (\id_X \otimes T_2(\un,\un) )=T(\id_{T(X)} \otimes \eta_\un)\Gamma_X$;
 \item $ \Gamma_X(\id_X \otimes \eta_\un)=\eta_{T(X)} \eta_X$.
\end{enumerate}
\end{lem}

\begin{proof}
Let us prove Part (a). We have:
\begin{align*}
\mu_X \Gamma_X&=  (\rev_{X} \otimes \id_{T(X)})\bigl (\id_X \otimes (s^r_X \otimes \mu_X)
T_2(\rdual{T(X)},T(X))T(\rcoev_{T(X)}) \bigr )\\
&=  (\rev_{X} \otimes \eta_X)(\id_X \otimes \rcoev_{X}T_0) \quad \text{by \eqref{rant2}}\\
&= \eta_X \otimes T_0.
\end{align*}

Let us prove Part (b). We have:
\begin{align*}
 T(&\mu_X) \Gamma_{T(X)} T_2(X,\un)\\
 &= (\rev_{T(X)}(\id_{T(X)} \otimes s^r_{T(X)}) \otimes T(\mu_X))\\
&\phantom{XXXXX}  (\id_{T(X)} \otimes
T_2(\rdual{T^2(X)},T^2(X))T(\rcoev_{T^2(X)}))T_2(X,\un)\\
&= (\rev_{T(X)}(\id_{T(X)} \otimes s^r_{T(X)}T(\rdual{\mu}_X)) \otimes \id_{T^2(X)})(\id_{T(X)} \otimes
T_2(\rdual{T(X)},T(X)))\\
&\phantom{XXXXX}  T_2(X,\rdual{T(X)}\otimes T(X))T(\id_X\otimes\rcoev_{T(X)})\\
&= (\rev_{T(X)}(\id_{T(X)} \otimes s^r_{T(X)}T(\rdual{\mu}_X))T_2(X,\rdual{T(X)}) \otimes \id_{T^2(X)})\\
&\phantom{XXXXX}  T_2(X\otimes\rdual{T(X)},T(X))T(\id_X\otimes\rcoev_{T(X)})\quad \text{by \eqref{comoneq1}}\\
&= (T_0T(\rev_X)T(\id_X \otimes \rdual{\eta}_X)\otimes \id_{T^2(X)})\\
&\phantom{XXXXX}  T_2(X\otimes\rdual{T(X)},T(X))T(\id_X\otimes\rcoev_{T(X)})\quad \text{by \eqref{rant1}}\\
&= (T_0\otimes T(\eta_X))T_2(\un,X)T((\rev_X \otimes \id_X)(\id_X\otimes\rcoev_{T(X)}))=T(\eta_X).
\end{align*}

Let us prove Part (c). Denote by $\lhs_X$ the left hand side of Part (c). Firstly, using the functoriality of $T_2$, we
have:
\begin{align*}
\lhs_X= &\bigl ( \rev_{X \otimes T(\un)} (\id_{X \otimes T(\un)} \otimes \alpha_X) \otimes \id_{T(T(X) \otimes T(\un))} \bigr )\\
& \bigl(\id_{X\otimes T(\un)} \otimes T_2(\rdual{T(\un)} \otimes \rdual{T(X)}, T(X) \otimes T(\un)) \bigr)\\
& \bigl(\id_X \otimes T_2(\un, \rdual{T(\un)} \otimes \rdual{T(X)}\otimes  T(X) \otimes T(\un)) T(\rcoev_{T(X) \otimes
T(\un)}) \bigr )
\end{align*}
where $\alpha_X=s^r_{X \otimes T(\un)} T(\rdual{T_2(X,T(\un))}) T(\rdual{\mu}_\un \otimes \id_{\rdual{T(X)}})$. Now
\begin{align*}
\alpha_X &= (s^r_{T(\un)} \otimes s^r_X)T_2(\rdual{T^2(\un)},\rdual{T(X)})T(\rdual{\mu}_\un \otimes \id_{\rdual{T(X)}})
\quad \text{by \eqref{rantcomult}}\\
&= (s^r_{T(\un)}T(\rdual{\mu}_\un) \otimes s^r_X)T_2(\rdual{T(\un)},\rdual{T(X)})
\end{align*}
and, using \eqref{comoneq1},
\begin{align*}
&\bigl(\id_{X \otimes T(\un)} \otimes   T_2(\rdual{T(\un)},\rdual{T(X)})   \otimes \id_{T(X) \otimes T(\un)}\bigr )\\
& \phantom{XXX}\bigl ( \id_{X \otimes T(\un)} \otimes T_2(\rdual{T(\un)} \otimes
\rdual{T(X)}, T(X) \otimes T(\un)) \bigr)\\
&  \phantom{XXXXXX}   \bigl ( \id_X \otimes T_2(\un,\rdual{T(\un)} \otimes \rdual{T(X)}\otimes T(X) \otimes T(\un)) \bigr )\\
   & \phantom{XXXXXXXX}=\bigl(\id_X \otimes T_2(\un,\rdual{T(\un)}) \otimes T_2(\rdual{T(X)}, T(X)  \otimes T(\un))\bigr) \\
   & \phantom{XXXXXXXXXXXXX}\bigl(\id_X \otimes T_2(\rdual{T(\un)},\rdual{T(X)}\otimes T(X)  \otimes T(\un)) \bigr).
\end{align*}
Therefore, since   $\rev_{X \otimes T(\un)}=\rev_X (\id_X \otimes \rev_{T(\un)} \otimes \id_{\rdual{X}})$,  we have:
\begin{align*}
&\lhs_X = \bigl ( \rev_X (\id_X \otimes s^r_X) \otimes \id_{T(T(X) \otimes T(\un))} \bigr )\\
& \phantom{X}\bigl ( \id_X \otimes \rev_{T(\un)}(\id_{T(\un)} \otimes s^r_{T(\un)}T(\rdual{\mu}_\un))
T_2(\un,\rdual{T(\un)}) \otimes
T_2(\rdual{T(X)},T(X) \otimes T(\un)) \bigl )\\
& \phantom{XXXX} \bigl ( \id_X \otimes T_2(\rdual{T(\un)},\rdual{T(X)} \otimes T(X) \otimes T(\un))T(\rcoev_{T(X) \otimes
T(\un)}) \bigr ).
\end{align*}
Now $\rev_{T(\un)}(\id_{T(\un)} \otimes s^r_{T(\un)}T(\rdual{\mu}_\un)) T_2(\un,\rdual{T(\un)})=T_0T(\rdual{\eta}_\un)$ by
\eqref{rant1}. Hence, using~\eqref{comoneq2},
\begin{align*}
\lhs_X &= \bigl ( \rev_X (\id_X \otimes s^r_X) \otimes T(\id_{T(X)} \otimes \eta_\un) \bigr )\\
& \phantom{XXXX}  \bigl ( \id_X \otimes T_2(\rdual{T(X)},T(X))T(\rcoev_{T(X)}) \bigr )  \\
& = T(\id_{T(X)} \otimes \eta_\un) \Gamma_X.
\end{align*}

Let us prove Part (d). We have:
\begin{align*}
\Gamma_X &(\id_X \otimes \eta_\un) \\
&=(\rev_{X}(\id_X \otimes s^r_X) \otimes \id_{T^2(X)})\bigl (\id_X \otimes
T_2(\rdual{T(X)},T(X))T(\rcoev_{T(X)}) \eta_\un \bigr )\\
&=(\rev_{X}(\id_X \otimes s^r_X) \otimes \id_{T^2(X)})\bigl (\id_X \otimes
T_2(\rdual{T(X)},T(X))\eta_{\rdual{T(X)} \otimes T(X)} \rcoev_{T(X)}  \bigr )\\
&=(\rev_{X}(\id_X \otimes s^r_X \eta_{\rdual{T(X)}}) \otimes \eta_{T(X)}) (\id_X \otimes \rcoev_{T(X)} )
\quad \text{by \eqref{bimonad3}}\\
&=(\rev_{X}(\id_X \otimes \rdual{\eta}_X) \otimes \eta_{T(X)}) (\id_X \otimes \rcoev_{T(X)} )
\quad \text{by \eqref{rsitha}}\\
&=\eta_{T(X)} \eta_X.
\end{align*}
\end{proof}

\begin{lem}\label{propgamr2}
For any $T$-module $(M,r)$, $T(r)\Gamma_M\co (M,r) \otimes (T(\un),\mu_\un) \to (T(M),\mu_M)$ is a morphism of $T$-modules.
\end{lem}
\begin{proof}
Since
\begin{align*}
T_2(&\rdual{T(M)},T(M)) T(\rcoev_{T(M)})\mu_\un \\
 & = T_2(\rdual{T(M)},T(M))\mu_{\rdual{T(M)}\otimes T(M)}T^2(\rcoev_{T(M)})\\
& = (\mu_{\rdual{T(M)}} \otimes \mu_{T(M)})T^2_2(\rdual{T(M)},T(M))T^2(\rcoev_{T(M)}) \quad \text{by \eqref{bimonad1},}
\end{align*}
we have
\begin{align*}
T(r)&  \Gamma_M  (r \otimes \mu_\un)\\
&= (\rev_{M}(\id_M \otimes s^r_M) \otimes T(r))\bigl (r \otimes
T_2(\rdual{T(M)},T(M))T(\rcoev_{T(M)}) \mu_\un \bigr )\\
&= (\rev_{M}(r \otimes s^r_M \mu_{\rdual{T(M)}}) \otimes T(r)\mu_{T(M)})\\
& \phantom{XXXXXXX} \bigl (\id_{T(M)} \otimes
T^2_2(\rdual{T(M)},T(M))T^2(\rcoev_{T(M)}) \bigr )\\
&= (\rev_{M}(r \otimes s^r_M T(s^r_{T(M)}) T^2(\rdual{\mu}_M) \otimes T(r)\mu_{T(M)})\\
& \phantom{XXXXXXX} \bigl (\id_{T(M)} \otimes
T^2_2(\rdual{T(M)},T(M))T^2(\rcoev_{T(M)}) \bigr )  \quad \text{by \eqref{rsithm}}\\
&= (\rev_{M}(r \otimes s^r_M T(s^r_{T(M)})  \otimes T(r)\mu_{T(M)}T^2(\mu_M))\\
& \phantom{XXXXXXX} \bigl (\id_{T(M)} \otimes T^2_2(\rdual{T^2(M)},T^2(M))T^2(\rcoev_{T^2(M)}) \bigr ).
\end{align*}
Now
\begin{align*}
(T&(s^r_{T(M)}) \otimes \id_{T^3(M)})T^2_2(\rdual{T^2(M)},T^2(M))T^2(\rcoev_{T^2(M)})\\
&=T_2(\rdual{T(M)},T^3(M)) T \bigl( (s^r_{T(M)} \otimes \id_{T^3(M)}) T_2(\rdual{T^2(M)},T^2(M))
T(\rcoev_{T^2(M)}) \bigr )\\
&=T_2(\rdual{T(M)},T^3(M)) T \bigl( (\id_{T(M)} \otimes \Gamma_{T(M)})(\rcoev_{T(M)} \otimes \id_{T(\un)}) \bigr )\\
&=(\id_{T(\rdual{T(M)})} \otimes T(\Gamma_{T(M)})) T_2(\rdual{T(M)},T(M) \otimes T(\un)) T(\rcoev_{T(M)} \otimes
\id_{T(\un)}),
\end{align*}
and, using \eqref{Tmoddef},
\begin{equation*}
T(r)\mu_{T(M)}T^2(\mu_M)=\mu_M T^2(r\mu_M)=\mu_M T^2(rT(r))=T(r)\mu_{T(M)}T^3(r).
\end{equation*}
Therefore, we get
\begin{align*}
T(r)&\Gamma_M (r \otimes \mu_\un)T_2(M,T(\un))\\
&= \bigl (\rev_{M}(r \otimes s^r_M)   \otimes T(r)\mu_{T(M)}T^3(r) T(\Gamma_{T(M)}) \bigr )\\
& \phantom{XX} \bigl (\id_{T(M)} \otimes T_2(\rdual{T(M)},T(M) \otimes T(\un)) T(\rcoev_{T(M)} \otimes
\id_{T(\un)}) \bigr ) T_2(M,T(\un))\\
&= \bigl (\rev_{M}(r \otimes s^r_M) T_2(M,\rdual{T(M)})  \otimes T(r)\mu_{T(M)}  T(T^2(r)\Gamma_{T(M)})  \bigr ) \\
& \phantom{XX} T_2(M \otimes \rdual{T(M)},T(M) \otimes T(\un)) T(\id_M \otimes \rcoev_{T(M)} \otimes
\id_{T(\un)}) \quad \text{by \eqref{comoneq1}}\\
&= \bigl (\rev_{M}(r \otimes s^r_M T(\rdual{r})) T_2(M,\rdual{M})  \otimes T(r)\mu_{T(M)}T(\Gamma_M) \bigr ) \\
& \phantom{XX} T_2(M \otimes \rdual{M},M \otimes T(\un)) T(\id_M \otimes \rcoev_{M} \otimes
\id_{T(\un)}) \quad \text{by functoriality of $\Gamma$}\\
&= \bigl (T_0T(\rev_{M})  \otimes T(r)\mu_{T(M)}T(\Gamma_M) \bigr ) \\
& \phantom{XX} T_2(M \otimes \rdual{M},M \otimes T(\un)) T(\id_M \otimes \rcoev_{M} \otimes
\id_{T(\un)}) \quad \text{by Theorem~\ref{dual-ant}(b)}\\
&= \bigl (T_0  \otimes \mu_M T^2(r)T(\Gamma_M) \bigr ) T_2(\un,M \otimes T(\un))\\
& \phantom{XX}  T\bigl ((\rev_M \otimes \id_{M \otimes T(\un)}) (\id_M \otimes \rcoev_{M} \otimes
\id_{T(\un)})\bigr) \\
&= \mu_MT(T(r)\Gamma_M) \quad \text{by \eqref{comoneq2}.}
\end{align*}
Hence $T(r)\Gamma_M$ is a morphism of $T$-modules.
\end{proof}

Now we are ready to prove Theorem~\ref{Hopfmodthm}. Let $(M,r,\rho)$ be a right Hopf $T$-module and $i\co M^{\mathrm{co} T}
\to M$ be an equalizer of the pair $(\rho,\id_M \otimes \eta_\un)$. We will show that $r T(i)$ is an isomorphism in $\cc$
(by constructing an inverse) and we will check that $r T(i)$ is a morphism of right Hopf $T$-modules. By
Remark~\ref{remisoHopf}, this will prove the theorem.

Set $ \psi_M=T(r ) \Gamma_M \rho \co M \to T(M)$.
\begin{lem}\label{proofHopfmodlem} The morphism $\psi_M$ enjoys the following properties:
\begin{enumerate}
  \renewcommand{\labelenumi}{{\rm (\alph{enumi})}}
  \item $r \psi_M =\id_M$;
  \item $\psi_M r=\mu_M T(\psi_M)$;
  \item $T(\rho)\psi_M=T(\id_M \otimes \eta_\un)\psi_M$;
  \item $\psi_M i= \eta_M i$.
\end{enumerate}
\end{lem}
\begin{proof}
Since $M$ is a $T$-module and a right $T$-comodule, we have, by Lemma~\ref{propgamr}(a),
 $r \psi_M = r T(r ) \Gamma_M
\rho = r  \mu_M \Gamma_M \rho = r (\eta_M \otimes T_0) \rho = r  \eta_M  =\id_M$. Hence Part~(a). Moreover, we have:
\begin{align*}
\psi_M r & = T(r) \Gamma_M \rho r \\
 &= T(r) \Gamma_M (r \otimes \mu_\un)T_2(M,T(\un)) T(\rho) \quad \text{by \eqref{defHopfmod}} \\
 &=  \mu_M T(T(r)\Gamma_M) T(\rho) \quad \text{by Lemma~\ref{propgamr2}} \\
 &= \mu_M T(\psi_M).
\end{align*}
Hence Part (b). Now, since $M$ is a right Hopf $T$-module, we have, by Lemma~\ref{propgamr}(c),
\begin{align*}
T(\rho)\psi_M & = T(\rho r) \Gamma_M \rho \\
 &= T (r \otimes \mu_\un)T(T_2(M,T(\un))) T^2(\rho) \Gamma_M \rho\\
 &= T (r \otimes \mu_\un)T(T_2(M,T(\un)))  \Gamma_{M \otimes T(\un)} (\rho \otimes \id_{T(\un)})\rho\\
 &= T (r \otimes \id_{T(\un)}) T\bigl((\id_{T(M)} \otimes \mu_\un)T_2(M,T(\un))\bigr )  \Gamma_{M \otimes T(\un)}
  (\id_M \otimes T_2(\un,\un))\rho\\
 &= T (r \otimes \id_{T(\un)}) T(\id_{T(M)} \otimes \eta_\un)\Gamma_M \rho\\
 &= T(\id_M \otimes \eta_\un)\psi_M.
\end{align*}
Hence Part (c). Lastly, we have
\begin{align*} \psi_M i &= T(r) \Gamma_M \rho i = T(r) \Gamma_M (\id_M \otimes \eta_\un) i
 \\&= T(r) \eta_{T(M)} \eta_M i \quad \text{by Lemma~\ref{propgamr} (d)}\\&=  \eta_M r \eta_M i=
  \eta_M i.
 \end{align*}
Hence Part (d).
\end{proof}

By Lemma~\ref{proofHopfmodlem}(c), $\psi_M$ equalizes the pair  $(T(\rho),T(\id_M \otimes \eta_\un))$. Since, by assumption,
$T(i)$ is an equalizer of the pair $(T(\rho),T(\id_M \otimes \eta_\un))$, there exists a (unique) map $\phi_M\co M \to
T(M^{\mathrm{co} T})$ such that $\psi_M =T(i)\phi_M$.

Let us check that $\phi_M$ is inverse to $rT(i)$. We have $r T(i) \phi_M= r \psi_M = \id_{T(M)}$ by
Lemma~\ref{proofHopfmodlem}(a). In order to show that $\phi_M r T(i)= \id_{M^{\mathrm{co} T}}$, it is enough to check that
$T(i)\phi_M r T(i)=T(i)$ because $T(i)$, being an equalizer, is a monomorphism. Now
\begin{align*}
T(i)\phi_M r T(i)&=\psi_M r T(i)=\mu_M T(\psi_M) T(i) \quad \text{by Lemma~\ref{proofHopfmodlem}(b)}\\
&= \mu_M T(\eta_M) T(i) \quad \text{by Lemma~\ref{proofHopfmodlem}(d)}\\
&= T(i).
\end{align*}
Hence $rT(i)$ is an isomorphism in $\cc$.

Finally, let us check  that $rT(i)$ is a morphism of right Hopf modules. Firstly, we have $r  T(r  T(i)) = r  T(r
) T^2(i)=r \mu_M T^2(i)= r  T(i) \mu_{M^{\mathrm{co} T}}$. Therefore $r  T(i)$ is a morphism of $T$-modules. Secondly, we
have:
\begin{align*}
\rho r  T(i) & =(r  \otimes \mu_\un) T_2(M,T(\un)) T(\rho) T(i) \\
& = (r  \otimes \mu_\un) T_2(M,T(\un)) T(\id_M \otimes \eta_\un) T(i) \\
 & = (r  \otimes \mu_\un T(\eta_\un)) T_2(M,\un) T(i) \\ &  = (r  T(i)\otimes \id_T(\un))T_2(M^{\mathrm{co} T},\un).
\end{align*}
So $r  T(i)$ is also a morphism of  right $T$-comodules. This completes the proof of Theorem~\ref{Hopfmodthm}.

\subsection{Proof of Theorem~\ref{Hopfmodcor}.\label{proofmodhopfcor}}
Firstly, by Lemma~\ref{lemexHopf} and functoriality of $\mu$ and $T_2$, the assignments $X \mapsto
\bigl(T(X),\mu_X,T_2(X,\un) \bigr)$ and $ f \mapsto T(f)$ define a functor, denoted $\tilde{T}$, from $\cc$ to the category
of right Hopf $T$-modules.

Assume $\tilde{T}$ is an equivalence. In particular $\tilde{T}$ is conservative. If $f$ is a morphism in $\cc$ such that
$T(f)$ is an isomorphism in $\cc$, then $\tilde{T}(f)$ is an isomorphism (by Remark~\ref{remisoHopf}) and so is $f$ (since
$\tilde{T}$ is conservative). Hence $T$ is conservative.

Let us prove the converse. Let $(M,r,\rho)$  be a  right Hopf $T$-module and $M^{{\mathrm{co} T}}$ be its coinvariant part,
which exists by assumption.  By the universal property of equalizers, any  morphism of right Hopf modules $f\co (M,r,\rho)
\to (M',r',\rho')$  induces a morphism $M^{\mathrm{co} T}\to M'^{\mathrm{co} T}$.  This defines a functor $?^{\mathrm{co}
T}$ from the category of right Hopf $T$-modules to $\cc$. By Theorem~\ref{Hopfmodthm}, the functor $?^{\mathrm{co} T}$ is a
right quasi-inverse of $\tilde{T}$. Assume now that $T$ is conservative. It is enough to prove that, for any object $X$ of
$\cc$, $\eta_X\co X \to T(X)$ is the coinvariant part of the right $T$-comodule $(T(X), T_2(X,\un))$. Indeed, if this is
true, then $\eta_X$ induces a functorial isomorphism $X \iso \tilde{T}(X)^{\mathrm{co} T}$, so that $?^{\mathrm{co} T}$ is
also a left quasi-inverse of $\tilde{T}$. We have the following lemma:

\begin{lem}\label{lemTiequal}
Let $T$ be a right Hopf monad on a right autonomous category $\cc$. Then $T(\eta_X)$ is an equalizer of the pair
$\bigl(T(T_2(X,\un)),T(\id_{T(X)} \otimes \eta_\un)\bigr)$.
\end{lem}
\begin{proof}
Let $f \co Y \to T^2(X)$ be a morphism in $\cc$ equalizing the morphisms $T(T_2(X,\un))$ and $T(\id_{T(X)} \otimes
\eta_\un)$. If there exists $g \co Y \to T^2(X)$ such that $f=T(\eta_X)g$, then $g=\mu_X T(\eta_X)g=\mu_X f$, and so $g$ is
unique. All we have to check is that $f=T(\eta_X)\mu_X f$. We have:
\begin{align*}
T^2(\eta_X)f & =T\bigl (T(\mu_X)\Gamma_{T(X)}\bigr )T(T_2(X,\un))f \quad \text{by Lemma~\ref{propgamr}(b)}\\
&=T\bigl(T(\mu_X)\Gamma_{T(X)}\bigr )T(\id_{T(X)} \otimes \eta_\un)f \quad \text{by assumption}\\
&=T\bigl(T(\mu_X)\eta_{T^2(X)}\eta_{T(X)}\bigr )f\quad \text{by Lemma~\ref{propgamr}(d)}\\
&=T(\eta_{T(X)}\mu_X\eta_{T(X)})f \\
&=T(\eta_{T(X)})\,f\quad \text{by \eqref{mon-u}.}
\end{align*}
Hence $f=\mu_{T(X)}T(\eta_{T(X)})f=\mu_{T(X)} T^2(\eta_X)f=T(\eta_X) \mu_X f$.
\end{proof}

Now let $X$ be an object of $\cc$.  The right Hopf $T$-module $(T(X), \mu_X, T_2(X,\un))$ admits a coinvariant part
$i\co T(X)^{\mathrm{co} T} \to T(X)$ (by assumption) which is an equalizer of the pair $(T_2(X,\un), \id_X \otimes
\eta_\un)$.  Since $\eta_X$ equalizes this pair by \eqref{bimonad3}, there exists a unique morphism $j\co X \to
T(X)^{\mathrm{co} T}$ such that $\eta_X= i j$. To prove that $\eta_X$ is an equalizer, we need to show that $j$ is an
isomorphism (since $i$ is an equalizer). Now, applying $T$ to this situation, we have two equalizers for the pair
$(T(T_2(X,\un)), T(\id_X \otimes \eta_\un))$, namely  $T(\eta_X)$ (by lemma \ref{lemTiequal}) and $T(i)$ (because $T$
preserves coinvariants of right Hopf $T$-modules).  Therefore, since $T(\eta_X)=T(i) T(j)$, the morphism $T(j)$ is an
isomorphism, and so is $j$ because $T$ is conservative. This completes the proof of Theorem~\ref{proofmodhopfcor}.

\section{Integrals}\label{sect-int}

In this section, we study integrals of Hopf monads. In particular, using the decomposition theorem for Hopf modules, we
prove a theorem on the existence of universal left and right integrals for a Hopf monad.

\subsection{Integrals}
Let $T$ be a bimonad on a monoidal category $\cc$ and $K$ be an endofunctor of $\cc$. A \emph{($K$-valued) left integral of
$T$} is a functorial morphism $c\co T \to K$ such that:
\begin{equation}\label{intaxiomleft}
(\id_{T(\un)} \otimes c_X)T_2(\un,X)=\eta_\un \otimes c_X.
\end{equation}
A \emph{($K$-valued) right integral of $T$} is a functorial morphism $c\co T \to K$ such that:
\begin{equation}\label{intaxiomright}
(c_X \otimes \id_{T(\un)})T_2(X,\un)=c_X \otimes \eta_\un.
\end{equation}

\begin{exa}
Let $A$ be a bialgebra in a braided category $\cc$. Consider the bimonad $T=A\otimes ?$ on $\cc$, see Example~\ref{Ex2}. Let
$\chi \co A \to k$ be a morphism in $\cc$. Set $K=k\, \otimes ?$ and define $c\co T \to K$ by $c_X=\chi \otimes \id_X$. Then
$c$ is a $K$-valued left (resp.\@ right) integral of $T$ if and only if $\chi$ is a $k$-valued left (resp.\@ right) integral
of $A$.
\end{exa}

Let $T$ be a bimonad on a monoidal category $\cc$. A left (resp.\@ right) integral $\lambda\co T \to I$ of $T$ is
\emph{universal} if, for any left (resp.\@ right) integral $c \co T \to K$ of~$T$, there exists a unique functorial morphism
$f\co I \to K$ such that $c=f\lambda$.

Note that a  universal left (resp.\@ right) integral of $T$ is unique up to unique functorial isomorphism.\\

\subsection{Existence of universal integrals}

Recall that, according to  Lemma~\ref{lemcomod}, if $T$ is a  comonoidal endofunctor of an autonomous category $\cc$ and $X$
an object of $\cc$, then we have a right $T$-comodule $\rdual{(T(X), T_2(\un,X))}$ and a left $T$-comodule $\ldual{(T(X),
T_2(X,\un))}$.

\begin{prop}\label{existint}
Let $T$ be a bimonad on an autonomous category $\cc$.
\begin{enumerate}
\renewcommand{\labelenumi}{{\rm (\alph{enumi})}}
\item If, for any object $X$ of $\cc$, the right $T$-module $\rdual{(T(X), T_2(\un,X))}$ admits coinvariants,
then $T$ admits a universal left integral $\lambda^l\co T \to I_l$, which is characterized by the fact that
$\rdual{(\lambda^l_X)}\co \rdual{I_l(X)} \to \rdual{T(X)}$ is the coinvariant part of $\rdual{(T(X), T_2(\un,X))}$ for all
object $X$ of $\cc$.
\item If, for any object $X$ of $\cc$, the left $T$-module $\ldual{(T(X), T_2(X,\un))}$ admits coinvariants,
then $T$ admits a universal right integral $\lambda^r \co T \to I_r$,  which is characterized by the fact that
$\ldual{(\lambda^r_X)} \co \ldual{I_r(X)} \to \ldual{T(X)}$ is the coinvariant part of $\ldual{(T(X), T_2(X,\un))}$ for all
object $X$ of $\cc$.
\end{enumerate}
\end{prop}

\begin{proof}
We prove Part (a), from which Part (b) can be deduced using the opposite bimonad. For an object $X$ of $\cc$, we have
$\rdual{(T(X), T_2(\un,X))}=(\rdual{T(X)}, \rho^r_X)$, with $\rho^r_X=\rdual{T_2(\un,X)}(\id_{\rdual{T(X)}} \otimes
\rcoev_{T(\un)})$.

Firstly,  observe that a functorial morphism $c\co T \to K$ is a left $K$-valued integral of $T$ if and only if, for any
object $X$ of $\cc$, the morphism $\rdual{c}_X \co \rdual{K(X)}\to \rdual{T(X)}$ equalizes the pair
$(\rho^r_X,\id_{\rdual{T(X)}} \otimes \eta_\un)$. Indeed, we have $(\id_{T(\un)} \otimes c_X)T_2(\un,X)=\eta_\un \otimes
c_X$ if and only~if
\begin{align*}
(\id_{\rdual{T(X)} \otimes T(\un)} &\otimes \rev_{K(X)})(\id_{T(\un)} \otimes c_X)T_2(\un,X)(\rcoev_{T(X)} \otimes
 \id_{\rdual{K(X)}})\\
&=(\id_{\rdual{T(X)} \otimes T(\un)} \otimes \rev_{K(X)})(\eta_\un \otimes c_X)(\rcoev_{T(X)} \otimes
 \id_{\rdual{K(X)}}),
\end{align*}
that is, if and only if $\rho^r_X\rdual{c}_X= \rdual{c}_X \otimes \eta_\un$.

Now assume that, for any object $X$ of $\cc$, $\rdual{(T(X), T_2(\un,X))}$ admits a coinvariant part $i_X\co E(X)\to
\rdual{T(X)}$. The morphism $i_X$ is an equalizer of the pair $(\rho^r_X,\id_{\rdual{T(X)}} \otimes \eta_\un)$. Define
$\lambda^l_X=\ldual{i_X}\co T(X)\cong \ldual{(\rdual{T(X)})} \to \ldual{E(X)}$. Using the universal property of equalizers,
one checks easily that the assignment $X \mapsto \ldual{E(X)}$ defines an endofunctor $I_l=\ldual{E}$ of $\cc$ and that
$\lambda^l\co T \to I_l$ is a functorial morphism. By the initial remark, $\lambda^l$ is a left integral for $T$.

Let $K$ be an endofunctor of $\cc$ and $c$ be a left $K$-valued integral. Again by the initial remark, there exists a unique
morphism $a_X \co \rdual{K(X)} \to E(X)$ such that $\rdual{c_X}=i_X a_X$. Using the universal property of equalizers, one
checks that $a \co \rdual{K} \to E$ is a functorial morphism. Dualizing, we obtain that there  exists a unique functorial
morphism $f= \ldual{a}\co I_l=\ldual{E} \to K$ such that $c = f \lambda^l$. Hence $\lambda^l$ is a universal left integral.
\end{proof}

Recall that, for any endofunctor $K$ of an autonomous category $\cc$, we form two endofunctors $\rexcla{K}=\rdual{?}\circ
K^\opp\circ \ldual{?}^\opp $ and $\lexcla{K}=\ldual{?}\circ K^\opp\circ ?^{\vee\opp}$, see Section~\ref{sect-autono}. This
defines two functors $\lexcla{?},\,\rexcla{?}\co \End(\cc)^\opp \to \End(\cc)$ such that $\lexcla{?}$ and $?^{!\opp}$ are
quasi-inverse.

\begin{thm}\label{Iequiv}
Let $T$ be a Hopf monad on an autonomous category $\cc$. Assume that left Hopf $T$-modules and right Hopf $T$-modules admit
coinvariants which are preserved by $T$ (such is the case if equalizers exist in $\cc$ and are preserved by $T$). Suppose
moreover that $T$ is conservative. Then there exist two auto-equivalences $I_l$ and $I_r$ of the category $\cc$, a universal
$I_l$-valued left integral of $T$, and a universal $I_r$\trait valued right integral of $T$. Moreover $I_l^!$ is
quasi-inverse to $I_r$ and $\lexcla{I}_r$ is quasi-inverse to $I_l$.
\end{thm}

\begin{exa}\label{exintA}
Let $A$ be a Hopf algebra, with invertible antipode $S_A$, in a braided autonomous category $\cc$. Consider the Hopf monad
$T=A\otimes ?$ on~$\cc$, and assume that $\cc$ splits idempotents as in Example~\ref{HM-exa}. Then Theorem~\ref{Iequiv}
applies, and there exists a universal left integral $\lambda^l\co T \to I_l$ and a universal right integral $\lambda^r\co T
\to I_r$ on $T$, where $I_l$ and $I_r$ are equivalences of $\cc$ such that $I_l^!$ is quasi-inverse to $I_r$. Moreover, by
Proposition~\ref{existint} and since $\otimes$ commutes with equalizers, there exists objects $k_l$ and $k_r$ and morphisms
$\int^l\co A \to k_l$ and $\int^r\co A \to k_r$ in $\cc$ such that $I_l=k_l\otimes?$, $\lambda^l_X=\int^l \otimes \id_X$,
$I_r=k_r \otimes ?$ and $\lambda^r_X=\int^r \otimes \id_X$. The morphisms $\int^l$ and $\int^r$ are universal left and right
integrals of the Hopf algebra $A$ respectively. Since $I_l^!=? \otimes \rdual{k_l}$ is quasi-inverse to $I_r=k_r \otimes ?$,
we see that $k_r \otimes \rdual{k_l} \cong \un$. Hence we may assume $k_r = k_l$ and this object, denoted $\mathrm{Int}$, is
$\otimes$-invertible. Let us summarize this discussion: there exists a $\otimes$-invertible object $\mathrm{Int}$ of $\cc$,
a universal left integral $\int^l\co A \to \mathrm{Int}$  and a universal right integral $\int^r\co A \to \mathrm{Int}$ on
$A$. This result was first proven in \cite{BKLT}.
\end{exa}

\begin{rem}
In Theorem~\ref{Iequiv}, the auto-equivalences $I_l$ and $I_r$ are in general not isomorphic to $1_\cc$. Such is already the
case in the setting of Example~\ref{exintA} (since in Example 3.1 of~\cite{BKLT} the object $\mathrm{Int}$ is not isomorphic
to $\un$).
\end{rem}

\begin{proof}[Proof of Theorem~\ref{Iequiv}]
By Proposition~\ref{existint},  $T$ admits universal left and right integrals, which we denote $\lambda^l \co T \to I_l$ and
$\lambda^r \co T \to I_r$  respectively.

Let $X$ be an object of $\cc$. By Lemma~\ref{lemexHopf}, $(T(X),\mu_X,T_2(X,\un))$ is a right Hopf $T$-module. So
$\ldual{(T(X),\mu_X,T_2(X,\un))}$ is a left Hopf $T$-module (by Lemma~\ref{leftrighthopfmod}) whose coinvariant part is
$\ldual{(\lambda^r_X)}\co \ldual{I_r(X)} \to \ldual{T(X)}$ (by Proposition~\ref{existint}(b)). Therefore, by
Theorem~\ref{Hopfmodthm}, we have an isomorphism $ T(\ldual{I_r(X)}) \to \ldual{T(X)}$ of left Hopf $T$\trait modules.
Applying the right dual functor of Lemma~\ref{leftrighthopfmod} to this isomorphism, we obtain an isomorphism of right Hopf
$T$-modules $T(X) \iso \rdual{T(\ldual{I_r(X)})}$. Likewise, using Lemma~\ref{leftrighthopfmod},
Proposition~\ref{existint}(a), and Theorem~\ref{Hopfmodthm}, we have an isomorphism of right Hopf $T$\trait modules between
$\rdual{T(\ldual{I_r(X)})}$ and $T\bigl (\rdual{I_l(\ldual{I_r(X)})} \bigr )$. By Theorem~\ref{Hopfmodcor}, we deduce from
this a functorial isomorphism
\begin{equation*}
X \simeq T(X)^{\mathrm{co} T} \iso \bigl( \rdual{T(\ldual{I_r(X)})}\bigr)^{\mathrm{co} T}\simeq
\rdual{I_l(\ldual{I_r(X)})}=I_l^!I_r(X).
\end{equation*}

Similarly, applying the previous construction to $T^\opp$, we obtain a functorial isomorphism $1_\cc \iso \lexcla{I}_r I_l$.
Hence, using Remark~\ref{rem-excla}, we obtain that $I_l$ and $I_r$ are auto-equivalences of the category $\cc$ such that
$I_l^!$ is quasi-inverse to $I_r$ and $\lexcla{I}_r$ is quasi-inverse to $I_l$.
\end{proof}

\subsection{Integrals and antipodes}
In this section we show, as in the classical case, that left (resp.\@ right) integrals are transported to right (resp.\@
left) integrals via the antipode. It turns out that this works only for integrals with values in endofunctors admitting a
right adjoint.

\begin{prop}\label{antipandintegs}
Let $T$ be a bimonad on an autonomous category $\cc$ and $J$, $K$ be endofunctors of $\cc$.
\begin{enumerate}
\renewcommand{\labelenumi}{{\rm (\alph{enumi})}}
\item
Assume $T$ is a right Hopf monad. Let $c\co T \to J$ be a left integral of $T$ and suppose we have a functorial morphism
$\varepsilon \co J \,\,\lexcla{K} \to 1_\cc$. For any object $X$ of~$\cc$, set:
\begin{equation*}
c^{(\varepsilon)}_X=s^r_{\ldual{K(X)}}T\bigl(\rexcla{c}_{K(X)} \,\rexcla{\varepsilon}_X\bigr)\co T(X) \to K(X).
\end{equation*}
Then the functorial morphism $c^{(\varepsilon)}\co T \to K$ is right integral of $T$.
\item
Assume $T$ is a left Hopf monad. Let $d\co T \to K$ be a right integral of $T$ and suppose we have a functorial morphism
$\varepsilon'\co K \rexcla{J} \to 1_\cc$. For any object $X$ of $\cc$, set:
\begin{equation*}
\leftidx{^{(\varepsilon')}}{d}{_X}=s^l_{\rdual{J(X)}}T\bigl(\lexcla{d}_{J(X)}\,\lexcla{\varepsilon'}_X\bigr)\co T(X) \to
J(X).
\end{equation*}
Then the functorial morphism $\leftidx{^{(\varepsilon')}}{d}{} \co T \to J$ is a left integral of $T$.
\item Assume that $T$ is a Hopf monad. Suppose that $\lexcla{K}$ is right adjoint to $J$, with adjunction morphisms
$\alpha \co J\,\lexcla{K} \to 1_\cc$ and $\beta \co 1_\cc \to \lexcla{K}J$. Then the assignment $c \mapsto c^{(\alpha)}$
defines a bijection between $J$-valued left integrals of $T$ and $K$-valued right integrals of $T$, whose inverse is given
by $d \mapsto \leftidx{^{(\rexcla{\beta})}}{d}{}$.
\end{enumerate}
\end{prop}
\begin{proof}
Let us prove Part (a).  Set $d=c^{(\varepsilon)}$. Let $X$ be an object of $\cc$ and set $Y= \ldual{K(X)}$ and $y=
\rexcla{c}_{K(X)} \rexcla{\varepsilon}_X$. By \eqref{intaxiomleft}, we have $\rdual{T_2(\un,Y)}(y \otimes
\id_{\rdual{T(\un)}})=y \otimes \rdual{\eta_\un}$. Therefore:
\begin{align*}
d_X & \otimes T(\rcoev_\un) \eta_\un\\
& =(d_X \otimes \eta_\un \rcoev_\un T_0)T_2(X,\un) \quad \text{by \eqref{comoneq2}}\\
& = \bigl (s^r_Y T(y)\otimes (s^r_\un \otimes\mu_\un) T_2(\rdual{T(\un)},T(\un))T(\rcoev_{T(\un)}) \bigr)
    T_2(X,\un)\quad \text{by \eqref{rant2}}\\
& = \bigl (s^r_Y \otimes (s^r_\un \otimes\mu_\un) T_2(\rdual{T(\un)},T(\un)) \bigr) T_2(\rdual{T(Y)},\rdual{T(\un)}\otimes
T(\un))
    T(y \otimes \rcoev_{T(\un)})\\
& = \bigl ((s^r_Y \otimes s^r_\un) T_2(\rdual{T(Y)},\rdual{T(\un)})\otimes\mu_\un \bigr) \\
& \phantom{XXXXXXX} T_2(\rdual{T(Y)}\otimes
\rdual{T(\un)}, T(\un)) T(y\otimes \rcoev_{T(\un)})\quad \text{by \eqref{comoneq1}}\\
& = \bigl (s^r_Y T(\rdual{T_2(\un,Y)})\otimes\mu_\un \bigr) T_2(\rdual{T(Y)}\otimes \rdual{T(\un)}, T(\un)) T(y\otimes
\rcoev_{T(\un)})\quad \text{by \eqref{rantcomult}}
\end{align*}
and so
\begin{align*}
d_X & \otimes T(\rcoev_\un) \eta_\un\\
& = (s^r_Y \otimes\mu_\un )T(\rdual{T_2(\un,Y)}) T_2(\rdual{T(Y)}, T(\un))\\
& \phantom{XXXXXXX}  T\bigl ((\rdual{T_2(\un,Y)} \otimes \id_{T(\un)})(y\otimes \rcoev_{T(\un)}) \bigr)\\
& = (s^r_Y \otimes\mu_\un )T(\rdual{T_2(\un,Y)}) T_2(\rdual{T(Y)}, T(\un))\\
& \phantom{XXXXXXX}  T(y \otimes (\rdual{\eta}_\un \otimes \id_{T(\un)})\rcoev_{T(\un)})
  \quad \text{by \eqref{intaxiomleft}}\\
& = (s^r_Y \otimes\mu_\un )T(\rdual{T_2(\un,Y)}) T_2(\rdual{T(Y)}, T(\un))T(y \otimes \eta_\un \rcoev_\un) \\
& = \bigl (s^r_Y T(y)\otimes\mu_\un T(\eta_\un)T(\rcoev_\un) \bigr)T(X,\un) \\
& = (d_X \otimes T(\rcoev_\un))T(X,\un) \quad \text{by \eqref{mon-u}.}
\end{align*}
Since $\rcoev_\un$ is an isomorphism, we get $d_X \otimes \eta_\un=(d_X \otimes \id_{T(\un)})T(X,\un)$. Hence $d$ is a
$K$-valued right integral of $T$.

Part (b) is obtained by applying Part (a) to the opposite Hopf monad. Let us prove Part (c). Let $c\co T \to J$ be a left
integral of $T$. Set $c'={}^{(\rexcla{\beta})}(c^{(\alpha)})$. Let us check that $c'=c$. For any object $X$ of $\cc$, we
have:
\begin{align*}
c'_X&=s^l_{\rdual{J(X)}}T\bigl( \ldual{(s^r_{\ldual{K(\rdual{J(X)})}}T(\rexcla{c}_{K(\rdual{J(X)})}
\,\rdual{\alpha}_{J(X)}))} \,\beta_X\bigr) \\
  &=s^l_{\rdual{J(X)}}T( \ldual{T(\rexcla{c}_{K(\rdual{J(X)})} \,\rdual{\alpha}_{J(X)})} \, \ldual{(
s^r_{\ldual{K(\rdual{J(X)})}} )} \,\beta_X)\\
  &=s^l_{\rdual{J(X)}}T\bigl(\ldual{T(\rdual{T(\beta_X)}
\,\rexcla{c}_{K(\rdual{J(X)})} \,\rdual{\alpha}_{J(X)}))} \,\ldual{ s^r_X }\bigr)\\
  &=\alpha_{J(X)} c_{\lexcla{K}J(X)} T(\beta_X) s^l_{\rdual{T(X)}} T\bigl(\ldual{s^r_{X}}\bigr)\\
&=\alpha_{J(X)} c_{\lexcla{K}J(X)}T(\beta_X) \quad \text{by Proposition~\ref{anti-inv}}\,,\\
&=\alpha_{J(X)} J(\beta_X) c_X=c_X \quad \text{by adjunction.}
\end{align*}
Hence ${}^{(\rexcla{\beta})}(c^{(\alpha)})=c$. Applying this to the opposite Hopf monad, we obtain that
${({}^{(\rexcla{\beta})}d)}^{(\alpha)}=d$ for any right integral $d \co T \to K$.
\end{proof}

\begin{prop}\label{corint}
Under the hypotheses of Theorem~\ref{Iequiv}, in the bijective correspondence of Proposition~\ref{antipandintegs}(c), a
universal left integral of $T$ is transformed into a universal right integral of $T$, and conversely.
\end{prop}
\begin{proof}
By Theorem~\ref{Iequiv},  there exist a universal left integral $\lambda^l\co T \to J$ of $T$ and a universal right integral
$\lambda^r\co T \to K$ of $T$ such that $\lexcla{K}$ is quasi-inverse (and, in particular, right adjoint) to $J$. Denote
$\alpha \co J\,\, \lexcla{K} \to 1_\cc$ and $\beta \co 1_\cc \to \lexcla{K}J$ the adjunction isomorphisms. Using
Proposition~\ref{antipandintegs}, define a right integral $c : T \to K$ and a left integral $d : T \to J$ by
$c={(\lambda^l)}^{(\alpha)}$ and  $d= {}^{(\rexcla{\beta})}{(\lambda^r)}$. We have to show that they are universal. Since
$\lambda^l$ and $\lambda^r$ are universal, there exist unique functorial morphisms $f\co K \to K$ and $g\co J\to J$ such
that $d=f \lambda^l$ and $c=g \lambda^r$. It is sufficient to prove that $f$ and $g$ are isomorphisms. Since $d^{(\alpha)}=
\lambda^r$ by Proposition~\ref{antipandintegs}(c), we have, for any object $X$ of $\cc$,
\begin{equation*}
\lambda^r_X= d^{(\alpha)}_X= {(\lambda^l)}^{(\alpha)}_X T(\rexcla{f}_{K(X)})=c_X T(\rexcla{f}_{K(X)}) =K(\rexcla{f}_{K(X)})
c_X=K(\rexcla{f}_{K(X)})g_X \lambda^r_X.
\end{equation*}
Thus $K(\rexcla{f}_{K})g =\id_{K}$ by the universal property of $\lambda^r$. By functoriality of $g$, we also have $g
K(\rexcla{f}_{K})= \id_{K}$. Hence $g$ is an isomorphism. Similarly one shows that $f$ is an isomorphism.
\end{proof}

\section{Semisimplicity}\label{sect-semsim}

In this section, we define semisimple and separable monads, and give a characterization of semisimple Hopf monads (which
generalizes Maschke's theorem).

\subsection{Semisimple monads}
Let $T$ be a monad on a category $\cc$. Recall that for any object $Y$ of $\cc$, $(T(Y),\mu_Y)$ is a $T$-module. Such a
$T$-module is said to be \emph{free}. If $(M, r)$ is a $T$-module, then $r$ is a $T$-linear morphism from the free module
$(T(M),\mu_M)$ to $(M,r)$. Note that $\eta_M\co M \to T(M)$ is a section of $r$ in $\cc$, but in general $\eta_M$ is not
$T$-linear.

\begin{prop}\label{equiv-ss-def}
Let $T$ be a monad on a category $\cc$. The following conditions are equivalent:
\begin{enumerate}
\renewcommand{\labelenumi}{{\rm (\roman{enumi})}}
\item For any $T$-module $(M, r)$, the $T$-linear morphism $r$ has a $T$-linear section;
\item Any $T$-linear morphism has a $T$-linear section if and only if the underlying morphism  in $\cc$  has a section;
\item Any $T$-module is a $T$-linear retract of a free $T$-module.
\end{enumerate}
\end{prop}
 A \emph{semisimple} monad is a monad satisfying the equivalent conditions of Proposition~\ref{equiv-ss-def}.
\begin{rem}
Assume that $\cc$ is abelian semisimple, and $T$ is additive. Then $T$ is semisimple if and only if the category $T\ti\cc$
of $T$-modules is abelian semisimple.
\end{rem}
\begin{proof}[Proof of Proposition~\ref{equiv-ss-def}]
We have (ii) implies (i) since $\eta_M$ is a section of $r$ in $\cc$. Clearly (i) implies (iii). Let us show that (iii)
implies (ii). Let $f \co (M,r) \to (N,s)$ a $T$-linear morphism between two $T$-modules and $g \co N \to M$ be a section of
$f$ in~$\cc$. By assumption, $(N,s)$ is a retract of $(T(X),\mu_X)$ for some object $X$ of $\cc$. Let $p\co T(X) \to N$ and
$i\co N \to T(X)$ be $T$-linear morphisms such that $pi=\id_N$. Set $g'=r T(g p\eta_X)i\co N \to M$. We have:
\begin{align*}
g's&=r T(g p\eta_X)is=r T(g p\eta_X)\mu_XT(i)=r \mu_M T^2(g p\eta_X)T(i)\\
&=rT(r)T(T(g p\eta_X)i)=rT(g')
\end{align*}
and $fg'=fr T(g p\eta_X)i=sT(fg p\eta_X)i=sT(p)T(\eta_X)i=p\mu_XT(\eta_X)i=pi=\id_N$. Hence $g'$ is a $T$-linear section of
$f$.
\end{proof}

\subsection{Separable monads}
Let $A$ be an algebra in a monoidal category. In particular $A \otimes A$ is a $A$-bimodule, and the multiplication $m \co A
\otimes A \to A$ of $A$ is a morphism of $A$-bimodules. Recall that $A$ is \emph{separable} if $m$ has a section
$\sigma\co A \to A \otimes A$  as a morphism of $A$-bimodules, which means that:
\begin{equation*}
(m \otimes \id_A)(\id_A \otimes \sigma)=(\id_A\otimes m)(\sigma \otimes \id_A) \quad \text{and} \quad m \sigma=\id_A.
\end{equation*}
In this case, set $\gamma= \sigma u \co \un \to A \otimes A$, where $u\co \un \to A$ is the unit of $A$. Then the morphism
$\gamma$ satisfies:
\begin{equation*}
(m \otimes \id_A) (\id_A \otimes \gamma)=(\id_A \otimes m) (\gamma \otimes \id_A) \quad \text{and} \quad m \gamma = u.
\end{equation*}
Conversely if $\gamma\co \un\to A \otimes A$ satisfies the above equation, then $A$ is separable and the section of $m$ is
$\sigma=(m \otimes \id_A) (\id_A \otimes \gamma)\co A \to A \otimes A$. We extend this notion to monads.
\begin{prop}\label{equiv-sep-def}
Let $T$ be a monad on a category $\cc$. The following conditions are equivalent:
\begin{enumerate}
\renewcommand{\labelenumi}{{\rm (\roman{enumi})}}
\item One may choose functorially for each $T$-module $(M,r)$ a $T$-linear section $\sigma_{(M,r)}$
of the morphism $r \co T(M) \to M$. Here `functorially' means that, for any $T$-linear morphism $f\co (M,r) \to (N,s)$, we
have:
\begin{equation*}
 \sigma_{(N,s)} f= T(f)  \sigma_{(M,r)};
\end{equation*}
\item There exists a functorial morphism $\varsigma\co T \to T^2$ such that:
\begin{equation*}
 T(\mu_X)  \varsigma_{T(X)}= \varsigma_X \mu_{X} = \mu_{T(X)}  T(\varsigma_X) \quad \text{and} \quad \mu_X \varsigma_X = \id_{T(X)};
\end{equation*}
\item There exists a functorial morphism $\gamma \co 1_\cc \to T^2$ such that
\begin{equation*}
 T(\mu_X) \gamma_{T(X)}=\mu_{T(X)} T(\gamma_X)  \quad \text{and} \quad \quad \mu_X \gamma_X= \eta_X.
\end{equation*}
\end{enumerate}
\end{prop}
A \emph{separable} monad is  a monad satisfying the equivalent conditions of Proposition~\ref{equiv-sep-def}.

\begin{proof}
Let us show that (i) implies (ii). Define $\varsigma = \sigma_{(T,\mu)}\co T \to T^2$, which is clearly a functorial
morphism such that $\varsigma_X \mu_X =\id_{T(X)}$. Since $\mu_X$ is $T$-linear, the functoriality of $\sigma$ gives
$\varsigma_X \mu_X = T(\mu_X) \varsigma_{T(X)}$. Finally, using the $T$-linearity of $\sigma$, we have $\varsigma_X \mu_X=
\mu_{T(X)} T(\varsigma_X)$.

Let us show that (ii) implies (iii). Set $\gamma= \varsigma \eta \co 1_\cc \to T^2$. Then
\begin{align*}
T(\mu_X) \gamma_{T(X)}& = T(\mu_X)  \varsigma_{T(X)} \eta_{T(X)}= \varsigma_X \mu_{X} \eta_{T(X)}\\
&= \varsigma_X \mu_{X} T(\eta_X)=\mu_{T(X)}  T(\varsigma_X)T(\eta_X)=\mu_{T(X)} T(\gamma_X),
\end{align*}
and $\mu_X \gamma_X= \mu_X \varsigma_X \eta_X= \eta_X$.

Let us show that (iii) implies (i). For any $T$-module $(M,r)$, set $\sigma_{(M,r)}=T(r) \gamma_M$. We have:
\begin{equation*}
\sigma_{(M,r)} r =T(r) \gamma_M r =T(r)T^2(r) \gamma_{T(M)} =T(r)T(\mu_M) \gamma_{T(M)}=T(r)\sigma_{(T(M),\mu_M)}
\end{equation*}
and $r\sigma_{(M,r)}=rT(r) \gamma_M =r\mu_M \gamma_M=r \eta_M =\id_M$. Therefore $\sigma_{(M,r)}$ is a $T$-linear section of
$r$. Finally, for any $T$-linear morphism $f\co (M,r) \to (N,s)$, we have:
\begin{equation*}
\sigma_{(N,s)} f= T(s) \gamma_N f= T(s)T^2(f) \gamma_M=T(sT(f)) \gamma_M =T(fr) \gamma_M=T(f)\sigma_{(M,r)}.
\end{equation*}
Hence $\sigma$ is functorial.
\end{proof}

\subsection{Cointegrals}\label{sect-cointeg}
Let $(T,\mu,\eta)$ be a bimonad on a monoidal category $\cc$. A \emph{cointegral} of $T$ is a morphism $\Lambda\co \un \to
T(\un)$ satisfying:
\begin{equation}
\mu_\un T(\Lambda)=\Lambda T_0.
\end{equation}
This condition means that $\Lambda$ is a morphisms of $T$-modules from $(\un,T_0)$ to $(T(\un),\mu_\un)$.

\begin{exa}\label{exintB}
Let $A$ be a bialgebra in a braided category $\cc$  and $\Lambda \co \un \to A$ be a morphism in $\cc$. Then $\Lambda$
is an cointegral of the bimonad $A \otimes ?$ (resp.\@ $? \otimes A$) of Example~\ref{Ex2} if and only if $\Lambda$ is a
left (resp.\@ right) integral \emph{in} $A$.
\end{exa}

\subsection{Maschke Theorem}
In this section, we extend the Theorem of Maschke, which characterizes semisimple Hopf algebras  in terms of
(co)integrals,  to the (non-linear) setting of Hopf monads.

\begin{thm}[Maschke Theorem for Hopf monads]\label{thmMaschke}
Let $T$ be a right Hopf monad on a right autonomous category. The following assertions are equivalent:
\begin{enumerate}
\renewcommand{\labelenumi}{{\rm (\roman{enumi})}}
 \item $T$ is semisimple;
 \item $T$ is separable;
 \item $T$ admits a cointegral $\Lambda\co \un \to T(\un)$ such that $T_0\Lambda=\id_\un$.
\end{enumerate}
\end{thm}
\begin{proof}
We have (ii) implies (i) by Propositions~\ref{equiv-ss-def} and \ref{equiv-sep-def}.

Let us show that (i) implies (iii). Consider the $T$-module $(\un, T_0)$. Since $T$ is semisimple, there exists a $T$-linear
morphism $\Lambda \co (\un, T_0) \to (T(\un),\mu_\un)$ such that $T_0 \Lambda = \id_\un$. The $T$-linearity of $\Lambda$
means $\mu_\un T(\Lambda)= \Lambda T_0$, that is, $\Lambda$ is  a cointegral.

Finally, let us show that (iii) implies (ii). Consider the morphisms $\Gamma_X\co X \otimes T(\un) \to T^2(X)$ as defined in
\eqref{defgamma}. Set $\gamma_X= \Gamma_X (\id_X \otimes \Lambda)\co X \to T^2(X)$. By Lemma~\ref{propgamr2} applied to the
$T$-module $(T(X),\mu_X)$, we have:
\begin{equation}\label{demMaskeq1}
T(\mu_X)\Gamma_{T(X)}(\mu_X \otimes \mu_\un)T_2(T(X),T(\un))=\mu_{T(X)} T\bigl(T(\mu_X)\Gamma_{T(X)}\bigr).
\end{equation}
 Composing the left hand side of \eqref{demMaskeq1} with $T(\eta_X \otimes \Lambda)$ gives:
\begin{align*}
 T(\mu_X)\Gamma_{T(X)}&(\mu_X \otimes \mu_\un)T_2(T(X),T(\un))T(\eta_X \otimes \Lambda)\\
 &= T(\mu_X)\Gamma_{T(X)}(\mu_X T(\eta_X) \otimes \mu_\un T(\Lambda))T_2(X,\un)\\
 &= T(\mu_X)\Gamma_{T(X)}(\id_{T(X)} \otimes \Lambda T_0)T_2(X,\un)\\
 &= T(\mu_X)\gamma_{T(X)}.
\end{align*}
 Composing the right hand side of \eqref{demMaskeq1} with $T(\eta_X \otimes \Lambda)$ gives:
\begin{align*}
 \mu_{T(X)} &T\bigl(T(\mu_X)\Gamma_{T(X)}\bigr)T(\eta_X \otimes \Lambda)\\
 &= \mu_{T(X)} T\bigl(T(\mu_X)T^2(\eta_X)\Gamma_X(\id_X \otimes \Lambda)\bigr)\\
 &= \mu_{T(X)} T(\gamma_X).
\end{align*}
Hence $T(\mu_X)\gamma_{T(X)}=\mu_{T(X)} T(\gamma_X)$. Moreover, using Lemma~\ref{propgamr}(a) and since
$T_0\Lambda=\id_\un$, we have $\mu_X \gamma_X= \mu_X\Gamma_X (\id_X \otimes \Lambda)=(\eta_X \otimes T_0\Lambda)=\eta_X$. We
conclude that $T$ is separable by Proposition~\ref{equiv-sep-def}.
\end{proof}

\section{Sovereign and involutory Hopf monads}\label{sect-sovinv}

In this section, we introduce and study sovereign and involutory Hopf monads.

\subsection{Sovereign categories}\label{sect-sovecat}
Let $\cc$ be a left autonomous category. Recall that the choice of a left dual $(\ldual{X},\lev_X, \lcoev_X)$ for each
object $X$ of $\cc$ defines a left dual functor $\ldual{?}\co \cc^\opp \to \cc$. This is a strong monoidal functor and it is
unique up to unique monoidal isomorphism. In particular, this leads to the \emph{double left dual functor} $\lldual{?}\co
\cc \to \cc$, defined by $X \mapsto \ldual{(\ldual{X})}$ and $f \mapsto \ldual{(\ldual{f})}$, which is a strong monoidal
functor.

A \emph{sovereign structure} on a left autonomous category $\cc$ consists in the choice of a left dual for each object of
$\cc$ (and hence a strong monoidal functor $\lldual{?}\co\cc \to \cc$) together with a monoidal morphism $\phi \co 1_\cc \to
\lldual{?}$.  Note that, by Lemma \ref{isomono}, a sovereign structure $\phi \co 1_\cc \to \lldual{?}$ is automatically an
isomorphism.  Two sovereign structures are \emph{equivalent} if the corresponding monoidal isomorphisms coincide via the
canonical identification of the double dual functors.

A \emph{sovereign category} is a left autonomous category endowed with an equivalence class of sovereign structures.

Let $\cc$ be a sovereign category, with chosen left duals $(\ldual{X}, \lev_X, \lcoev_X)$ and sovereign isomorphisms
$\phi_X\co X \iso \lldual{X}$. For each object $X$ of $\cc$, set:
\begin{align*}
&\rev_X=\lev_{\ldual{X}}(\phi_X \otimes \id_{\ldual{X}}) \co X \otimes \ldual{X} \to \un,\\
&\rcoev_X=(\id_{\ldual{X}} \otimes \phi^{-1}_X)\lcoev_{\ldual{X}}\co \un \to \ldual{X} \otimes X.
\end{align*}
Then $(\ldual{X}, \rev_X, \rcoev_X)$ is a right dual of $X$. Therefore $\cc$ is autonomous. Moreover the right dual functor
$\rdual{?} \co \cc^\opp \to \cc$ defined by this choice of right duals coincides with $\ldual{?}$ as a strong monoidal
functor. However we will not necessarily make this choice of duals.

\begin{rem}\label{sovstructiso}
Let $\cc$ be a sovereign category, with sovereign structure $\phi\co 1_\cc \to \lldual{?}$. Since $\cc$ is autonomous, we
have $\phi^{-1}= \lexcla{\phi}= \rexcla{\phi}$ by Lemma \ref{isomono}. Explicitly, we have $\phi_X^{-1}=
\ldual{(\phi_{\rdual{X}})}= \rdual{(\phi_{\ldual{X}})}$ and $\lldual{(\phi_{\rrdual{X}})}=
\rrdual{(\phi_{\lldual{X}})}=\phi_X$ for all object $X$ of $\cc$ (up to the canonical isomorphisms of Remark~\ref{qieq}).
\end{rem}

\subsection{Sovereign functors}
Let $\cc$, $\dd$ be sovereign categories, with sovereign structure $\phi$ and $\phi'$ respectively. Let $F \co \cc \to \dd$
be a strong monoidal functor. Recall (see Section~\ref{strongfunctdual}) that $F$ defines a functorial isomorphism
$F^l_1(X)\co F(\ldual{X}) \to \ldual{F(X)}$.  Hence a functorial isomorphism
$F^{ll}_1(X)=\ldual{(F^l_1(X)^{-1})}F^l_1(\ldual{X})\co F(\lldual{X}) \to \lldual{F(X)}$. We will say that $F$ is
\emph{sovereign} if
\begin{equation}\label{deffunctsove}
F(\phi_X)F^{ll}_1(X)=\phi'_{F(X)}
\end{equation}
for all object $X$ of $\cc$.

\subsection{Square of the antipode}\label{sect-S2}
Let $\cc$ be a sovereign category, with  sovereign structure $\phi\co 1_\cc \to \lldual{?}$, and $T$ be a Hopf monad
on~$\cc$. Define $S^2 \in \Nat(T,T)$ by
\begin{equation}\label{defSsquare}
S^{2}_X=\phi^{-1}_{T(X)} s^l_{\ldual{T(X)}}T(\ldual{(s^l_X)}) T(\phi_X)
\end{equation}
for any object $X$ of $\cc$. We call $S^2$ the \emph{square of the antipode} of $T$.

Note that $S^2$ depends, in general, on the sovereign structure of $\cc$.

\begin{exa}
Let $A$ be a Hopf algebra in a braided autonomous category $\cc$, with braiding by $\tau$ and sovereign structure $\phi\co
1_\cc \to \lldual{?}$. Then the square of the antipode of the left Hopf monad $A \otimes ?$ on~$\cc$ (see
Example~\ref{ExHopf}) is given by $S^2_X=\phi_A^{-1} \uu_A(S_A)^2 \otimes \id_X$ for any object $X$ of~$\cc$, where $S_A$ is
the antipode of $A$ and $\uu_A = (\lev_A \tau_{A,\ldual{A}} \otimes \id_{\lldual{A}})(\id_A \otimes \lcoev_{\ldual{A}})\co A
\to \lldual{A}$ is  the Drinfeld isomorphism (see Section~\ref{sect-braidtwistcat}). Note that if $\cc$ is ribbon with twist
$\theta$ (see Section~\ref{sect-braidtwistcat}), then $\phi_A=\uu_A\theta_A$ and so $S^2_X=\theta_A^{-1}(S_A)^2 \otimes
\id_X$.   In particular, if $H$ is a finite-dimensional Hopf algebra over a field $\kk$, then the square of the antipode of
the left Hopf monad $H \otimes_\kk ?$ on~$\vect(\kk)$  is given by $S^2_X=(S_H)^2 \otimes \id_X$ for any finite-dimensional
$\kk$-vector space $X$.
\end{exa}

\begin{prop}\label{propS2auto}
The functorial morphism $S^2\co T \to T$ is an automorphism of the Hopf monad $T$ (see Section~\ref{sect-morphhopfmon}).
Moreover the inverse of $S^2$, denoted $S^{-2}$, is given by:
\begin{equation*}
S^{-2}_X=\phi_{\rrdual{T(X)}} s^r_{\rdual{T(X)}}T(\rdual{(s^r_X)}) T(\phi_{\rrdual{X}}^{-1})
\end{equation*}
for all object $X$ of $\cc$ (up to the canonical isomorphisms of Remark~\ref{qieq}).
\end{prop}

\begin{rem}
 Recall that, in Section~\ref{sect-convoantip}, we defined an anti-automorphism $S$ of the monoid
$(\Hom(1_\cc,T),*,\eta)$. Nevertheless, the notations are not in conflict since $S^{2}f=(S)^2(f)$ and
$S^{-2}(f)=(S^{-1})^2(f)$ for all $f\in\Nat(1_\cc,T)$.  In particular $S^2f$ does not depend on the sovereign structure of
$\cc$ (unlike $S^2$).
\end{rem}

\begin{proof}[Proof of Proposition~\ref{propS2auto}]
Let $(M,r)$ be a $T$-module. By Theorem~\ref{dual-ant}(a), we have $\lldual{(M,r)}=(\lldual{M},\lldual{r}\Sigma_M)$, where
$\Sigma_M=s^l_{\ldual{T(M)}}T(\ldual{(s^l_M)})$. Also, if $\varphi\co N \to M$ is an isomorphism in $\cc$, then
$(M,r)^\varphi=(N,\varphi^{-1}rT(\varphi))$ is a $T$-module. Define a functor $F\co T\ti\cc \to T\ti\cc$ by
\begin{equation*}
(M,r) \mapsto \bigl(\lldual{(M,r)}\bigr)^{\phi_M}=\bigl(M, \phi_M^{-1} \lldual{r}\Sigma_M T(\phi_M) \bigr )= (M, rS^2_M),
\quad f \mapsto f.
\end{equation*}
By the preliminary remarks, $F$ is well-defined. Since $\phi\co 1_\cc \to \lldual{?}$ is monoidal, the functor $F$ is
monoidal strict. Also $U_TF=U_T$. Therefore $S^2$ is a morphism of bimonads by Lemma~\ref{lemmorphbimon}, and so of Hopf
monads.

Let us show that $S^2$ is an automorphism. Remark that if $(N,s)$ is a $T$-module and $\varphi\co N \to M$ is an isomorphism
in $\cc$, then $(N,s)_\varphi=(M,\varphi s T(\varphi^{-1}))$ is a $T$\ti module. Also $((M,r)^\varphi)_\varphi=(M,r)$ and
$((N,s)_\varphi)^\varphi=(N,s)$ for any $T$-modules $(M,r)$, $(N,s)$ and any isomorphism $\varphi\co N \to M$ in $\cc$.
Therefore $F$ is an autofunctor of $T\ti\cc$ with inverse given by $F^{-1}(M,r)=\rrdual{((M,r)_{\phi_M})}$ (up to the
canonical isomorphisms of Remark~\ref{qieq}). Now, by Theorem~\ref{dual-ant}(b),
$\rrdual{(M,r)}=(\rrdual{M},\rrdual{r}\Sigma'_M)$, where $\Sigma'_M=s^r_{\rdual{T(M)}}T(\rdual{(s^r_M)})$. Therefore:
\begin{align*}
F^{-1}(M,r) & =\rrdual{(\lldual{M},\phi_M r T(\phi^{-1}_M))}\\
&=\bigl(M, \rrdual{\phi}_{M} \rrdual{r} \rrdual{T(\phi_M^{-1})} \Sigma'_{\lldual{M}} \bigr )\\
&=\bigl(M, r\rrdual{\phi}_{T(M)} \Sigma'_{M} T(\rrdual{(\phi_M^{-1})}) \bigr )\\
&=\bigl(M, r  \phi_{\rrdual{T(M)}} \Sigma'_{M} T(\phi_{\rrdual{M}}^{-1})\bigr ) \quad \text{by Remark~\ref{sovstructiso}}\\
&=(M,rS^{-2}_M),
\end{align*}
where $S^{-2}_M=\phi_{\rrdual{T(M)}} s^r_{\rdual{T(M)}}T(\rdual{(s^r_M)}) T(\phi_{\rrdual{M}}^{-1})$. Again by
Lemma~\ref{lemmorphbimon}, we get that $S^{-2}\co T \to T$ is a morphism of Hopf monads and is an inverse of $S^2$.
\end{proof}

\begin{lem}\label{lemS2ad}
Let $\cc$ be a sovereign category, with sovereign structure $\phi \co 1_\cc \to \lldual{?}$, and $T$ be a Hopf monad on
$\cc$. Let $a\in\Nat(1_\cc, T)$ and $a^\sharp\in\Nat(U_T,U_T)$ as in Lemma~\ref{key-lemma}. The following conditions are
equivalent:
\begin{enumerate}
\renewcommand{\labelenumi}{{\rm (\roman{enumi})}}
\item $L_a=R_a S^2$, where $L_a$ and $R_a$ are defined as in~\eqref{defLR};
\item $\phi a^\sharp \in \Nat(U_T,\lldual{?} U_T)=\Nat(U_T, U_T\lldual{?}_{T\ti\cc})$ lifts to $\Nat(1_{T\ti\cc},\lldual{?}_{T\ti\cc})$.
\end{enumerate}
\end{lem}
\begin{proof}
Let $(M,r)$ be a $T$-module. Recall that $(\phi a^\sharp)_{(M,r)}=\phi_M r a_M$. Also, by Theorem~\ref{dual-ant}(a), we
have:
\begin{equation*}
\lldual{(M,r)}=\bigl(\lldual{M},\lldual{r}s^l_{\ldual{T(M)}}T(\ldual{(s^l_M)}) \bigr )=\bigl(\lldual{M}, \phi_M r S^2_M
T(\phi_M^{-1})\bigr ).
\end{equation*}
Therefore $\phi a^\sharp$ lifts to a functorial morphism $1_{T\ti\cc} \to \lldual{?}_{T\ti\cc}$ if and only if, for any
$T$-module $(M,r)$, we have $\phi_M r a_M r = \phi_M r S^2_M T(\phi_M^{-1}) T \bigl( \phi_M r a_M \bigr)$ or, equivalently,
$r\mu_M a_{T(M)} = r \mu_M T(a_M)S^2_M$ since $\phi$ is an isomorphism and $rT(r)=r\mu_M$. By Lemma~\ref{key-lemma},  this
last condition is equivalent to $\mu_X a_{T(X)} = \mu_X T(a_X)S^2_X$ for all object $X$ of $\cc$, that is, $L_a=R_aS^2$.
\end{proof}

\subsection{Sovereign Hopf monads}\label{sect-sovhopfmon}
Let $\cc$ be a sovereign category. A Hopf monad $T$ on  $\cc$ is \emph{sovereign} if it is endowed with a grouplike
element~$G$ (which is $*$-invertible by Lemma~\ref{gl-lem2}), called the \emph{sovereign element} of $T$, satisfying:
\begin{equation}
S^2=\ad_G. \label{sov3}
\end{equation}
Here $S^2$ is the square of the antipode of $T$ (see Section~\ref{sect-S2}) and $\ad$ is the adjoint action of $T$ (see
Section~\ref{sect-adj}).

\begin{prop}\label{equivstructsove}
Let $\cc$ be a sovereign category and $T$ be a Hopf monad on $\cc$. Then sovereign elements of $T$ are in bijection with
sovereign structures on~$T\ti\cc$.
\end{prop}

\begin{proof}
Denote $\phi\co 1_\cc \to \lldual{?}$ the sovereign structure of $\cc$. Suppose that $G$ is a sovereign element of $T$.
Since $S^2=\ad_G$ and so $L_G=R_GS^2$, the functorial morphism $\phi G^\sharp$ lifts to a functorial morphism $\Phi\co
1_{T\ti\cc} \to \lldual{?}_{T\ti\cc}$ by Lemma~\ref{lemS2ad}. Since $\phi$ and $G^\sharp$ are monoidal (see
Lemma~\ref{gl-lem1}), so is the lift $\Phi$ of $\phi G^\sharp$, which hence defines a sovereign structure on $T\ti\cc$.

Conversely, let $\Phi\co 1_{T\ti\cc} \to \lldual{?}_{T\ti\cc}$ be a sovereign structure on $T\ti\cc$. Since the functorial
morphism $\phi^{-1}U_T(\Phi)$ is monoidal, there exists a (unique) grouplike element $G$ of $T$ such that
$\phi^{-1}U_T(\Phi)=G^\sharp$ (by Lemma~\ref{gl-lem1}). Since $\phi G^\sharp=U_T(\Phi)$ lifts to the functorial morphism
$\Phi$, we have (by Lemma~\ref{lemS2ad}) that $L_G=R_GS^2$, that is $S^2=\ad_G$. Hence $G$ is a sovereign element of $T$.
\end{proof}

\subsection{Involutory Hopf monads}\label{sect-invol}
A Hopf monad $T$ on a sovereign category $\cc$ is \emph{involutory} if it satisfies $S^2=\id_T$, where $S^2$ denotes the
square of the antipode as defined in Section~\ref{sect-S2}. Note that this notion depends on the choice of a sovereign
structure on $\cc$.

\begin{prop}\label{equivetreinvol}
Let $\cc$ be a sovereign category and $T$ a Hopf monad on $\cc$. The following conditions are equivalent:
\begin{enumerate}
 \renewcommand{\labelenumi}{{\rm (\roman{enumi})}}
\item $T$ is involutory;
\item $\eta$ is a sovereign structure on $\cc$;
\item There exists a sovereign structure on $T\ti\cc$ such that the forgetful functor $U_T\co T\ti\cc \to \cc$ is sovereign;
\item We have $s^r_X=\phi_{\rdual{X}}^{-1}\, s^l_X \, T(\phi_{\rdual{T(X)}})$ for any object $X$ of $\cc$
(up to the canonical isomorphisms of Remark~\ref{qieq}), where $\phi$ is the sovereign structure of $\cc$.
\end{enumerate}
\end{prop}

\begin{proof}
Clearly (i) implies (ii)  since $\eta$ is grouplike and $\ad_\eta=\id_T$. Assume (ii) and equip $T\ti\cc$ with the sovereign
structure defined by $\eta$. Then  the forgetful functor $U_T$ is sovereign. Hence (ii) implies (iii).

Let us prove that (iii) implies (iv). By Theorem~\ref{dual-ant}, we have preferred choices of left and right duals of
$(M,r)$, namely $\ldual{(M,r)}=(\ldual{M}, s^l_M T(\ldual{r}))$ and $\rdual{(M,r)}=(\rdual{M}, s^r_M T(\rdual{r}))$. With
this choice of duals, $(U_T)^l_1(M,r)=\id_{\ldual{M}}$. Let $\Phi$ be a sovereign structure on $T\ti\cc$ such that $U_T$ is
sovereign. We have $U_T(\Phi)=\phi_{U_T}$ by~\eqref{deffunctsove}. Let $(M,r)$ be a $T$-module. We have
$\Phi_{\rdual{(M,r)}} s^r_M T(\rdual{r})= s^l_M T(\ldual{r}) T(\Phi_{\rdual{(M,r)}})$ since $\Phi_{\rdual{(M,r)}}$ is
$T$-linear, and so $\phi_{\rdual{M}} s^r_M T(\rdual{r})= s^l_M T(\phi_{\rdual{T(M)}})T(\rdual{r})$. Hence (iv) by
Lemma~\ref{key-lemma}.

Finally, let us prove that (iv) implies (i). For any object $X$ of $\cc$, we have:
\begin{align*}
S^2_X &= \phi^{-1}_{T(X)} s^l_{\ldual{T(X)}}T(\ldual{(s^l_X)}) T(\phi_X) \\
 & = \phi^{-1}_{T(X)}  \phi_{T(X)}s^r_{\ldual{T(X)}} T(\phi_{\rdual{T(\ldual{T(X)})}}^{-1})
      T(\ldual{(s^l_X)}) T(\phi_X) \\
 &=  s^r_{\ldual{T(X)}} T(\rdual{(s^l_X)}) = \id_{T(X)} \quad \text{by Proposition~\ref{anti-inv}}
\end{align*}
Hence  $T$ is involutory.
\end{proof}

\section{Quasitriangular and ribbon Hopf monads}\label{sect-quasirib}

In this section, we define \Rt matrices and twists for a Hopf monad. They encode the facts that the category of modules over
the Hopf monad is braided or ribbon. We first review some well-known properties of braided and ribbon categories.

\subsection{Braided categories, twists, and ribbon categories}\label{sect-braidtwistcat}
Recall that a \emph{braiding} on a monoidal category $\cc$ is a functorial isomorphism $\tau \in \Nat(\otimes,\otimes^\opp)$
such that:
\begin{align}
& \tau_{X,Y\otimes Z}=(\id_Y \otimes \tau_{X,Z})(\tau_{X,Y} \otimes \id_Z) ; \label{braid1a} \\
& \tau_{X \otimes Y,Z}=(\tau_{X,Z} \otimes \id_Y)(\id_X \otimes \tau_{Y,Z}); \label{braid1b}
\end{align}
for all objects $X,Y,Z$ of $\cc$. A \emph{braided} category is a monoidal category endowed with a braiding. If $\tau$ is a
braiding on $\cc$, then so is  its \emph{mirror} $\overline{\tau}$  defined by $\overline{\tau}_{X,Y}=\tau^{-1}_{Y,X}$.

If $\cc$ is a braided category, with braiding $\tau$, and if $(X,Y,e,h)$ is a duality in $\cc$, then $(Y,X, e \tau_{X,Y},
\tau^{-1}_{Y,X} h)$ is a duality too. In particular, a braided category which is left (resp.\@ right) autonomous is also
right (resp.\@ left) autonomous, and so is autonomous.

Let $\cc$ be a braided autonomous category. Let $\uu\co 1_\cc \to \lldual{?}$ be the functorial morphism defined, for any
object $X$ of $\cc$, by:
\begin{equation}\label{defdriniso}
\uu_X = (\lev_X \tau_{X,\ldual{X}} \otimes \id_{\lldual{X}})(\id_X \otimes \lcoev_{\ldual{X}}).
\end{equation}
\begin{lem}\label{lemdriniso}
The morphism $\uu$ enjoys the following properties:
\begin{enumerate}
  \renewcommand{\labelenumi}{{\rm (\alph{enumi})}}
  \item $\uu_{X \otimes Y}= (\uu_X \otimes \uu_Y)\tau^{-1}_{X,Y}\tau^{-1}_{Y,X}$ for all objects $X,Y$ of $\cc$;
  \item $\uu_\un = (\lldual{?})_0 \co \un \iso \lldual{\un}$;
  \item $\uu$ is an isomorphism and, for any object $X$ of $\cc$,
  \begin{equation*}
      \uu^{-1}_X=(\lev_{\ldual{X}} \otimes \id_X)(\id_{\lldual{X}} \otimes \tau^{-1}_{\ldual{X},X}\lcoev_{X}).
   \end{equation*}
\end{enumerate}
\end{lem}
We will refer to $\uu$ as the \emph{Drinfeld isomorphism} of $\cc$.

\begin{rem}
The Drinfeld isomorphism $\uu$ is monoidal (and so is  a sovereign structure on $\cc$) if and only if the braiding $\tau$ is
a \emph{symmetry}, that is, $\overline{\tau}=\tau$.
\end{rem}

Recall that a \emph{twist} on a braided category $\cc$, with braiding $\tau$,  is a functorial isomorphism $\Theta \in
\Nat(1_\cc,1_\cc)$ satisfying:
\begin{equation}
\Theta_{X\otimes Y}=(\Theta_X \otimes \Theta_Y)\tau_{Y,X}\tau_{X,Y} \label{ribbon1a}
\end{equation}
for all objects $X,Y$ of $\cc$. If $\cc$ is autonomous, then a twist $\Theta$ on $\cc$ is said to be \emph{self-dual} if it
satisfies
\begin{equation}
\ldual{?}\, \Theta=\Theta \ldual{?} \quad \text{(or, equivalently,\; $\rdual{?}\Theta=\Theta\, \rdual{?}$).}
\label{ribbon1b}
\end{equation}
A \emph{ribbon category} is a braided autonomous category endowed with a self-dual twist.

The following proposition  establishes a correspondence (via the Drinfeld isomorphism) between the notions of sovereign
structure and twist in the context of an autonomous braided category.

\begin{prop}[Deligne] \label{proptwistsinbraid}
Let $\cc$ be an autonomous braided category, and denote $\uu$ its Drinfeld isomorphism. The assignment $\Theta \mapsto \uu
\circ\Theta$ defines a bijection between twists on $\cc$ and sovereign structures on $\cc$.
\end{prop}

\subsection{Quasitriangular bimonads}
Let $T$ be a monad on a monoidal category~$\cc$. Recall (see Section \ref{convol}) that a functorial morphism $R\in
\Nat(\otimes, T \otimes^\opp T) =\Nat(\otimes,\otimes^\opp \circ T^{\times 2})$ is $*$-invertible if there exists a
(necessarily unique) functorial morphism $R^{*-1} \in \Nat(\otimes^\opp, T \otimes T) =\Nat(\otimes^\opp, \otimes \circ
T^{\times 2}) $ such that $R^{*-1}
* R= \eta \otimes \eta$ and $R * R^{*-1}= \eta\otimes^\opp \eta$, where $*$ is  the convolution product as defined in
 \eqref{convoprodini}.

An \emph{\Rt matrix} for a bimonad $(T,\mu,\eta)$ on a monoidal category $\cc$ is a $*$-invertible functorial morphism $R\in
\Nat(\otimes, T \otimes^\opp T)$ such that:
\begin{align}
& (\mu_Y \otimes \mu_X)R_{T(X),T(Y)}T_2(X,Y)=(\mu_Y \otimes \mu_X)T_2(T(Y),T(X))T(R_{X,Y}); \label{rmat1}\\
\begin{split}
& (\id_{T(Z)}  \otimes T_2(X,Y))R_{X \otimes Y,Z} \\
& \phantom{XXXXXX}=(\mu_Z \otimes \id_{T(X) \otimes T(Y)}) (R_{X,T(Z)} \otimes \id_{T(Y)} )(\id_X \otimes R_{Y,Z}) ;
\end{split}\label{rmat2}\\
\begin{split}
& (T_2(Y,Z) \otimes \id_{T(X)})R_{X,Y \otimes Z} \\
& \phantom{XXXXXX}=( \id_{T(Y) \otimes T(Z)} \otimes \mu_X) ( \id_{T(Y)} \otimes R_{T(X),Z} ) (R_{X,Y} \otimes \id_Z);
\end{split}\label{rmat3}
\end{align}
for all objects $X,Y,Z$ of $\cc$.  A \emph{quasitriangular bimonad} is a bimonad equipped with an \Rt matrix.

\begin{exa}\label{exaquasik}
Let $H$ be a bialgebra over a field $\kk$. Let $r=\sum_i a_i \otimes b_i \in H \otimes_\kk H$. For any \kt vector spaces $X$
and $Y$, set:
\begin{equation*}
R_{X,Y}(x \otimes y)=\sum_i  b_i \otimes y \otimes a_i \otimes x \in H \otimes_\kk Y \otimes_\kk H \otimes_\kk X.
\end{equation*}
Then $R$ is an \Rt matrix for the bimonad $H\otimes_\kk ?$ on $\Vect(\kk)$ if and only if $r$ is an \Rt matrix for $H$ (in
the usual sense).
\end{exa}

\begin{thm}\label{braidthm}
Let $T$ be a bimonad on a monoidal category~$\cc$. Any \Rt matrix $R$ for $T$ yields a braiding $\tau$ on $T\ti\cc$ as
follows:
\begin{equation*}
\tau_{(M,r),(N,s)}=(s \otimes t)R_{M,N}\co (M,r) \otimes (N,s) \to (N,s) \otimes (M,r)
\end{equation*}
for any $T$-modules $(M,r)$ and $(N,s)$. This assignment gives a bijection between \Rt matrices for~$T$ and braidings on
$T\ti\cc$.
\end{thm}

\begin{proof}
Let $R\in\Nat(\otimes,T \otimes^\opp T)$ and set $\tau=R^\sharp$, where the canonical bijection $ ?^\sharp\co\Nat(\otimes,T
\otimes^\opp T) \rightarrow\Nat(U_T \otimes U_T,U_T \otimes^\opp U_T) $  of Lemma~\ref{key-lemmaN}  is given by
$f^\sharp_{(M,r),(N,s)}=(s \otimes r)f_{M,N}$ for all $T$-modules $(M,r)$ and $(N,s)$. In this correspondence, $\tau$ is an
isomorphism if and only if $R$ is $*$-invertible, and $\tau$ is $T$\ti linear in each variable (and so lifts to an element
of $\Nat(\otimes_{T\ti\cc}, \otimes_{T\ti\cc}^\opp)$) if and only if $R$ satisfies \eqref{rmat1}. Moreover, $\tau$ satisfies
\eqref{braid1a} and \eqref{braid1b} if and only if $R$ satisfies \eqref{rmat2} and \eqref{rmat3}. Hence the bijection
between \Rt matrices and braidings.
\end{proof}

\begin{cor}\label{corR1}
If $R$ is an \Rt matrix for a bimonad $T$, then $R^{*-1}_{21}=R^{*-1} \sigma_{\cc,\cc}$ is also an \Rt matrix for $T$.
Moreover, if $\tau$ is the braiding on $T\ti\cc$ associated with $R$, then its mirror $\overline{\tau}$ is the braiding on
$T\ti\cc$ associated with $R^{*-1}_{21}$.
\end{cor}
\begin{proof}
Let $R$, $R'$ be two \Rt matrices for $T$ and let $\tau$, $\tau'$ be their associated braidings on $T\ti\cc$ (see
Theorem~\ref{braidthm}). Given two $T$-modules $(M,r)$ and $(N,s)$, we have
\begin{align*}
\tau_{(N,s),(M,r)} \tau'_{(M,r),(N,s)}&=(r \otimes s) R_{N,M} (s \otimes r) R'_{M,N}\\
&=(rT(r) \otimes sT(s)) R_{T(N),T(M)} R'_{M,N}\\
&=(r\mu_r \otimes s\mu_s) (R_{2,1})_{T(M),T(N)} R'_{M,N} \quad \text{by \eqref{Tmoddef}}\\
&=(r \otimes s) (R_{2,1}*R')_{M,N}=(R_{2,1}*R')^\sharp_{(N,s),(M,r)}.
\end{align*}
As a result, by Lemma~\ref{key-lemmaN}, $\tau'=\overline{\tau}$ if and only if $R'=R^{*-1}_{2,1}$.
\end{proof}

\begin{cor}\label{corR2}
Let $T$ be a quasitriangular bimonad on a monoidal category~$\cc$. Then its \Rt matrix $R$ verifies $(\id \otimes
T_0)R_{\un,X}=\eta_X=(T_0 \otimes \id)R_{X,\un}$ as well as the following Yang-Baxter equation:
\begin{align*}
(\mu_Z & \otimes \mu_Y \otimes \mu_X)(R_{T(Y),T(Z)} \otimes \id_{T^2(X)})(\id_{T(Y)} \otimes R_{T(X),Z})(R_{X,Y} \otimes
\id_Z)\\
 & = (\mu_Z \otimes \mu_Y \otimes \mu_X)(\id_{T^2(Z)} \otimes R_{T(X),T(Y)})(R_{X,T(Z)} \otimes \id_{T(Y)})(\id_X \otimes
 R_{Y,Z}).
\end{align*}
Moreover, if $\cc$ is left autonomous and $T$ has a left antipode $s^l$, then
\begin{align*}
& R^{*-1}_{X,Y} = \bigl (\id_{T(X) \otimes T(Y)}\otimes\lev_X (s^l_X \otimes \id_X)\bigr )\\
& \phantom{XXXXXX} (\id_{T(X)} \otimes
R_{\ldual{T(X)},Y} \otimes \id_X ) (\lcoev_{T(X)} \otimes\id_{Y \otimes X});\\
& \ldual{R_{X,Y}}=(s^l_Y \otimes s^l_X)R_{\ldual{T(X)},\ldual{T(Y)}}.
\end{align*}
Likewise, if $\cc$ is right autonomous and $T$ has a right antipode $s^r$, then
\begin{align*}
& R^{*-1}_{X,Y} =\bigl(\rev_Y (\id_Y \otimes s^r_Y )  \otimes\id_{T(X) \otimes T(Y)}\bigr )\\
& \phantom{XXXXXX} (\id_Y \otimes R_{X,\rdual{T(Y)}} \otimes \id_{T(Y)} ) (\id_{Y \otimes X}
\otimes\rcoev_{T(Y)});\\
&\rdual{R}_{X,Y}=(s^r_Y \otimes s^r_X)R_{\rdual{T(X)},\rdual{T(Y)}}.
\end{align*}
\end{cor}
\begin{proof}
The corollary results, by standard application of  Lemma \ref{key-lemmaN}, from the facts that a braiding $\tau$
satisfies $\tau_{X,\un}=\id_X=\tau_{\un, X}$, the Yang Baxter equation,  $\tau^{-1}_{X,Y}=(\id \otimes \lev_X)(\id \otimes
\tau_{\ldual{X}, Y} \otimes \id)(\lcoev_X \otimes \id) $ and $\ldual{(\tau_{X,Y})}=\tau_{\ldual X, \ldual Y}$ when $\cc$ is
left autonomous, and $\tau^{-1}_{X,Y}=(\rev_Y \otimes \id)(\id \otimes \tau_{X,\rdual{Y}} \otimes \id)(\id \otimes
\rcoev_Y)$ and $\rdual{(\tau_{X,Y})}=\tau_{\rdual X, \rdual Y}$ when $\cc$ is right autonomous.
\end{proof}

\begin{cor}\label{corR3}
Let $T$ be a quasitriangular bimonad on an autonomous category $\cc$. If $T$ has a left (resp.\@ right) antipode, then $T$
has also a right (resp.\@ left) antipode, and so is a quasitriangular Hopf monad.
\end{cor}
\begin{proof}
Suppose that $T$ has a left antipode. Since $T\ti\cc$ is left autonomous (by Theorem~\ref{dual-ant}(a)) and braided (by
Theorem~\ref{braidthm}), $T\ti\cc$ is also right autonomous. Hence $T$ has a right antipode (by Theorem~\ref{dual-ant}(b)).
\end{proof}

\subsection{Drinfeld elements}
In this section, $T$ is a quasitriangular Hopf monad on a sovereign category~$\cc$ (see Section~\ref{sect-sovecat}). Let
$\phi\co 1_\cc \to \lldual{?}$ be the sovereign structure of $\cc$ and $R$ be the \Rt matrix of $T$.

The \emph{Drinfeld element} of $T$ is the functorial morphism $u\in \Nat(1_\cc,T)$ defined, for any object $X$ of $\cc$, by:
\begin{equation}\label{DrElteq}
u_X=\bigl (\lev_{T(X)}(s^l_{T(X)} \otimes \id_{T(X)})R_{X,\ldual{T^2(X)}} \otimes \mu_X \phi^{-1}_{T^2(X)}\bigr)(\id_X
\otimes \lcoev_{\ldual{T^2(X)}}).
\end{equation}

\begin{exa}\label{exadrielt}
Let $H$ be a finite-dimensional quasitriangular Hopf algebra over a field $\kk$. Recall that $H \otimes_\kk ?$ is a
quasitriangular Hopf monad on $\vect(\kk)$ (see Example~\ref{exaquasik}). Then the Drinfeld element $u$ of $H\otimes_\kk ?$
is given by $u_X(x)=d \otimes x$, where $d$ is the (usual) Drinfeld element of $H$. Recall that $d=\sum_i S(b_i)a_i$, where
$r=\sum_i a_i \otimes b_i \in H \otimes H$ is the \Rt matrix of $H$.
\end{exa}

\begin{lem}\label{lemDrinelt}
We have $U_T(\uu)=\phi \, u^\sharp$, where $U_T\co T\ti\cc \to \cc$ is the forgetful functor, $?^\sharp \co \Nat(1_\cc,T)
\to \Nat(U_T,U_T)$ is the canonical bijection of Lemma~\ref{key-lemma}, and $\uu$ is the Drinfeld isomorphism of $T\ti\cc$
(see Section~\ref{sect-braidtwistcat}).
\end{lem}
\begin{proof}
Let $(M,r)$ be a $T$-module. By Theorems~\ref{dual-ant} and \ref{braidthm}, we have
\begin{align*}
  U_T(\uu_{(M,r)})
 & =U_T\bigl((\lev_{(M,r)} \tau_{(M,r),\ldual{(M,r)}} \otimes \id_{\lldual{(M,r)}})(\id_{(M,r)} \otimes \lcoev_{\ldual{(M,r)}})\bigr)\\
 & =(\lev_M (s^l_M T(\ldual{r}) \otimes r)R_{M,\ldual{M}} \otimes \id_{\lldual{M}})(\id_M \otimes \lcoev_{\ldual{M}})\\
 & =(\lev_{T(M)} (s^l_M T(\ldual{(rT(r))}) \otimes \id_{T(M)})R_{M,\ldual{M}} \otimes \id_{\lldual{M}})(\id_M \otimes \lcoev_{\ldual{M}})\\
& =(\lev_{T(M)} (s^l_M T(\ldual{(r\mu_M)}) \otimes \id_{T(M)})R_{M,\ldual{M}} \otimes \id_{\lldual{M}})(\id_M \otimes \lcoev_{\ldual{M}})\\
& =(\lev_{T(M)} (s^l_M \otimes \id_{T(M)})R_{M,\ldual{T^2(M)}} \otimes \lldual{(r\mu_M)})(\id_M \otimes
\lcoev_{\ldual{T^2(M)}}).
\end{align*}
Since $\phi^{-1}_M\lldual{(r\mu_M)}=r\mu_M\phi^{-1}_{T^2(M)}$, we get $U_T(\uu_{(M,r)})=\phi_Mru_M=\phi_M u^\sharp_{(M,r)}$.
\end{proof}

\begin{prop}\label{propDrielt1}
The Drinfeld element $u$ of $T$ satisfies:
\begin{enumerate}
 \renewcommand{\labelenumi}{{\rm (\alph{enumi})}}
 \item $T_2 u_\otimes = (u \otimes u)*R^{*-1}*R^{*-1}_{21}$, where $(T_2 u_\otimes)_{X,Y}=T_2(X,Y) u_{X \otimes Y}$;
 \item $T_0 u_\un=\id_\un$;
 \item $u$ is $*$-invertible and, for any object $X$ of $\cc$,
     \begin{equation*}
     u^{*-1}_X =(\lev_{\ldual{X}} \otimes \id_{T(X)})\bigl(\phi_X \otimes (s^l_X \otimes \mu_X)R_{T(X),T(X)^*} \lcoev_{T(X)} \bigr );
     \end{equation*}
 \item $S^2=\ad_u$, where $S^2$ and $\ad_u$ are as in \eqref{defad} and \eqref{defSsquare} respectively.
\end{enumerate}
\end{prop}

\begin{proof}
Denote $\tau$ the braiding of $T\ti\cc$ induced by $R$. Let $\uu$ be the Drinfeld isomorphism of $T\ti\cc$ (see
Section~\ref{sect-braidtwistcat}) and $?^\sharp\co \Nat(1_\cc,T) \to \Nat(U_T,U_T)$  be the canonical bijection of
Lemma~\ref{key-lemma}.  Recall that $U_T(\uu)=\phi \, u^\sharp$ by Lemma~\ref{lemDrinelt}.

Let us prove Part (a). Let $(M,r)$ and $(N,s)$ be $T$-modules. By Lemma~\ref{lemdriniso}(a),
\begin{equation*}
\uu_{(M,r) \otimes (N,s)}= (\uu_{(M,r)} \otimes \uu_{(N,s)})\tau^{-1}_{(M,r),(N,s)}\tau^{-1}_{(N,s),(M,r)}.
\end{equation*}
Evaluating with $U_T$, and since $U_T(\uu)_{(M,r)}=\phi_M u^\sharp_{(M,r)}= \phi_M ru_M$, we get:
\begin{equation*}
  \phi_{M \otimes N} (r \otimes s)T_2(M,N)u_{M \otimes N}
  = (\phi_M ru_M \otimes \phi_N su_N) (r \otimes s)R^{*-1}_{M,N} (s \otimes r)R^{*-1}_{N,M}.
\end{equation*}
Therefore, since $\phi$ is a monoidal isomorphism,
\begin{align*}
  (r \otimes s)&T_2(M,N)u_{M \otimes N}\\
  &= (ru_M \otimes su_N) (r T(r) \otimes sT(s))R^{*-1}_{T(M),T(N)} R^{*-1}_{N,M}\\
  &= (rT(r)u_{T(M)} \otimes sT(s)u_{T(N)}) (\mu_M \otimes \mu_N)R^{*-1}_{T(M),T(N)} R^{*-1}_{N,M}
      \quad \text{by \eqref{Tmoddef}}\\
  &= (r \mu_M u_{T(M)} \otimes s \mu_N u_{T(N)}) (\mu_M \otimes\mu_N)R^{*-1}_{T(M),T(N)} R^{*-1}_{N,M}
      \quad \text{by \eqref{Tmoddef}}\\
  &= (r \otimes s)\bigl( (u \otimes u)*R^{*-1}*R^{*-1}_{21}\bigr)_{M,N}.
\end{align*}
Hence Part (a) by  Lemma~\ref{key-lemmaN}.

Let us prove Part (b). We have:
\begin{align*}
T_0u_\un &= T_0 u^\sharp_{(T(\un),\mu_\un)}\eta_\un \quad \text{by Lemma~\ref{key-lemma}}\\
& = T_0\phi_{T(\un)}^{-1}U_T(\uu_{(T(\un),\mu_\un)})\eta_\un\\
& = \phi_{\un}^{-1} \lldual{T_0} U_T(\uu_{(T(\un),\mu_\un)})\eta_\un\\
& = \phi_{\un}^{-1} U_T(\uu_{(\un,T_0)})T_0\eta_\un
    \quad \text{since $T_0$ is $T$-linear by \eqref{bimonad2}}\\
& = \phi_{\un}^{-1} U_T(\uu_{(\un,T_0)})
    \quad \text{by \eqref{bimonad4}}\\
&= \id_\un  \quad \text{by Lemma~\ref{lemdriniso}(b).}
\end{align*}

Let us prove Part (c). Set $u'=\bigl(U_T(\uu^{-1})\phi \bigr)^\flat \in \Nat(1_\cc,T)$, that is,
\begin{equation*}
u'_X= (\lev_{\ldual{X}} \otimes \id_{T(X)})\bigl(\phi_X \otimes (s^l_X \otimes \mu_X)R_{T(X),T(X)^*} \lcoev_{T(X)} \bigr ).
\end{equation*}
Then $u'^\sharp u^\sharp=U_T(\uu^{-1})\phi \phi^{-1}U_T(\uu)=\id_{U_T}$ and $u^\sharp u'^\sharp=
\phi^{-1}U_T(\uu)U_T(\uu^{-1})\phi=\id_{U_T}$. Therefore $u'*u=\eta=u*u'$ by Lemma~\ref{key-lemma}, that is, $u$ is
$*$-invertible with inverse $u'$.

Finally, let us prove Part (d). The functorial morphism $\phi  u^\sharp\in\Nat(U_T,\lldual{?}U_T)$ lifts to the functorial
morphism $\uu\in\Nat(1_{T\ti\cc},\lldual{?}_{T\ti\cc})$ by Lemma~\ref{lemDrinelt}. Therefore $L_u=R_uS^2$ by
Lemma~\ref{lemS2ad}, and so $S^2=\ad_u$ since $u$ is $*$-invertible.
\end{proof}

\subsection{Ribbon Hopf monads}
Let $T$ be a monad on a monoidal category~$\cc$. Recall (see Section \ref{convol}) that $\theta\in \Nat(1_\cc T)$ is
$*$-invertible if there exists a (necessarily unique) functorial morphism $\theta^{*-1} \in \Nat(1_\cc, T)$ such that
$\theta^{*-1} * \theta= \eta =\theta * \theta^{*-1}$, where $*$  is the convolution product as defined in
\eqref{prodconvo1T}. Recall also (see Section~\ref{sect-centralelts}) that $\theta\in \Nat(1_\cc T)$ is central if
$\mu_X\theta_{T(X)}=\mu_X T(\theta_X)$ for all object $X$ of $\cc$.

A \emph{twist} for a quasitriangular bimonad $T$ on a monoidal category $\cc$ is a  central and  $*$-invertible functorial
morphism $\theta \co 1_\cc \to T$ such that:
\begin{align}
T_2 \theta_\otimes =(\theta \otimes \theta) * R_{21} * R ,\label{ribT2}
\end{align}
where $R$ is the \Rt matrix of $T$ and $R_{21}=R \sigma_{\cc,\cc}$. Explicitly, \eqref{ribT2} means
\begin{equation*}
T_2(X,Y) \theta_{X \otimes Y}=(\mu_X \theta_{T(X)} \mu_X \otimes \mu_Y \theta_{T(Y)} \mu_Y)
  R_{T(Y),T(X)} R_{X,Y}
\end{equation*}
for all objects $X,Y$ of $\cc$.

A twist of a quasitriangular Hopf monad on an autonomous category is said to be \emph{self-dual} if it satisfies:
\begin{equation}\label{ribT3}
 \quad S(\theta)=\theta,
\end{equation}
where  $S\co \Nat(1_\cc T) \to \Nat(1_\cc T)$ is the map defined in  \eqref{defSL}. Explicitly, \eqref{ribT3} means
that $\ldual{\theta_X}=s^l_X\theta_{\ldual{T(X)}}$ (or $\rdual{\theta}_X=s^r_X\theta_{\rdual{T(X)}}$) for all object $X$ of
$\cc$.

A \emph{ribbon Hopf monad} is a quasitriangular Hopf monad on an autonomous category endowed with a self-dual twist.

\begin{exa}\label{exaribk}
Let $H$ be a finite-dimensional quasitriangular Hopf algebra over a field $\kk$. Then $H \otimes_\kk ?$ is a quasitriangular
monad on $\vect(\kk)$, see Example~\ref{exaquasik}. Let $v\in H$ and set  $\theta_X(x)= v \otimes x$  for any
finite-dimensional \kt vector space $X$ and $x \in X$. Then $\theta$ is self-dual twist for $H \otimes_\kk ?$ if and only if
$v$ is a ribbon element for~$H$.
\end{exa}

\begin{thm}\label{ribthm}
Let $T$ be a quasitriangular Hopf monad on an autonomous category~$\cc$. Any twist $\theta$ for $T$ yields a twist $\Theta$
on $T\ti\cc$ as follows:
\begin{equation*}
\Theta_{(M,r)}=r \theta_M\co (M,r) \to (M,r)
\end{equation*}
for any $T$-module $(M,r)$. This assignment gives a bijection between twists for $T$ and twists on $T\ti\cc$. Moreover, in
this correspondence, $\theta$ is self-dual (and so $T$ is ribbon) if and only if $\Theta$ is self-dual (and so $T\ti\cc$ is
ribbon).
\end{thm}
\begin{proof}
Let $\theta\in\Nat(1_\cc,T)$ and set $\Theta=\theta^\sharp$, where $?^\sharp$ is the canonical bijection $?^\sharp\co
\Nat(1_\cc,T) \rightarrow\Nat(U_T,U_T)$ given by $f^\sharp_{(M,r)}=rf_M$ for all $T$-module $(M,r)$. In this correspondence,
$\Theta$ is an isomorphism if and only if $\theta$ is $*$-invertible, and $\Theta$ is $T$\ti linear (and so lifts to a
functorial morphism $1_{T\ti\cc} \to 1_{T\ti\cc}$) if and only if $\theta$ is central (by Lemma~\ref{lemdefcentral}).
Moreover $\Theta$ satisfies \eqref{ribbon1a} if and only if $\theta$ satisfies \eqref{ribT2}. Hence the bijection between
twists for $T$ and twists on $T\ti\cc$. Finally, $\Theta$ satisfies \eqref{ribbon1b} if and only if $\theta$ satisfies
\eqref{ribT3}  by definition of $S$.
\end{proof}

\subsection{Ribbon and sovereign Hopf monads}
By Theorem~\ref{ribthm}, given a ribbon Hopf monad $T$ on a autonomous category $\cc$, the category $T\ti\cc$ of $T$-modules
is ribbon and so sovereign. Nevertheless $\cc$ itself is not necessarily sovereign. In this section, we study the case when
$\cc$ is sovereign, which allows to encode the  sovereign structure on $T\ti\cc$ by a sovereign element of $T$.
Recall that a sovereign element $G$ of $T$ is a grouplike element of $T$ satisfying $S^2=\ad_G$, where $S^2$ is the square
of the antipode of $T$ (see Section~\ref{sect-S2}) and $\ad$ is the adjoint action of $T$ (see Section~\ref{sect-adj}).

\begin{thm}\label{thmGtu}
Let $T$ be a quasitriangular Hopf monad on an sovereign category~$\cc$. Let $u$ be the Drinfeld element of $T$. Then the map
$\theta \mapsto G=u*\theta$ defines a bijection between twists of $T$ and sovereign elements of $T$. In this correspondence,
a twist $\theta$ is self-dual (and so $T$ is ribbon) if and only if the sovereign element $G=u*\theta$ satisfies
$S(u)=G^{*-1}*u*G^{*-1}$.
\end{thm}
\begin{proof}
Let $\phi$ be the sovereign structure on $\cc$, $\uu$ be the Drinfeld isomorphism of $T\ti\cc$ (see
Section~\ref{sect-braidtwistcat}), $?^\sharp \co \Nat(1_\cc,T) \to \Nat(U_T,U_T)$ be the canonical bijection of
Lemma~\ref{key-lemma}, and $?^\flat$ be the inverse of $?^\sharp$. Recall that $U_T(\uu)=\phi u^\sharp$ by
Lemma~\ref{lemDrinelt}. By Proposition \ref{proptwistsinbraid}, the assignment $\Theta \mapsto \uu\Theta$ defines a
bijection between twists on $T\ti\cc$ and sovereign structures on $T\ti\cc$. By Theorem~\ref{ribthm}, twists $\Theta$ on
$T\ti\cc$ are in bijection with twists $\theta$ for $T$. By Proposition~\ref{equivstructsove}, sovereign structures on
$T\ti\cc$ are in bijection with sovereign elements $G$ of $T$. Hence a bijection between twists $\theta$ for $T$ and
sovereign elements $G$ of $T$, which is given by:
\begin{equation*}
\theta \mapsto G=\bigl(\phi^{-1} U_T(\uu \Theta)\bigr)^\flat=(u^\sharp \theta^\sharp)^\flat=u*\theta.
\end{equation*}
Via this correspondence, we have $S(\theta)=\theta$ if and only if $S(u^{*-1}*G)=u^{*-1}*G$ or, equivalently (see
Lemmas~\ref{lemsrsl} and \ref{gl-lem2}), $S(u)=G^{*-1}*u*G^{*-1}$.
\end{proof}

\begin{cor}\label{corrib1}
Let $T$ be a ribbon Hopf monad on a sovereign category $\cc$, with twist $\theta$ and Drinfeld element $u$. Then
$\theta^{*-2}=u*S(u)=S(u)*u$.
\end{cor}
\begin{proof}
Since $G=u*\theta$ is grouplike by Theorem~\ref{thmGtu}, we have $G^{*-1}=S(G)$ by Lemma~\ref{gl-lem2}. Now
$S(G)=S(u*\theta)=S(\theta)*S(u)=\theta *S(u)$ by Lemma~\ref{lemsrsl} and \eqref{ribT3}. Therefore
$\theta^{*-1}*u^{*-1}=\theta *S(u)$ and so $\theta^{*-2}=S(u)*u$. Likewise, since we also have $G=\theta*u$ (because
$\theta$ is central), we get $\theta^{*-2}=u*S(u)$.
\end{proof}

\begin{cor}\label{corrib2}
Let $T$ be a quasitriangular Hopf monad on a sovereign category~$\cc$. Let $u$ be the Drinfeld element of $T$. Suppose that
$T$ is involutory (see Section~\ref{sect-invol}). Then $u^{*-1}$ is a twist for $T$, which is self-dual if and only if
$S(u)=u$.
\end{cor}
\begin{proof}
Results directly from Theorem~\ref{thmGtu} since, when $T$ is involutory, the unit $\eta$ of $T$ is a sovereign element of
$T$ (by Proposition~\ref{equivetreinvol}).
\end{proof}

\section{Examples and applications}\label{sect-perspec}

Throughout the previous sections, we chose the Hopf monad associated with a Hopf algebra as a paradigm of the notion of Hopf
monad. In this section, we give other examples of Hopf monads, so as to illustrate the generality of the notion.

\subsection{Monoidal adjunctions and Hopf monads}\label{sect-gen-ex}
Let $\cc$ and $\dd$ be categories. It is a standard fact (see \cite{ML1}) that if $(F\co \cc \to \dd, U\co \dd \to \cc)$ is
a pair of adjoint functors, with adjunction morphisms $\eta \co 1_\cc \to UF$ and $\varepsilon \co FU \to 1_\dd$, then $T=
UF$ is a monad on~$\cc$, with product $\mu=U (\varepsilon_F) \co T^2 \to T$ and unit $\eta$. Also there exists a unique
functor $K \co \dd \to T\ti\cc$ such that $U_T K = U$ and
$KF=F_T$. The functor $K$ is given by  $A \mapsto \bigl(U(A), U(\varepsilon_A)\bigr)$.
This fact admits the following monoidal version:

\begin{thm}\label{thm-gen-ex}
Let $\cc,\dd$ be two monoidal categories and $U\co\dd \to \cc$ be a strong monoidal functor. Assume that the functor $U$ has
a left adjoint $F$. Then $F$ is a comonoidal functor and the monad $T=UF$ on $\cc$ has a canonical structure of a bimonad.
The canonical functor $K \co \dd \to T\ti \cc$ is strong monoidal and satisfies $U_T K=U$ as monoidal functors and $KF=F_T$
as comonoidal functors. Furthermore, if $\cc$ and $\dd$ are  left (resp.\@ right) autonomous, then $T$ is a left (resp.\@
right) Hopf monad.
\end{thm}

\begin{rem}
Any bimonad or Hopf monad $T$ is of the form of Theorem~\ref{thm-gen-ex}, since the forgetful functor
$U_T$ is strong monoidal, $F_T$ is left adjoint to $U_T$, and $T= U_T F_T$.
\end{rem}
\begin{rem}
In the last assertion of Theorem~\ref{thm-gen-ex}, it is sufficient to assume that~$\cc$ is left (resp.\@ right) autonomous
and that, for any object~$X$ of $\cc$, $F(X)$ has a left (resp.\@ right) dual in $\dd$. Indeed, this can be seen by applying
Theorem~\ref{thm-gen-ex} to the restriction of $U$ to the full subcategory of $\dd$ made of objects having a left (resp.\@
right) dual.
\end{rem}

\begin{proof}[Proof of Theorem~\ref{thm-gen-ex}]
Denote $\eta \co 1_\cc \to UF$ and $\varepsilon \co FU \to 1_\dd$ the adjunction morphisms. Define $F_2\co F\otimes \to F
\otimes F$ by setting, for any objects $X,Y$ of $\cc$,
\begin{equation*}
F_2(X,Y)=\varepsilon_{F(X) \otimes F(Y)} F\bigl(U_2(F(X),F(Y))\bigr) F(\eta_X \otimes \eta_Y)
\end{equation*}
and set $F_0=\varepsilon_{\un_\dd} F(U_0)\co F(\un_\cc) \to \un_\dd$. One verifies that $(F,F_2,F_0)$ is a comonoidal
functor. Since $U$ is strong monoidal, we may also view it as a strong comonoidal functor (with  comonoidal structure
defined by $U_2^{-1}$ and~$U_0^{-1}$). Therefore both $T=UF$ and $FU$ are comonoidal functors by Lemma~\ref{compo}(b). One
checks that $\eta\co 1_\cc \to T$, $\varepsilon \co FU \to 1_\dd$, and $\mu=U (\varepsilon_F) \co T^2 \to T$ are comonoidal
morphisms. As a result, the monad $T=UF$ is a bimonad.

For any object $A$ of $\dd$, we have $K(A)=\bigl (U(A), U(\varepsilon_A) \bigr)$. For all objects $A,B$ of~$\dd$, the
morphism $U_2(A,B)\co U(A) \otimes U(B) \to U(A \otimes B)$ lifts to a ($T$-linear) morphism $K_2(A,B) \co K(A) \otimes K(B)
\to K(A \otimes B)$. Likewise, $U_0\co \un_\cc \to U(\un_\dd)$ lifts to a ($T$-linear) morphism $K_0
\co\un_{T\ti\cc}=(\un_\cc,T_0) \to K(\un_\dd)$. Moreover, $(K,K_2,K_0)$ is a strong monoidal functor such that $U_T K= U$ as
monoidal functors, because $U_T$ is strict monoidal, faithful, and conservative. We also have $K F= F_T$ as comonoidal
functors, since $K F= \bigl(UF, U(\varepsilon_{F})\bigr)=(T,\mu)=F_T$, $(KF)_2=U_2^{-1}(F,F)\,U(F_2)=T_2=(F_T)_2$, and
$(KF)_0=U_0^{-1}U(F_0)=T_0=(F_T)_0$.

Finally, assuming $\cc$ and $\dd$ are both left autonomous, we construct a left antipode $s^l$ for the bimonad $T$ as
follows (the right case can be done similarly). Let $A$ be an object of $\dd$. Since $A$ has a left dual in $\dd$ and $K$ is
strong monoidal, the $T$\trait module $K(A)$ has a left dual in $T\ti\cc$. By Lemma~\ref{prefer}, we may choose this left
dual of the form $\bigl((\ldual{U(A)},\rho_A),\lev_{U(A)}, \lcoev_{U(A)}\bigr)$, with $\rho_A \co T(\ldual{U(A)}) \to
\ldual{U(A)}$ uniquely determined. One verifies that $\rho_A$ is functorial in $A$. Note that the $T$-linearity of
$\lev_{U(A)}$ and $\lcoev_{U(A)}$ translate respectively as:
\begin{align}
& T_0 T(\lev_{U(A)})= \lev_{U(A)} (\rho_{A} \otimes U(\varepsilon_A)) T_2(\ldual{U(A)},U(A)),\label{linev}\\
&\lcoev_{U(A)} T_0 =(U(\varepsilon_A) \otimes \rho_A)T_2(U(A), \ldual{U(A)})T(\lcoev_{U(A)}).\label{lincoev}
\end{align}
Now, for any object $X$ of $\cc$, set $s^l_X= \ldual{\eta_X}\rho_{F(X)} \co T(\ldual{T(X)}) \to \ldual{X}$. Clearly $s^l$ is
functorial. Fix an object $X$ of $\cc$. By~\eqref{mon-u} and the functoriality of $\rho$, we have:
\begin{align*}
\rho_{F(X)} &=\ldual{\eta_{T(X)}}\ldual{\mu_X}\rho_{F(X)}=\ldual{\eta_{T(X)}}\ldual{U(\varepsilon_{F(X)})}\rho_{F(X)}\\
& =\ldual{\eta_{T(X)}}\rho_{FT(X)} T(\ldual{U(\varepsilon_{F(X)})}) =s^l_{T(X)}T(\ldual{\mu_X}).
\end{align*}
Therefore we get:
\begin{align*}
T_0 T & (\lev_X) T(\ldual{\eta_X} \otimes \id_X) =T_0 T(\lev_{T(X)})T(\id_{\ldual{T(X)}} \otimes \eta_X)\\
&=\lev_{T(X)} (s^l_{T(X)}T(\ldual{\mu_X}) \otimes \mu_X) T_2(\ldual{T(X)},T(X))T(\id_{\ldual{T(X)}} \otimes \eta_X)
  \quad \text{by \eqref{linev}}\\
&=\lev_{T(X)} (s^l_{T(X)}T(\ldual{\mu_X}) \otimes \id_{T(X)}) T_2(\ldual{T(X)},X) \quad \text{by \eqref{mon-u}.}
\end{align*}
Likewise we have:
\begin{align*}
(\eta_X \otimes & \id_{\ldual{X}}) \lcoev_X T_0 =(\id_{T(X)} \otimes \ldual{\eta_X})\lcoev_{T(X)}T_0\\
&=(\mu_X \otimes \ldual{\eta_X} s^l_{T(X)}T(\ldual{\mu_X})) T_2(T(X),\ldual{T(X)}) T(\lcoev_{T(X)})
\quad \text{by \eqref{lincoev}} \\
&=(\mu_X \otimes s^l_X) T_2(T(X),\ldual{T(X)}) T(\lcoev_{T(X)}) \quad \text{by \eqref{mon-u}.}
\end{align*}
Hence $s^l$ satisfies \eqref{lant1} and \eqref{lant2}, that is, $s^l$ is a left antipode for $T$.
\end{proof}

\subsection{Tannaka reconstruction} Fiber functors are an interesting source of examples of Hopf monads.

Let $\kk$ be a field. Given a \kt algebra $B$, we denote $\biMod{B}{B}$ the category of $B$\trait bimodules, $\bimod{B}{B}$
the category of  finitely generated $B$-bimodules, $\biproj{B}{B}$ the category of finitely generated projective
$B$-bimodules, and $\lmod{B}$ the category of finitely generated  left $B$-modules.

A \emph{tensor category over $\kk$} is an autonomous category endowed with a structure of $\kk$-linear abelian category such
that $\otimes$ is bilinear and $\End(\un)=\kk$. Let $\cc$ be a tensor category and $B$ be a $\kk$-algebra. A \emph{$B$-fiber
functor for $\cc$} is a $\kk$-linear exact strong monoidal functor $\cc \to \biMod{B}{B}$. A $B$-fiber functor takes values
in $\biproj{B}{B}$ (because it preserves duality) and it is faithful if $B$ is non-trivial (that is, $B \neq 0$).

We say that a tensor category  $\cc$ is \emph{bounded} if there exists a $k$-linear equivalence of categories $\Xi \co\cc \to
\lmod{E}$  for some finite-dimensional $\kk$-algebra $E$.  By Proposition~2.14 of \cite{De} due to O. Gabber, this is
equivalent to the  more intrinsic following  conditions:
\begin{itemize}
\item in $\cc$, all objects have
finite length and  $\Hom$ spaces are finite-dimensional;
\item $\cc$  admits a generator, i.e., an object $X$ such that any object is a subquotient
of $X^{\oplus n}$ for some integer $n$.
\end{itemize}

\begin{thm}\label{prop-fiber}
Let $\cc$ be a bounded tensor category over a field $\kk$, $B$ be a  non-trivial  finite-dimensional $\kk$-algebra, and
$\omega$ be a $B$-fiber functor for~$\cc$. Then the functor~$\omega$, viewed as a functor $\cc \to \bimod{B}{B}$, admits a
left adjoint $F$. The endofunctor $T=\omega F$ is a bimonad on $\bimod{B}{B}$ and induces, by restriction, a Hopf monad
$T_0$ on $\biproj{B}{B}$. The categories $T_0\ti \biproj{B}{B}$ and $T\ti \bimod{B}{B}$ are isomorphic (as $\kk$-linear
monoidal categories). Furthermore, the canonical functor $\cc \to T \ti \bimod{B}{B}\cong T_0 \ti \biproj{B}{B}$ is a
$\kk$-linear strong monoidal equivalence.
\end{thm}

\begin{proof}
Let $R$ and $E$ be two finite-dimensional $\kk$-algebras. Let $G\co \lmod{E} \to \lmod{R}$ be a $\kk$-linear exact functor.
Being right exact, $G$ is of the form $\leftidx{_R}{M}{} \otimes_E ?$ for some $R$-$E$-bimodule $M$. Since $G$ is left
exact, $M$ is flat, and it is also of finite type, hence projective. Let $\ldual{M}$ be the $E$-$R$-bimodule
$\Hom_R(\leftidx{_R}{M}{_E},\leftidx{_R}{R}{_R})$. Then $F=\leftidx{_E}{\ldual{M}}{}\otimes_R ?$ is left adjoint to $G$. So
$T=GF$ is a monad on $\lmod{R}$. Moreover, the canonical functor $K \co \lmod{E} \to T \ti \lmod{R}$ is an equivalence if
$G$ is faithful (that is, if $M$ is faithfully flat as an $E$-module).

Via a $\kk$-linear equivalence $\Xi \co \cc \to \lmod{E}$, this applies to $\omega\co \cc \to \bimod{B}{B}$ (with $R=B
\otimes_\kk B^{\opp}$),  and shows that $\omega$ has a left adjoint $F$. Hence $T=\omega F$ is a bimonad on $\bimod{B}{B}$
by Theorem~\ref{thm-gen-ex}. The canonical functor $K\co \cc \to T\ti \bimod{B}{B}$ is a $\kk$-linear strong monoidal
equivalence. Now $\omega$, and so $T$, takes  values in $\biproj{B}{B}$. Denote $\omega_0\co \cc \to \biproj{B}{B}$ and $F_0
\co \biproj{B}{B}\to \cc$ the restrictions of $\omega$ and $F$. Then $F_0$ is left adjoint to $\omega_0$ and $T_0 = F_0
\omega_0$ is a Hopf monad on $\biproj{B}{B}$ (by Theorem~\ref{thm-gen-ex}). Lastly, consider a $T$-module $(N,r)$, where $N$
is an object of $\bimod{B}{B}$.  Since $K$ is an equivalence, $(N,r)$ is isomorphic to $K(Y)$, for some $Y$ in $\cc$, and so
$N \simeq \omega(Y)$. In particular $N$ is in $\biproj{B}{B}$. Therefore we have $T\ti\bimod{B}{B}=T_0\ti\biproj{B}{B}$,
hence the theorem.
\end{proof}

\begin{cor} Let $\cc$ be a semisimple tensor category over a field $\kk$.  Assume the set $\Lambda$ of isomorphy classes of
simple objects of $\cc$ is finite and  $\End(V)=\kk$ for each simple object $V$ of $\cc$.  Let $B$ be the $\kk$-algebra
$\kk^\Lambda$. Then there exist a $\kk$-linear Hopf monad $T$ on $\bimod{B}{B}$ and a $\kk$-linear strong monoidal
equivalence $\cc \to T\ti \bimod{B}{B}$.
\end{cor}
\begin{proof}
By Theorem 4.1 of \cite{Hay}, we have a canonical $B$-fiber functor $\cc \to \bimod{B}{B}$. Hence the corollary by
Theorem~\ref{prop-fiber}, noticing that $\bimod{B}{B}=\biproj{B}{B}$ because  $B$ is semisimple.
\end{proof}

\subsection{Doubles of Hopf monads}
Let $\cc$ be an autonomous category and $T$ be a Hopf monad on $\cc$. Assume that the coend
\begin{equation*}
D_T(X)=\int^{Y \in \cc} \mspace{-25mu} \ldual{T(Y)} \otimes X \otimes Y
\end{equation*}
exists for all object $X$ of $\cc$. Then $D_T$ is a Hopf monad on $\cc$ (see \cite{BV2} for details), which is called the
\emph{double} of~$T$.

Consider in particular the double $D$ of the trivial Hopf monad $1_\cc$. Note that it always exists when, for example, $\cc$
has a generator or $\cc$ is finitely semisimple. Then $D$ is a quasitriangular Hopf monad, whose category $D\ti\cc$ of
modules is canonically isomorphic (as braided categories) to the center $Z(\cc)$ of $\cc$, see \cite{BV2}. Note that the
forgetful functor $Z(\cc) \to \cc$ is then monadic.

\bibliographystyle{amsalpha}

\end{document}